\documentclass[12pt,oneside,english]{amsart}
\usepackage[T1]{fontenc}
\usepackage[latin9]{inputenc}
\usepackage{amsthm}
\usepackage{amssymb}

\makeatletter

\providecommand{\tabularnewline}{\\}

\numberwithin{equation}{section}
\numberwithin{figure}{section}
\theoremstyle{plain}
\newtheorem{thm}{\protect\theoremname}
  \theoremstyle{definition}
  \newtheorem{defn}[thm]{\protect\definitionname}
  \theoremstyle{plain}
  \newtheorem*{thm*}{\protect\theoremname}
  \theoremstyle{plain}
  \newtheorem{cor}[thm]{\protect\corollaryname}

\usepackage{tikz}
\usetikzlibrary{backgrounds}
\usetikzlibrary{decorations.fractals}
\usetikzlibrary{calc,intersections,through,backgrounds}

\makeatother

\usepackage{babel}
  \providecommand{\corollaryname}{Corollary}
  \providecommand{\definitionname}{Definition}
  \providecommand{\theoremname}{Theorem}
\providecommand{\theoremname}{Theorem}

\begin{document}

\title{The compactifications of moduli spaces of burniat surfaces with $2\leq K^{2}\leq5$}

\author{Xiaoyan Hu}
\begin{abstract}
We describe compactifications of moduli spaces of Burniat surfaces
with $2\leq K_{X}^{2}\leq5$ obtained by adding KSBA surfaces, i.e.
slc surfaces $X$ with ample canonical class $K_{X}$.

\tableofcontents{}
\end{abstract}
\maketitle

\section*{Introduction}

Burniat surfaces are special cases of surfaces of general type with
$p_{g}=q=0$, $2\le K_{X}^{2}\le6$. They were first introduced by
Burniat in \cite{Bu66}. Peters \cite{Pet77} reinterpreted Burniat's
construction using the modern language of branched abelian covers.
In \cite{LP01}, Lopes and Pardini proved that a minimal surface $S$
of general type with $p_{g}(S)=0,\mbox{ }K_{S}^{2}=6$, and bicanonical
map of degree 4 is a Burniat surface. Moreover, they showed that minimal
surfaces $S$ with $p_{g}=0,\mbox{ }K_{S}^{2}=6$ and bicanonical
map of degree 4 form a four-dimensional irreducible component of the
moduli space of surfaces of general type. 

In \cite{KSB88}, Kollár and Shepherd-Barron introduced stable surfaces
and proposed a way to compactify the moduli space of surfaces of general
type by adding stable surfaces (also called KSBA surfaces). They showed
that the appropriate singularities to permit for the surfaces at the
boundaries of moduli spaces are semi log canonical (slc) and classified
all the semi log canonical surface singularities. The boundedness
of slc surfaces with a fixed canonical class $K^{2}$ was settled
in \cite{Ale94}. In \cite{Ale96a,Ale96b}, Alexeev extended Kollár
and Shepherd-Barron's construction to stable pairs and stable maps.

In \cite{AP09}, Alexeev and Pardini constructed an explicit compactification
of the moduli space of Burniat surfaces with $K^{2}=6$ by adding
KSBA surfaces, i.e. slc surfaces $X$ with ample canonical class $K_{X}$,
on the boundary. They also gave a constructive algorithm for computing
all stable Burniat surfaces (not necessarily from degenerations of
smooth surfaces), which reduced them to computing certain tilings
by matroid polytopes. 

The aim of this paper is to extend the results and methods in \cite{AP09}
from the case $K^{2}=6$ to all the remaining cases $2\le K^{2}\leq5$.
The moduli space $M_{Bur}^{d}$ of Burniat surfaces with $K^{2}=d$
is a subset of dimension $d-2$ in the moduli space $\mathfrak{M}^{can}$
of canonical surfaces, where a point in $M_{Bur}^{d}$ corresponds
to the canonical model of a smooth Burniat surface. When $d=6,5$,
the moduli space $M_{Bur}^{d}$ is an irreducible component in $\mathfrak{M}^{can}$.
Bauer and Catanese \cite{BC10b} showed that $M_{Bur}^{4}$ is a union
of two irreducible subvarieties $M_{Bur,1}^{4}$ and $M_{Bur,2}^{4}$,
where a general point of $M_{Bur,1}^{4}$ corresponds to a smooth
Burniat surface, while a general surface in $M_{Bur,2}^{4}$ has an
$A_{1}$-singularity (nodal case). Moreover, $M_{Bur,1}^{4}$ is an
irreducible component in $\mathfrak{M}^{can}$, whereas $M_{Bur,2}^{4}$
is contained in an irreducible component of dimension 3 in $\mathfrak{M}^{can}$.
The moduli space $M_{Bur}^{3}$ is irreducible and is contained in
an irreducible component of dimension 4 in $\mathfrak{M}^{can}$.
$M_{Bur}^{2}$ is just one point so already compact. Thus we will
restrict ourselves to compactifying the moduli space $M_{Bur}^{d},3\leq d\leq5$. 

We reduce the problem of compactifying $M_{Bur}^{d}$ to the one of
compactifying the moduli space of certain stable pairs $(Y,\frac{1}{2}D)$.
A point in $M_{Bur}^{d}$ corresponds to a Burniat surface $X$ with
$K^{2}=d$, that is the canonical model of a $\mathbb{Z}_{2}^{2}$-cover
of $Y=Bl_{9-d}\mathbb{P}^{2}$ branched along $12+d$ irreducible
curves consisting of 9 strict preimages of lines and $3+d$ exceptional
divisors. The branch data is encoded in the Hurwitz divisor $D$ (see
Section \ref{subsecabcover}). An \emph{abelian cover} of a variety
$Y$ with group $G$ or a $G$-cover is a finite map $\pi:X\rightarrow Y$
together with a faithful action of a finite abelian group $G$ on
$X$ such that $\pi$ exhibits $Y$ as the quotient of $X$ by $G$.
In the case $Y$ is smooth and $X$ is normal, Pardini in \cite{Par91}
described the general structure of abelian covers $\pi:X\rightarrow Y$
using the building data which we will discuss in Section \ref{subsecabcover}.
The work was extended to the case of non-normal abelian covers in
\cite{AP12}. In Sections \ref{degree5},\ref{degree4},\ref{degree3},
we list all the interesting degenerate configurations of stable pairs
$(Y,\frac{1}{2}D)$ with $K^{2}=3,4,5$, up to symmetry, and find
their canonical models using the minimal model program for 3-folds.
Here, interesting degenerate configurations are the ones with reducible
canonical models.

The stable surfaces appearing on the boundary are quite nontrivial
and provide examples of many interesting features of the general case.
The construction of the compactified coarse moduli spaces $\overline{M}_{Bur}^{d}$
of Burniat surfaces is an application of \cite{Ale08}, which provides
a stable pair compactification $\overline{M}_{\beta}(r,n)$ for the
moduli space of weighted hyperplane arrangements $(\mathbb{P}^{r-1},\sum b_{i}B_{i})$
with arbitrary weight $\beta=(b_{1},...,b_{n})$, $0\leq b_{i}\leq1$
and $b_{i}\in\mathbb{Q}$. In this paper, we apply \cite{Ale08} in
the case of $\mathbb{P}^{2}$ and $n=9$ with $\beta=(\frac{1}{2},...,\frac{1}{2})$. 

Several new phenomena happen in the case $K^{2}\leq5$ as compared
to the case $K^{2}=6$ in \cite{AP09}. Most importantly, when running
the minimal model program, in addition to divisorial contractions
occurring in the case $K^{2}=6$, flips and flops also occur. It is
also surprising that some non log canonical degenerations in the case
$K^{2}=6$ correspond to log canonical degenerations in the cases
$K^{2}\leq5$ . 

We first study degenerations of stable pairs $(Y,\frac{1}{2}D)$ and
apply the minimal model program to find the stable limit. We summarize
our main results below. 
\begin{thm}
\label{maintheory}The main component of the compactified coarse moduli
space $\overline{M}_{Bur}^{d}$ of stable Burniat surfaces, or equivalently,
of stable pairs $\left(Y,\frac{1}{2}D\right)$, is of dimension $d-2$,
irreducible for $d\neq4$, and with two components for $d=4$. The
types of degenerations, up to symmetry, are listed as below.

(i) There are 6 types of degenerate configurations of stable pairs
with reduced log canonical models in the moduli space of stable pairs
$\left(Y,\frac{1}{2}D\right)$ for $K^{2}=5$ case up to the symmetry
group $\mathbb{Z}_{6}$ described in Section \ref{degree5}.

(ii) There are 5 types of degenerations with reducible lc models in
the moduli space of stable pairs $\left(Y,\frac{1}{2}D\right)$ for
$K^{2}=4$ nodal case and 3 types of degenerations for $K^{2}=4$
non-nodal case up to the symmetry group $\mathbb{Z}_{2}$ described
in Section \ref{degree4}.

(iii) There are only 2 types of degenerations with reducible lc models
in the moduli space of stable pairs $\left(Y,\frac{1}{2}D\right)$
for $K^{2}=3$ described in Section \ref{degree3}.\emph{ }\\

\end{thm}
According to the general theory of \cite{Ale08}, the unweighted stable
hyperplane arrangements are described by matroid tilings of the hypersimplex
$\triangle(r,n)$. Their weighted counterparts are described by partial
tilings of the hypersimplex $\triangle(r,n)$ that cover a $\beta$-cut
hypersimplex $\triangle_{\beta}(r,n)$.

The polytope $\triangle_{Bur}^{d},d\leq6$ is the polytope in $\mathbb{R}^{12}$
that corresponds to the stable pairs $(Y,\frac{1}{2}D)$ with $K^{2}=d$,
where $Y=Bl_{9-d}\mathbb{P}^{2}$. Inductively, we restrict the matroid
tilings of the polytope $\triangle_{Bur}^{d}$ for each $d=6,5,4$
to the polytope $\triangle_{Bur}^{d-1}$ and find all possible stable
pairs in the main component of the compactified moduli space of stable
pairs with $K^{2}=d-1$. The possible surfaces produced by this computation
exactly coincide with the degenerations listed in Sections \ref{degree5},\ref{degree4},\ref{degree3}.
This also shows that the stable pairs listed in Sections \ref{degree5},\ref{degree4},\ref{degree3}
are all the degenerations for the main components of the compactified
moduli space of stable pairs with $K^{2}=d\leq5$. All the tilings
of $\triangle_{Bur}^{d-1}$ corresponding to degenerations are restrictions
of some tilings of $\triangle_{Bur}^{d}$. However, not all restrictions
of the tilings of $\triangle_{Bur}^{d}$ to $\triangle_{Bur}^{d-1}$
correspond to degenerations of Burniat surfaces with $K^{2}=d-1$.
For example, the tiling \#1 of the polytope $\triangle_{Bur}^{6}$
listed in table 2 \cite{AP09} is a tiling of the polytope $\triangle_{Bur}^{5}$
as well. Tiling \#1 corresponds to the non log canonical degeneration
Case 1 with $K^{2}=6$ in \cite{AP09}, but it does not correspond
to any non log canonical degenerations in the case $K^{2}=5$. \\

\textbf{Acknowledgments. }I would like to express my deepest gratitude
to my advisor, Prof. Valery Alexeev, for his excellent guidance, help
and his patience. Thank R. Pignatelli for helpful comments.

\section{Burniat Surfaces}

\subsection{Preliminaries}

We say that a variety is \emph{d.c.} (\emph{double crossings}) if
every point is either smooth or analytically isomorphic to $xy=0$.
We say that a variety is \emph{g.d.c.} (has \emph{generically double
crossings}) if it is d.c. outside a closed subset of codimension $\geq2$.

Let $X$ be a projective variety. Let $B=\sum b_{i}B_{i}$ be a linear
combination of effective divisors, where $b_{i}$ is the weight of
$B_{i}$ which is allowed to be an arbitrary rational number with
$0<b_{i}\leq1$. The divisors $B_{i}$'s are possibly reducible and
possibly have irreducible components in common. We recall some basic
definitions.\\

\begin{defn}
Assume that $X$ is a normal variety. A pair $(X,B)$ is called \emph{log
canonical} \emph{(lc)} if 
\end{defn}
(1) $m(K_{X}+B)$ is a Cartier divisor for some integer $m>0$,

(2) for every proper birational morphism $\pi:X'\rightarrow X$ with
normal $X'$, 
\begin{eqnarray*}
K_{X'}+\pi_{*}^{-1}B & = & \pi^{*}(K_{X}+B)+\sum a_{i}E_{i}
\end{eqnarray*}
one has $a_{i}\geq-1$ . Here the $E_{i}$'s are the irreducible exceptional
divisors of $\pi$, and the pullback $\pi^{*}$ is defined by extending
$\mathbb{Q}$-linearly the pullback on Cartier divisors; $\pi_{*}^{-1}B$
is the strict preimage of $B$. \\

\begin{defn}
A pair $(X,B)$ is called \emph{semi log canonical}\textbf{ }\emph{(slc)
}if
\end{defn}
(1) $X$ satisfies Serre's condition $S_{2}$,

(2) $X$ is g.d.c., and no divisor $B_{i}$ contains any component
of the double locus of $X$,

(3) $m(K_{X}+B)$ is a Cartier divisor for some integer $m>0$,

(4) for the normalization $\nu:X^{\nu}\rightarrow X$ , the pair $(X^{\nu},(\mbox{double locus})+\nu_{*}^{-1}B)$
is log canonical.\\

\begin{defn}
Let $(X,B)$ be a semi log canonical pair and $f:X\rightarrow S$
a proper morphism. A pair $(X^{c},B^{c})$ sitting in a diagram \\
\[
\begin{array}{ccccc}
X &  & \overset{\phi}{\dashrightarrow} &  & X^{c}\\
 & _{f}\searrow &  & \swarrow_{f^{c}}\\
 &  & S
\end{array}
\]
 is called a \emph{log canonical model} of $(X,B)$ if
\end{defn}
(1) $f^{c}$ is proper,

(2) $\phi$ is a birational contraction, that is, $\phi^{-1}$ has
no exceptional divisors,

(3) $B^{c}=\phi_{*}B$,

(4) $K_{X^{c}}+B^{c}$ is $\mathbb{Q}$-Cartier and $f^{c}$-ample,
and 

(5) $a(E,X,B)\leq a(E,X^{c},B^{c})$ for every $\phi$-exceptional
divisor $E\subset X$.\\

\begin{defn}
The pair $(X,B)$ is called \emph{stable} if it satisfies the following
conditions
\end{defn}
(1) on singularities: the pair $(X,B)$ is semi log canonical, and 

(2) numerical: the divisor $K_{X}+B$ is ample.\\

Let $\beta=(b_{1},...,b_{n})$, $0<b_{i}\leq1$, $b_{i}\in\mathbb{Q}$
be a weight. A hyperplane arrangement is a pair $\left(\mathbb{P}^{r-1},\sum b_{i}B_{i}\right)$
with weight $\beta$, where $B_{1},...,B_{n}$ are hyperplanes in
$\mathbb{P}^{r-1}$. The pair $(\mathbb{P}^{r-1},\sum b_{i}B_{i})$
is \emph{lc }if for each intersection $\cap_{i\in I}B_{i}$ of codimension
$k$, one has $\sum_{i\in I}b_{i}\leq k$, where $I\subset\{1,...,n\}$.
The pair $(\mathbb{P}^{r-1},\sum b_{i}B_{i})$ is stable if and only
if it is\emph{ }lc (slc being an analog of lc for nonnormal pairs)
and $\left|\beta\right|=\sum_{i=1}^{n}b_{i}>r$.

\subsection{Abelian covers\label{subsecabcover}}

We will recall some definitions and theorems from \cite{Par91,AP12}
first.
\begin{defn}
Let $G$ be a finite abelian group. An \emph{abelian cover} with Galois
group $G$, or a $G$-cover, is a finite morphism $\pi:X\rightarrow Y$
of varieties which is the quotient map for a generically faithful
action of a finite abelian group $G$. 

An \emph{isomorphism of $G$-covers} $\pi_{1}:X_{1}\rightarrow Y$,
$\pi_{2}:X_{2}\rightarrow Y$ is an isomorphism $\phi:X_{1}\rightarrow X_{2}$
such that $\pi_{1}=\pi_{1}\circ\phi$.
\end{defn}
Let $Y$ be a smooth variety and $X$ be a normal variety. Let $G$
be a finite abelian group and $G^{*}=\mbox{Hom}(G,\mathbb{C}^{*})$
is the group of characters of $G$. The $G$-action on $X$ with $X/G=Y$
is equivalent to the decomposition: 
\begin{eqnarray*}
\pi_{*}\mathcal{O}_{X} & = & \underset{\chi\in G^{*}}{\oplus}\mathcal{L}_{\chi}^{-1},\mbox{ }\mbox{ }\mathcal{L}_{1}=\mathcal{O}_{Y}
\end{eqnarray*}
where the $\mathcal{L}_{\chi}$ are line bundles on $Y$ and $G$
acts on $\mathcal{L}_{\chi}^{-1}$ via the character $\chi$. 

In this paper we will only discuss the case when $G=\mathbb{Z}_{2}^{r}$.
A set of \emph{building data} $(L_{\chi},D_{g})$ for the case $G=\mathbb{Z}_{2}^{r}$
described in \cite{Par91} can be simplified as 
\begin{itemize}
\item effective Cartier divisors $D_{g}$, $g\in G\setminus\{0\}$ (possibly
not distinct),
\item line bundles $L_{\chi}$, $\chi\in G^{*}$.
\end{itemize}
Moreover the building data for the case $G=\mathbb{Z}_{2}^{r}$ need
only satisfy the \emph{fundamental relations}:

\begin{eqnarray*}
L_{\chi}+L_{\chi'} & \equiv & L_{\chi\chi'}+\sum_{g\in G}\epsilon_{g}^{\chi,\chi'}D_{g}
\end{eqnarray*}
where $\epsilon_{g}^{\chi,\chi'}=1$ if both $\chi(g)=\chi'(g)=-1$
and $\epsilon_{g}^{\chi,\chi'}=0$ otherwise.

In particular, let $G=\mathbb{Z}_{2}^{2}=\{e,a,b,c\}$ and $G^{*}=\{\chi_{0},\chi_{1},\chi_{2},\chi_{3}\}$
be the character group with $\chi_{0}\equiv1$ , $\chi_{1}(b)=\chi_{1}(c)=-1$,
$\chi_{2}(a)=\chi_{2}(c)=-1$, $\chi_{3}(a)=\chi_{3}(b)=-1$, and
assume that $\mbox{Pic }Y$ has no 2-torsion. Then the building data
only needs to satisfy 

\begin{eqnarray*}
2L_{\chi_{1}} & = & D_{b}+D_{c}\\
2L_{\chi_{2}} & = & D_{a}+D_{c}\\
2L_{\chi_{3}} & = & D_{a}+D_{b}
\end{eqnarray*}

The general theory of abelian covers was extended to the case of non-normal
$X$ in \cite{AP12}; it is used in \cite{AP09}. For details of the
abelian covers for the case of non-normal $X$ we will refer to \cite{AP12}.
Now we will recall a theorem in \cite{AP12} which is needed for our
paper. 

For every building data $(L_{\chi},D_{g})$, \cite[Def. 2.2]{Par91}
defines a \emph{standard abelian cover} explicitly, by equations.
\begin{defn}
For a standard $G$-cover $\pi:X\rightarrow Y$, the \emph{Hurwitz
divisor }of $\pi$ is the $\mathbb{Q}$-divisor $D_{Hur}:=\sum_{i}\frac{m_{i}-1}{m_{i}}D_{i}$,
where $m_{i}$ is the ramification index of $D_{i}$. 

The Hurwitz formula 

\begin{eqnarray*}
K_{X} & \sim_{\mathbb{Q}} & \pi^{*}(K_{Y}+D_{Hur})
\end{eqnarray*}
shows that $X$ is of general type if and only if $K_{Y}+D_{Hur}$
is big.\end{defn}
\begin{thm*}
\cite[Proposition 2.5]{AP12}. Let $\pi:X\rightarrow Y$ be a $G$-cover
and let $D$ be the Hurwitz divisor of $\pi$. Then 

(i) The divisor $K_{X}$ is $\mathbb{Q}$-Cartier if and only if $K_{Y}+D$
is $\mathbb{Q}$-Cartier.

(ii) $K_{X}=\pi^{*}(K_{Y}+D)$.

(iii) The variety $X$ is slc if and only if so is the pair $\left(Y,D\right)$.\end{thm*}
\begin{cor}
For a $G$-cover $\pi:X\rightarrow Y$ with Hurwitz divisor $D$,
$X$ is stable if and only if the pair $(Y,D)$ is stable. 
\end{cor}

\subsection{The construction of the compactified moduli space $\overline{M}_{Bur}^{d}$. }

The compactified moduli space $\overline{M}_{Bur}^{6}$ is constructed
in \cite[Section 5.3]{AP09} as an adaption of the construction of
the moduli space $\overline{M}_{\mathbf{b}}(3,9)$ of weighted hyperplane
arrangements of 9 lines in $\mathbb{P}^{2}$ with weight $\mathbf{b}=\left(\frac{1}{2},...,\frac{1}{2}\right)$.
This construction carries over verbatim to the $K^{2}\leq5$ case.
We refer to \cite{AP09} for details.

Fix weight $\mathbf{b}=\left(\frac{1}{2},...,\frac{1}{2}\right)$
and a polytope $\triangle_{Bur}^{d}$ (see Section \ref{matroid}).
We define $ $$\overline{M}_{Bur}^{d}$ to be the moduli space of
stable toric varieties over $\mbox{G}_{Bur,\mathbf{b}}^{d}$ of topological
type $\triangle_{Bur}^{d}$, where $\mbox{G}_{Bur,\mathbf{b}}^{d}$
is the $\mathbf{b}$-cut of certain subvariety $\mbox{G}_{Bur}^{d}\subset\mbox{G}(3,9)$
(see \cite[Section 5.3]{AP09}). Thus $\overline{M}_{Bur}^{d}$ parametrizes
stable toric varieties $Z\rightarrow\mbox{G}_{Bur,\mathbf{b}}^{d}$,
and the moment polytopes of the irreducible components of $Z=\cup Z_{s}$
give a tiling of $\triangle_{Bur}^{d}$. For a stable toric variety
$Z\rightarrow\mbox{G}_{Bur,\mathbf{b}}^{d}$, one recovers the stable
pair $(Y,\dfrac{1}{2}D)$ as a GIT quotient $P_{Bur,Z}^{d}//_{\mathbf{b}}T$,
where $P_{Bur,Z}^{d}=P\times_{\mbox{G}_{Bur,\mathbf{b}}^{d}}Z$ is
the pullback of the universal family $P$.

\section{Burniat surfaces with $K^{2}=5$ \label{degree5}}

\subsection{Burniat surfaces with $K^{2}=5$.}

We will use the construction of Burniat surfaces in \cite{Pet77}.
To construct a Burniat surface $X$ with $K_{X}^{2}=5$, we start
with an arrangement of $9$ distinct lines $A_{0}$, $A_{1}$, $A_{2}$,
$B_{0}$, $B_{1}$, $B_{2}$, $C_{0}$, $C_{1}$, $C_{2}$ in $\mathbb{P}^{2}$.
The lines $A_{0}$, $B_{0}$, $C_{0}$ form a non-degenerate triangle
with the vertices $P_{A}$, $P_{B}$, $P_{C}$. Lines $A_{1},A_{2}$
pass through $P_{B},$ $B_{1},B_{2}$ pass through $P_{C}$, and $C_{1},C_{2}$
pass through $P_{A}$. Moreover, $A_{1},B_{1},C_{1}$ meet at one
point $P$. The other lines are in general position otherwise.

Blow up $\mathbb{P}^{2}$ at $P_{A},P_{B},P_{C},P$. We denote the
exceptional divisors on $Bl_{4}\mathbb{P}^{2}$ by $A_{3},B_{3},C_{3},E$
and by $A_{i},B_{i},C_{i},\mbox{ }i=0,1,2$ the strict preimages of
$A_{i},B_{i},C_{i}$ on $\mathbb{P}^{2}$ . The blowup morphism is
as follows\\

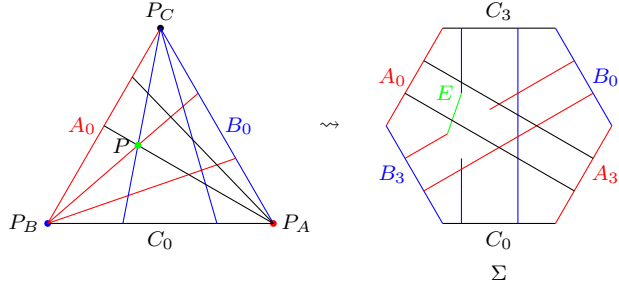
\begin{center}
\begin{tikzpicture}[xscale=1.5,yscale=1.5][font=\tiny]
      \coordinate (P_A) at (1,0);
      \node[right] at (P_A){$P_A$};
      \path[fill=red](P_A)circle[radius=0.03];
      \coordinate (P_B) at (-1,0);
      \node[left] at (P_B){$P_B$};
      \path[fill=blue](P_B)circle[radius=0.03];
      \coordinate (P_C) at (0,3^0.5);
      \node[above] at (P_C){$P_C$};
      \path[fill=black](P_C)circle[radius=0.03];
      
      \draw[black](P_A)--(P_B);
      \node[below,black]at($(P_A)!1/2!(P_B)$){$C_0$};
      \draw[red](P_B)--(P_C);
      \node[left,red]at($(P_B)!1/2!(P_C)$){$A_0$};
      \draw[blue](P_C)--(P_A);
      \node[right,blue]at($(P_C)!1/2!(P_A)$){$B_0$};
      \draw [name path=A_1,red]($(P_C)!1/3!(P_A)$)--(P_B);
      \draw [name path=A_2,red]($(P_C)!2/3!(P_A)$)--(P_B);
      \draw [name path=B_1,blue]($(P_A)!2/3!(P_B)$)--(P_C);
      \draw [name path=B_2,blue]($(P_A)!1/4!(P_B)$)--(P_C);
      \draw [name path=C_1,black]($(P_B)!1/2!(P_C)$)--(P_A);
      \draw [name path=C_2,black]($(P_B)!3/4!(P_C)$)--(P_A);
      \path [name intersections={of=A_1 and B_1}]; 
      \coordinate [label=left:$P$] (P) at (intersection-1);
      \path[fill=green](P)circle[radius=0.03];
\node at (1.5,3^0.5/2){$\rightsquigarrow$};
      \coordinate (f) at (1/2+3,0);
      \coordinate (a) at (1+3,3^0.5/2);
      \coordinate (b) at (1/2+3,3^0.5);
      \coordinate (c) at (-1/2+3,3^0.5);
      \coordinate (d) at (-1+3,3^0.5/2);
      \coordinate (e) at (-1/2+3,0);
      \draw[blue](a)--(b);
      \draw[black](b)--(c);
      \draw[red](c)--(d);
      \draw[blue](d)--(e);
      \draw[black](e)--(f);
      \draw[red](f)--(a);
      \coordinate (g) at ($(a)!1/3!(b)$);
      \coordinate (h) at ($(a)!2/3!(b)$);
      \coordinate (i) at ($(b)!1/3!(c)$);
      \coordinate (j) at (-1/3+3,3^0.5);
      \coordinate (k) at ($(c)!1/3!(d)$);
      \coordinate (l) at ($(c)!2/3!(d)$);
      \coordinate (m) at ($(d)!1/3!(e)$);
      \coordinate (n) at ($(d)!2/3!(e)$);
      \coordinate (o) at (-1/3+3,0);
      \coordinate (p) at ($(e)!2/3!(f)$);
      \coordinate (q) at ($(m)!1/4!(h)$);
      \coordinate (t) at ($(m)!1/2!(h)$);
      \coordinate (r) at ($(j)!1/3!(o)$);
      \coordinate (s) at ($(j)!2/3!(o)$);
      \coordinate (u) at ($(a)!1/3!(f)$);
      \coordinate (v) at ($(a)!2/3!(f)$);
      \draw[blue] (j)--(r);
      \draw[blue](s)--(o);
      \draw[blue](i)--(p);
      \draw[red](m)--(q);
      \draw[red](t)--(h);
      \draw[red](n)--(g);
      \draw[green](q)--(r);
      \draw[black](k)--(u);
      \draw[black](l)--(v);
      \node[above,black]at($(b)!1/2!(c)$){$C_3$};
      \node[below,black]at($(e)!1/2!(f)$){$C_0$};
      \node[right,blue]at($(a)!1/2!(b)$){$B_0$};
      \node[left,blue]at($(d)!1/2!(e)$){$B_3$};
      \node[left,red]at($(c)!1/2!(d)$){$A_0$};
      \node[right,red]at($(a)!1/2!(f)$){$A_3$};
      \node[left,green]at($(j)!1/3!(o)$){$E$};
\node[below,black]at($(e)+(0.5,-0.3)$){$\Sigma$};
\end{tikzpicture}

\par\end{center}
\begin{defn}
A Burniat surface $X$ with $K_{X}^{2}=5$ is a $\mathbb{Z}_{2}^{2}$-cover
of $\Sigma=\mbox{Bl}_{4}\mathbb{P}^{2}$ for the building data $D_{a}=\sum_{i=0}^{3}A_{i}$,
$D_{b}=\sum_{i=0}^{3}B_{i}$, $D_{c}=\sum_{i=0}^{3}C_{i}$, where
$a,b,c$ are the $3$ nonzero elements of $\mathbb{\mathbb{Z}}_{2}^{2}$. 
\end{defn}
For a $\mathbb{Z}_{2}^{2}$-cover, the Hurwitz divisor is $D=\dfrac{1}{2}\left(D_{a}+D_{b}+D_{c}\right)$.
Using the Riemann-Hurwitz formula $K_{X}=\pi^{*}(K_{\Sigma}+D)$,
we have 
\[
K_{X}^{2}=\left(\pi^{*}(K_{\Sigma}+D)\right)^{2}=4\cdot(K_{\Sigma}+D)^{2}=4\left(-\dfrac{1}{2}K_{\Sigma}\right)^{2}=5.
\]

By the theorem in Section \ref{subsecabcover}, we can reduce the
problem of compactifying the moduli space of stable Burniat surfaces
with $K^{2}=5$ to compactifying the moduli space of stable pairs
$\left(\Sigma,D\right)$ described above.

\subsection{Degenerations of Burniat surfaces with $K^{2}=5$. }

We consider degenerations of Burniat arrangements of curves on $\Sigma=\mbox{Bl}_{4}\mathbb{P}^{2}$.
When the arrangement on $\Sigma$ is not log canonical, choose a generic
one-parameter family of degenerating arrangements on $\Sigma$ degenerating
to it. Then the limit stable surface splits into several irreducible
components. Below, we consider such generic degenerations. Let $\mathcal{Y}$
be the total space of the one parameter family of surfaces isomorphic
to $\Sigma$ with the central fiber being the degenerating arrangement.
Write $\Sigma_{0}$ for the central fiber of $\mathcal{Y}$.\\

\textbf{Case 1}. When the curve $A_{2}$ degenerates to $A_{0}+C_{3}$,
$B_{2}$ degenerates to $A_{0}+B_{3}$, and $C_{2}$ degenerates to
$B_{3}+C_{0}$ (the first figure below). Let $L_{P}$ be the curve
in $\mathcal{Y}$ consisting of the points $P$ in each fiber, which
is the intersection of the curves $A_{1},B_{1},C_{1}$. We first blow
up the total space along $L_{P}$, then blow up the resulting total
space along $A_{0}$ in the central fiber. The central fiber $\Sigma_{0}$
becomes $Bl_{4}\mathbb{P}^{2}\cup\mathbb{F}_{1}$ (the second figure
below), where $A_{0}$ is the (-1)-curve in $\mathbb{F}_{1}$. Finally
we blow up the total space along the proper transform of $B_{3}$
in the component $Bl_{4}\mathbb{P}^{2}$ of $\Sigma_{0}$. The resulting
central fiber is a union of three components $Bl_{4}\mathbb{P}^{2}\cup Bl_{1}\mathbb{F}_{1}\cup\mathbb{F}_{0}$. 

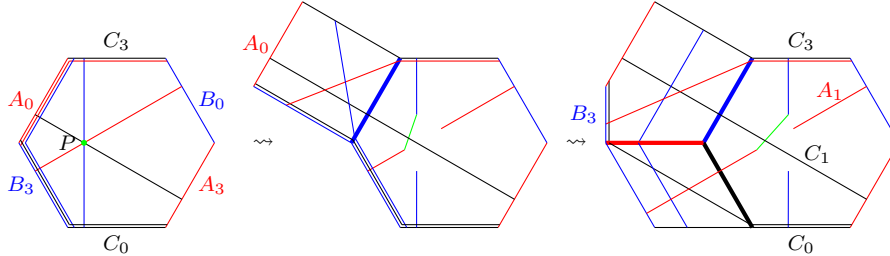
\begin{tikzpicture}[xscale=1.3,yscale=1.3][font=\tiny]
      \coordinate (a) at (1,3^0.5/2);
      \coordinate (b) at (1/2,3^0.5);
      \coordinate (c) at (-1/2,3^0.5);
      \coordinate (d) at (-1,3^0.5/2);
      \coordinate (e) at (-1/2,0);
      \coordinate (f) at (1/2,0);
      \draw[blue](a)--(b);
      \draw[black](b)--(c);
      \draw[red](c)--(d);
      \draw[blue](d)--(e);
      \draw[black](e)--(f);
      \draw[red](f)--(a);
      \coordinate (g) at ($(a)!1/3!(b)$);
      \coordinate (h) at ($(a)!2/3!(b)$);
      \coordinate (i) at ($(b)!1/3!(c)$);
      \coordinate (j) at (-1/3,3^0.5);
      \coordinate (k) at ($(c)!1/3!(d)$);
      \coordinate (l) at ($(c)!2/3!(d)$);
      \coordinate (m) at ($(d)!1/3!(e)$);
      \coordinate (n) at ($(d)!2/3!(e)$);
      \coordinate (o) at (-1/3,0);
      \coordinate (p) at ($(e)!2/3!(f)$);
      \coordinate (q) at ($(m)!1/4!(h)$);
      \coordinate (t) at ($(m)!1/2!(h)$);
      \coordinate (r) at ($(j)!1/3!(o)$);
      \coordinate (s) at ($(j)!2/3!(o)$);
      \coordinate (u) at ($(a)!1/3!(f)$);
      \coordinate (v) at ($(a)!2/3!(f)$);
      \draw[blue,name path=B_1] (j)--(o);
      \draw[red,name path=A_1](m)--(h);
      \draw[black,name path=C_1](l)--(v);
      \node[above,black]at($(b)!1/2!(c)$){$C_3$};
      \node[below,black]at($(e)!1/2!(f)$){$C_0$};
      \node[right,blue]at($(a)!1/2!(b)$){$B_0$};
      \node[left,blue]at($(d)!1/2!(e)$){$B_3$};
      \node[left,red]at($(c)!1/2!(d)$){$A_0$};
      \node[right,red]at($(a)!1/2!(f)$){$A_3$};
      \path [name intersections={of=A_1 and B_1}]; 
      \coordinate (P) at (intersection-1);
      \path[fill=green](P)circle[radius=0.03];
      \node[left,black] at(P){$P$};
      \draw[red]($(d)!1/30!(e)$)--($(c)!1/60!(f)$);
      \draw[red]($(c)!1/60!(f)$)--($(b)!1/30!(a)$);
      \draw[blue]($(e)!2/30!(f)$)--($(d)!2/60!(a)$);
      \draw[blue]($(d)!2/60!(a)$)--($(c)!2/30!(b)$);
      \draw[black]($(d)!1/30!(c)$)--($(e)!1/60!(b)$);
      \draw[black]($(e)!1/60!(b)$)--($(f)!1/30!(a)$);
  \node at ($(a)+(0.5,0)$){$\rightsquigarrow$};
      \coordinate (a) at (1+3.4,3^0.5/2);
      \coordinate (b) at (1/2+3.4,3^0.5);
      \coordinate (c) at (-1/2+3.4,3^0.5);
      \coordinate (d) at (-1+3.4,3^0.5/2);
      \coordinate (e) at (-1/2+3.4,0);
      \coordinate (f) at (1/2+3.4,0);
      \draw[blue](a)--(b);
      \draw[black](b)--(c);
      \draw[ultra thick,blue](c)--(d);
      \draw[blue](d)--(e);
      \draw[black](e)--(f);
      \draw[red](f)--(a);
      \coordinate (g) at ($(a)!1/3!(b)$);
      \coordinate (h) at ($(a)!2/3!(b)$);
      \coordinate (i) at ($(b)!1/3!(c)$);
      \coordinate (j) at (-1/3+3.4,3^0.5);
      \coordinate (k) at ($(c)!1/3!(d)$);
      \coordinate (l) at ($(c)!2/3!(d)$);
      \coordinate (m) at ($(d)!1/3!(e)$);
      \coordinate (n) at ($(d)!2/3!(e)$);
      \coordinate (o) at (-1/3+3.4,0);
      \coordinate (p) at ($(e)!2/3!(f)$);
      \coordinate (q) at ($(m)!1/4!(h)$);
      \coordinate (t) at ($(m)!1/2!(h)$);
      \coordinate (r) at ($(j)!1/3!(o)$);
      \coordinate (s) at ($(j)!2/3!(o)$);
      \coordinate (u) at ($(a)!1/3!(f)$);
      \coordinate (v) at ($(a)!2/3!(f)$);
      \draw[blue] (j)--(r);
      \draw[blue](s)--(o);
      \draw[red](m)--(q);
      \draw[red](t)--(h);
      \draw[green](q)--(r);
      \draw[black](l)--(v);
      \draw[red]($(c)!1/30!(d)$)--($(b)!1/30!(a)$);
      \draw[blue]($(e)!2/30!(f)$)--($(d)!2/30!(c)$);
      \draw[black]($(d)!1/30!(c)$)--($(e)!1/60!(b)$);
      \draw[black]($(e)!1/60!(b)$)--($(f)!1/30!(a)$);
      \coordinate (w) at ($(a)!1+2/3!(c)$);
      \coordinate (x) at ($(f)!1+2/3!(d)$);
      \draw[black](c)--(w);
      \draw[black](l)--($(v)!1+2/3!(l)$);
      \draw[black]($(d)!1/30!(c)$)--($(x)!1/30!(w)$);
      \draw[red](w)--(x);
      \draw[red]($(c)!1/30!(d)$)--($(d)!2/3!(x)$);
      \draw[blue]($(d)!2/30!(c)$)--($(c)!2/3!(w)$);
      \draw[blue](d)--(x);
 \node at ($(a)+(0.3,0)$){$\rightsquigarrow$};
 \node[left,red] at ($(w)!1/2!(x)$){$A_0$};
      \coordinate (a) at (1+3+4,3^0.5/2);
      \coordinate (b) at (1/2+3+4,3^0.5);
      \coordinate (c) at (-1/2+3+4,3^0.5);
      \coordinate (d) at (-1+3+4,3^0.5/2);
      \coordinate (e) at (-1/2+3+4,0);
      \coordinate (f) at (1/2+3+4,0);
      \draw[blue](a)--(b);
      \draw[black](b)--(c);
      \draw[ultra thick,blue](c)--(d);
      \draw[ultra thick,black](d)--(e);
      \draw[black](e)--(f);
      \draw[red](f)--(a);
      \coordinate (g) at ($(a)!1/3!(b)$);
      \coordinate (h) at ($(a)!2/3!(b)$);
      \coordinate (i) at ($(b)!1/3!(c)$);
      \coordinate (j) at (-1/3+3+4.2,3^0.5);
      \coordinate (k) at ($(c)!1/3!(d)$);
      \coordinate (l) at ($(c)!2/3!(d)$);
      \coordinate (m) at ($(d)!1/3!(e)$);
      \coordinate (n) at ($(d)!2/3!(e)$);
      \coordinate (o) at (-1/3+3+4.2,0);
      \coordinate (p) at ($(e)!2/3!(f)$);
      \coordinate (q) at ($(m)!1/4!(h)$);
      \coordinate (t) at ($(m)!1/2!(h)$);
      \coordinate (r) at ($(j)!1/3!(o)$);
      \coordinate (s) at ($(j)!2/3!(o)$);
      \coordinate (u) at ($(a)!1/3!(f)$);
      \coordinate (v) at ($(a)!2/3!(f)$);
      \draw[blue,name path=B_1] (j)--(r);
      \draw[blue,name path=B_1'](s)--(o);
      \draw[red,name path=A_1](m)--(q);
      \draw[red,name path=A_1'](t)--(h);
      \draw[green,name path=E](q)--(r);
      \draw[black,name path=C_1](l)--(v);
      \draw[red]($(c)!1/30!(d)$)--($(b)!1/30!(a)$);
      \draw[black]($(e)!1/30!(d)$)--($(f)!1/30!(a)$);
      \node[black,above] at ($(b)!1/2!(c)$){$C_3$};
      \node[black,below] at ($(e)!1/2!(f)$){$C_0$};
      \node[black,above] at ($(v)!1/3!(l)$){$C_1$};
      \node[red,above] at ($(h)!1/4!(m)$){$A_1$};
      \coordinate (w) at ($(a)!1+2/3!(c)$);
      \coordinate (x) at ($(f)!1+2/3!(d)$);
      \coordinate (y) at ($(a)!1+1/2!(d)$);
      \draw[black](c)--(w);
      \draw[black](l)--($(v)!1+2/3!(l)$);
      \draw[black]($(x)!2/30!(w)$)--($(y)!1/30!(d)$);
      \draw[blue]($(c)!1/2!(w)$)--($(d)!2/3!(y)$);
      \draw[red](w)--(x);
      \draw[red]($(c)!1/30!(d)$)--($(x)!2/3!(y)$);
      \draw[red,ultra thick](d)--(y);
      \draw[blue,name path=B_3](x)--(y);
      \node[blue,left] at ($(x)!1/2!(y)$){$B_3$};
      \coordinate (z) at ($(f)!2!(e)$);
      \draw[blue](y)--(z);
      \draw[black](e)--(z);
      \draw[black]($(y)!1/30!(d)$)--($(e)!1/30!(d)$);
      \draw[blue]($(d)!2/3!(y)$)--($(z)!1/3!(e)$);
      \draw[red]($(m)$)--($(h)!1+1/2!(m)$);      
\end{tikzpicture}     

We can use the triple point formula to compute $\left(K_{\Sigma_{0}}+D\right)|_{Y_{i}}.C$,
where $Y_{i}$ is a component of $\Sigma_{0}$ and $C$ is a curve
in the component $Y_{i}$. 

Let us recall the triple point formula: let $\Sigma_{0}=\cup Y_{i}$
be the central fiber in a smooth one-parameter family, and assume
that $\Sigma_{0}$ is reduced and has simple normal crossing. Let
$C$ be the intersection $Y_{i}\cap Y_{j}$ and assume that it is
a smooth curve. Denote by $p_{3}$ the number of the triple points
of $\Sigma_{0}$ contained in $C$, then 
\begin{eqnarray*}
(C|_{Y_{i}})^{2}+(C|_{Y_{j}})^{2}+p_{3} & = & 0.
\end{eqnarray*}
By the adjunction formula, we also have 
\begin{eqnarray*}
\left(K_{\Sigma_{0}}+D\right)|_{Y_{i}} & = & K_{Y_{i}}+D|_{Y_{i}}+(\mbox{the double locus}).
\end{eqnarray*}

The intersection number $\left(K_{\Sigma_{0}}+D\right)|_{Bl_{4}\mathbb{P}^{2}}.C$
is 0 when the curve $C$ is $A_{1}$, $C_{0}$, $C_{1}$ or $C_{3}$,
and positive for the other curves in $Bl_{4}\mathbb{P}^{2}$. In the
component $Bl_{1}\mathbb{F}_{1}$, $\left(K_{\Sigma_{0}}+D\right)|_{Bl_{1}\mathbb{F}_{1}}.B_{3}=0$
and $\left(K_{\Sigma_{0}}+D\right)|_{Bl_{1}\mathbb{F}_{1}}.C>0$ for
the other curves. We also have $\left(K_{\Sigma_{0}}+D\right)|_{\mathbb{F}_{0}}.C>0$
for all the curves in the component $\mathbb{F}_{0}$. Thus $K_{\Sigma_{0}}+D$
is big, nef and vanishes on $B_{3},C_{1}$ and $C_{3}$. The 3-fold
is the minimal model of the degenerate family. Using the inversion
of adjunction in \cite{Ka07}, we see that the pair $\left(\mathcal{Y},\mathfrak{\mathcal{D}}\right)$
is log terminal and $\mathcal{D}$ is an effective divisor on $\mathfrak{X}$
such that $K_{\mathcal{Y}}+\mathcal{D}$ is nef and big. By Base Point
Free theorem, the linear system $\left|n\left(K_{\mathfrak{\mathcal{Y}}}+\mathcal{D}\right)\right|$
is base point free for all sufficiently large $n\in\mathbb{N}$. Then
we can define a birational morphism by the linear system $\left|n\left(K_{\mathfrak{\mathcal{Y}}}+\mathcal{D}\right)\right|$,
which contracts $A_{1}B_{3},C_{0},C_{1},C_{3}$ labeled in the third
figure above. The image of the birational morphism is the lc model
of the degenerate family. \\

The surface $Bl_{4}\mathbb{P}^{2}$ becomes $\mathbb{P}^{2}$ after
contracting $A_{1},C_{0},C_{1},C_{3}$. The component $Bl_{1}\mathbb{F}_{1}$
becomes $\mathbb{F}_{0}$ after contracting $B_{3}$. The central
fiber of the resulting log canonical model is $\mathbb{F}_{0}\cup\mathbb{F}_{0}\cup\mathbb{P}^{2}$,
which is the first figure below. For $\mathbb{P}^{1}\times\mathbb{P}^{1}$,
there is a further degeneration that splits to $\mathbb{P}^{2}\cup\mathbb{P}^{2}$.
We list the three possible further degenerations below which are the
rest three figures. The second and third figures differ only by a
permutation of colors. Thus there are only two different degenerations,
we call them Case 2 and Case 3. 

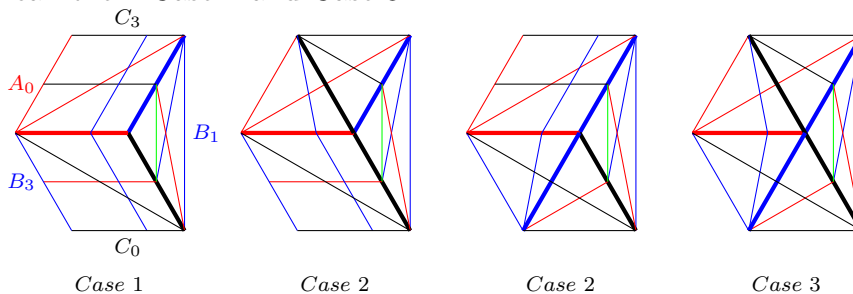
\begin{tikzpicture}[xscale=1.5,yscale=1.5][font=\tiny]
      \coordinate (O) at (0,0);
      \coordinate (a) at (1/2,3^0.5/2);
      \coordinate (b) at (-1/2,3^0.5/2);
      \coordinate (c) at (-1,0);
      \coordinate (d) at (-1/2,-3^0.5/2);
      \coordinate (e) at (1/2,-3^0.5/2);
      \draw[blue,ultra thick] (O)--(a);
      \draw[blue] ($(O)!1/2!(e)$)--(a);
      \draw[blue] (a)--(e);
      \draw[black,ultra thick] (O)--(e);
      \draw[red] ($(O)!1/2!(a)$)--(e);
      \draw[green] ($(O)!1/2!(a)$)--($(O)!1/2!(e)$);
      \node [blue,right] at ($(a)!1/2!(e)$){$B_1$};
      \draw[red,ultra thick] (O)--(c);
      \draw[red] (a)--(c);
      \draw[red] (c)--(b);
      \draw[black] ($(O)!1/2!(a)$)--($(c)!1/2!(b)$);
      \draw[black] (a)--(b);
      \draw[blue] ($(O)!1/3!(c)$)--($(a)!1/3!(b)$); 
      \node [black,above] at ($(a)!1/2!(b)$){$C_3$};
      \node [red,left] at ($(b)!1/2!(c)$){$A_0$};
      \draw[red] ($(c)!1/2!(d)$)--($(O)!1/2!(e)$);
      \draw[blue] (c)--(d);
      \draw[blue] ($(O)!1/3!(c)$)--($(e)!1/3!(d)$);
      \draw[black] (c)--(e);
      \draw[black] (d)--(e);
      \node [blue,left]at ($(c)!1/2!(d)$){$B_3$};
      \node [black,below] at ($(d)!1/2!(e)$){$C_0$};
      \node [black,below] at ($(d)+(1/3,-1/3)$){$Case$ $1$};
      \coordinate (O) at (0+2,0);
      \coordinate (a) at (1/2+2,3^0.5/2);
      \coordinate (b) at (-1/2+2,3^0.5/2);
      \coordinate (c) at (-1+2,0);
      \coordinate (d) at (-1/2+2,-3^0.5/2);
      \coordinate (e) at (1/2+2,-3^0.5/2);
      \draw[blue,ultra thick] (O)--(a);
      \draw[blue] ($(O)!1/2!(e)$)--(a);
      \draw[blue] (a)--(e);
      \draw[black,ultra thick] (O)--(e);
      \draw[red] ($(O)!1/2!(a)$)--(e);
      \draw[green] ($(O)!1/2!(a)$)--($(O)!1/2!(e)$);
      \draw[red,ultra thick] (O)--(c);
      \draw[red] (a)--(c);
      \draw[red] (c)--(b);
      \draw[black] ($(O)!1/2!(a)$)--(b);
      \draw[black] (a)--(b);
      \draw[black,ultra thick] (O)--(b);
      \draw[blue] ($(O)!1/3!(c)$)--(b); 
      \draw[red] ($(c)!1/2!(d)$)--($(O)!1/2!(e)$);
      \draw[blue] (c)--(d);
      \draw[blue] ($(O)!1/3!(c)$)--($(e)!1/3!(d)$);
      \draw[black] (c)--(e);
      \draw[black] (d)--(e);
      \node [black,below] at ($(d)+(1/3,-1/3)$){$Case$ $2$};
      \coordinate (O) at (0+4,0);
      \coordinate (a) at (1/2+4,3^0.5/2);
      \coordinate (b) at (-1/2+4,3^0.5/2);
      \coordinate (c) at (-1+4,0);
      \coordinate (d) at (-1/2+4,-3^0.5/2);
      \coordinate (e) at (1/2+4,-3^0.5/2);
      \draw[blue,ultra thick] (O)--(a);
      \draw[blue] ($(O)!1/2!(e)$)--(a);
      \draw[blue] (a)--(e);
      \draw[black,ultra thick] (O)--(e);
      \draw[red] ($(O)!1/2!(a)$)--(e);
      \draw[green] ($(O)!1/2!(a)$)--($(O)!1/2!(e)$);
      \draw[red,ultra thick] (O)--(c);
      \draw[red] (a)--(c);
      \draw[red] (c)--(b);
      \draw[black] ($(O)!1/2!(a)$)--($(c)!1/2!(b)$);
      \draw[black] (a)--(b);
      \draw[blue] ($(O)!1/3!(c)$)--($(a)!1/3!(b)$); 
      \draw[red] (d)--($(O)!1/2!(e)$);
      \draw[blue] (c)--(d);
      \draw[blue] ($(O)!1/3!(c)$)--(d);
      \draw[blue,ultra thick] (O)--(d);
      \draw[black] (c)--(e);
      \draw[black] (d)--(e);
      \node [black,below] at ($(d)+(1/3,-1/3)$){$Case$ $2$};
      \coordinate (O) at (0+6,0);
      \coordinate (a) at (1/2+6,3^0.5/2);
      \coordinate (b) at (-1/2+6,3^0.5/2);
      \coordinate (c) at (-1+6,0);
      \coordinate (d) at (-1/2+6,-3^0.5/2);
      \coordinate (e) at (1/2+6,-3^0.5/2);
      \draw[blue,ultra thick] (O)--(a);
      \draw[blue] ($(O)!1/2!(e)$)--(a);
      \draw[blue] (a)--(e);
      \draw[black,ultra thick] (O)--(e);
      \draw[red] ($(O)!1/2!(a)$)--(e);
      \draw[green] ($(O)!1/2!(a)$)--($(O)!1/2!(e)$);
      \draw[red,ultra thick] (O)--(c);
      \draw[red] (a)--(c);
      \draw[red] (c)--(b);
      \draw[black] ($(O)!1/2!(a)$)--(b);
      \draw[black] (a)--(b);
      \draw[black,ultra thick] (O)--(b);
      \draw[blue] ($(O)!1/3!(c)$)--(b); 
      \draw[red] (d)--($(O)!1/2!(e)$);
      \draw[blue] (c)--(d);
      \draw[blue] ($(O)!1/3!(c)$)--(d);
      \draw[blue,ultra thick] (O)--(d);
      \draw[black] (c)--(e);
      \draw[black] (d)--(e);
      \node [black,below] at ($(d)+(1/3,-1/3)$){$Case$ $3$};
\end{tikzpicture}     \\

Case 1 could be obtained from another degeneration when $B_{1}$ goes
to $A_{0}+B_{3}$ and $C_{1}$ degenerates to $B_{3}+C_{0}$ (the
first figure below). We first blow up the total space $\mathcal{Y}$
along the line $B_{3}$ and then blow up along the strict image of
$A_{0}$ in the component $Bl_{3}\mathbb{P}^{2}$ of the central fiber.
Finally we blow up the resulting total space along the proper transform
$\tilde{L}_{P}$ of the line $L_{P}$. The central fiber becomes $Bl_{3}\mathbb{P}^{2}\cup\mathbb{F}_{0}\cup Bl_{2}\mathbb{F}_{0}$
(the second figure blow). Running the minimal model program, we obtain
the lc model with the central fiber $\mathbb{P}^{2}\cup\mathbb{F}_{0}\cup\mathbb{F}_{0}$
(the third figure below), which is the same as Case 1 above by changing
the color of the building data due to the symmetry. Both degenerations
could come from Case 2 for $K^{2}=6$ in \cite{AP09}, with $A_{1},B_{1},C_{1}$
meeting at a point $P$. Case 2 could be obtained from the degeneration
of Case 7 for $K^{2}=6$ with the point $P$ on the boundary of the
hexagon. 

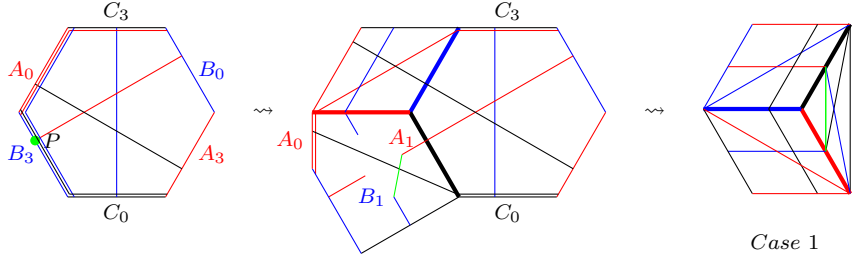
\begin{tikzpicture}[xscale=1.3,yscale=1.3][font=\tiny]
      \coordinate (a) at (1,3^0.5/2);
      \coordinate (b) at (1/2,3^0.5);
      \coordinate (c) at (-1/2,3^0.5);
      \coordinate (d) at (-1,3^0.5/2);
      \coordinate (e) at (-1/2,0);
      \coordinate (f) at (1/2,0);
      \draw[blue](a)--(b);
      \draw[black](b)--(c);
      \draw[red](c)--(d);
      \draw[blue](d)--(e);
      \draw[black](e)--(f);
      \draw[red](f)--(a);
      \coordinate (g) at ($(a)!1/3!(b)$);
      \coordinate (h) at ($(a)!2/3!(b)$);
      \coordinate (i) at ($(b)!1/3!(c)$);
      \coordinate (j) at ($(b)!1/2!(c)$);
      \coordinate (k) at ($(c)!1/3!(d)$);
      \coordinate (l) at ($(c)!2/3!(d)$);
      \coordinate (m) at ($(d)!1/3!(e)$);
      \coordinate (n) at ($(d)!2/3!(e)$);
      \coordinate (o) at ($(e)!1/2!(f)$);
      \coordinate (p) at ($(e)!2/3!(f)$);
      \coordinate (q) at ($(m)!1/4!(h)$);
      \coordinate (t) at ($(m)!1/2!(h)$);
      \coordinate (r) at ($(j)!1/3!(o)$);
      \coordinate (s) at ($(j)!2/3!(o)$);
      \coordinate (u) at ($(a)!1/3!(f)$);
      \coordinate (v) at ($(a)!2/3!(f)$);
      \draw[blue] (j)--(o);
      \draw[red](m)--(h);
      \draw[black](l)--(v);
      \node[above,black]at($(b)!1/2!(c)$){$C_3$};
      \node[below,black]at($(e)!1/2!(f)$){$C_0$};
      \node[right,blue]at($(a)!1/2!(b)$){$B_0$};
      \node[left,blue]at($(d)!1/2!(e)$){$B_3$};
      \node[left,red]at($(c)!1/2!(d)$){$A_0$};
      \node[right,red]at($(a)!1/2!(f)$){$A_3$};
      \node[right]at($(m)$){$P$};
      \path[fill=green](m)circle[radius=0.05];
      \draw[red]($(d)!1/30!(e)$)--($(c)!1/60!(f)$);
      \draw[red]($(c)!1/60!(f)$)--($(b)!1/30!(a)$);
      \draw[blue]($(e)!2/30!(f)$)--($(d)!2/60!(a)$);
      \draw[blue]($(d)!2/60!(a)$)--($(c)!2/30!(b)$);
      \draw[black]($(d)!1/30!(c)$)--($(e)!1/60!(b)$);
      \draw[black]($(e)!1/60!(b)$)--($(f)!1/30!(a)$);
  \node at ($(a)+(0.5,0)$){$\rightsquigarrow$};
      \coordinate (a) at (1+4,3^0.5/2);
      \coordinate (b) at (1/2+4,3^0.5);
      \coordinate (c) at (-1/2+4,3^0.5);
      \coordinate (d) at (-1+4,3^0.5/2);
      \coordinate (e) at (-1/2+4,0);
      \coordinate (f) at (1/2+4,0);
      \draw[blue](a)--(b);
      \draw[black](b)--(c);
      \draw[ultra thick,blue](c)--(d);
      \draw[ultra thick,black](d)--(e);
      \draw[black](e)--(f);
      \draw[red](f)--(a);
      \coordinate (g) at ($(a)!1/3!(b)$);
      \coordinate (h) at ($(a)!2/3!(b)$);
      \coordinate (i) at ($(b)!1/3!(c)$);
      \coordinate (j) at (-1/3+4.2,3^0.5);
      \coordinate (k) at ($(c)!1/3!(d)$);
      \coordinate (l) at ($(c)!2/3!(d)$);
      \coordinate (m) at ($(d)!1/3!(e)$);
      \coordinate (n) at ($(d)!2/3!(e)$);
      \coordinate (o) at (-1/3+4.2,0);
      \coordinate (p) at ($(e)!2/3!(f)$);
      \coordinate (q) at ($(m)!1/4!(h)$);
      \coordinate (t) at ($(m)!1/2!(h)$);
      \coordinate (r) at ($(j)!1/3!(o)$);
      \coordinate (s) at ($(j)!2/3!(o)$);
      \coordinate (u) at ($(a)!1/3!(f)$);
      \coordinate (v) at ($(a)!2/3!(f)$);
      \draw[blue] (j)--(o);
      \draw[red](m)--(h);
      \draw[black](l)--(v);
      \draw[red]($(c)!1/30!(d)$)--($(b)!1/30!(a)$);
      \draw[black]($(e)!1/30!(d)$)--($(f)!1/30!(a)$);
      \node[black,above] at ($(b)!1/2!(c)$){$C_3$};
      \node[black,below] at ($(e)!1/2!(f)$){$C_0$};
      \coordinate (w) at ($(a)!1+2/3!(e)$);
      \coordinate (x) at ($(b)!1+2/3!(d)$);
      \coordinate (y) at ($(a)!1+1/2!(d)$);
      \coordinate (w') at ($(e)!1/2!(w)$);
      \coordinate (y') at ($(d)!2/3!(y)$);
      \draw[black,name path=C_0](e)--(w);
      \draw[red,name path=A_2]($(x)!2/30!(w)$)--($(y)!1/30!(d)$);
      \draw[red,name path=A_1](m)--($(h)!1+1/6!(m)$);
      \draw[red,name path=A_1']($(h)!1+5/12!(m)$)--($(h)!1+2/3!(m)$);
      \draw[blue,name path=B_1](w')--($(w')!1/4!(y')$);
      \draw[blue,name path=B_1'](y')--($(y')!1/5!(w')$);  
      \draw[blue](w)--(x);
      \draw[black,name path=C_1]($(e)!1/30!(d)$)--($(x)!2/3!(y)$);
      \draw[red,ultra thick](d)--(y);
      \draw[red,name path=A_0](x)--(y);
      \draw[green,name path=E]($(h)!1+1/6!(m)$)--($(w')!1/4!(y')$);
      \node[red,left] at ($(x)!1/2!(y)$){$A_0$};
      \node[blue,left] at ($(w')!1/4!(y')$){$B_1$};
      \node[red,above] at ($(h)!1+1/6!(m)$){$A_1$};
      \coordinate (z) at ($(b)!2!(c)$);
      \draw[red,name path=A_0](y)--(z);
      \draw[black,name path=C_3](c)--(z);
      \draw[red,name path=A_2]($(y)!1/30!(d)$)--($(c)!1/30!(d)$);
      \draw[blue,name path=B_1]($(d)!2/3!(y)$)--($(z)!1/3!(c)$);
      \draw[black](l)--($(v)!1+1/2!(l)$);      
\node at ($(a)+(0.5,0)$){$\rightsquigarrow$};
      \coordinate (O) at (0+7,0+0.9);
      \coordinate (a) at (1/2+7,3^0.5/2+0.9);
      \coordinate (b) at (-1/2+7,3^0.5/2+0.9);
      \coordinate (c) at (-1+7,0+0.9);
      \coordinate (d) at (-1/2+7,-3^0.5/2+0.9);
      \coordinate (e) at (1/2+7,-3^0.5/2+0.9);
      \draw[black,ultra thick] (O)--(a);
      \draw[black] ($(O)!1/2!(e)$)--(a);
      \draw[black] (a)--(e);
      \draw[red,ultra thick] (O)--(e);
      \draw[blue] ($(O)!1/2!(a)$)--(e);
      \draw[green]($(O)!1/2!(e)$)--($(O)!1/2!(a)$);
      \draw[blue,ultra thick] (O)--(c);
      \draw[blue] (a)--(c);
      \draw[blue] (c)--(b);
      \draw[red] ($(O)!1/2!(a)$)--($(c)!1/2!(b)$);
      \draw[red] (a)--(b);
      \draw[black] ($(O)!1/3!(c)$)--($(a)!1/3!(b)$); 
      \draw[blue] ($(c)!1/2!(d)$)--($(O)!1/2!(e)$);
      \draw[black] (c)--(d);
      \draw[black] ($(O)!1/3!(c)$)--($(e)!1/3!(d)$);
      \draw[red] (c)--(e);
      \draw[red] (d)--(e);
      \draw[green]($(O)!1/2!(e)$)--($(O)!1/2!(a)$);
      \node [black,below] at ($(d)+(1/3,-1/3)$){$Case$ $1$};
\end{tikzpicture}     

For the first figure above, if moreover the curve $C_{2}$ degenerates
to $C_{3}+B_{0}$, then the lc model is the same as Case 3 with $K^{2}=5$.
\\

\textbf{Case 4.} When the curve $A_{2}$ degenerates to $A_{0}+C_{3}$
and $B_{2}$ degenerates to $B_{3}+A_{0}$. We first blow up the total
space $\mathcal{Y}$ along the line $A_{0}$ in the central fiber,
then blow up along the curve $\tilde{L}_{P}$, which is the proper
transform of $L_{P}$. The central fiber $\Sigma_{0}$ becomes $Bl_{4}\mathbb{P}^{2}\cup\mathbb{F}_{1}$
and the 3-fold is the minimal model of the degenerate family. This
case could be obtained from Case 6 for $K^{2}=6$ in \cite{AP09}
with $A_{1},B_{1},C_{1}$ meeting at a point $P$.

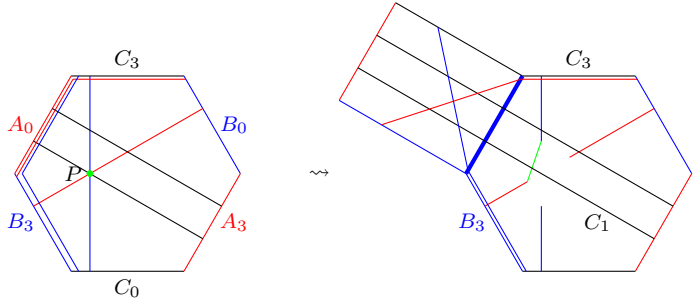
\begin{center}
\begin{tikzpicture}[xscale=1.5,yscale=1.5][font=\tiny]
      \coordinate (a) at (1,3^0.5/2);
      \coordinate (b) at (1/2,3^0.5);
      \coordinate (c) at (-1/2,3^0.5);
      \coordinate (d) at (-1,3^0.5/2);
      \coordinate (e) at (-1/2,0);
      \coordinate (f) at (1/2,0);
      \draw[blue](a)--(b);
      \draw[black](b)--(c);
      \draw[red](c)--(d);
      \draw[blue](d)--(e);
      \draw[black](e)--(f);
      \draw[red](f)--(a);
      \coordinate (g) at ($(a)!1/3!(b)$);
      \coordinate (h) at ($(a)!2/3!(b)$);
      \coordinate (i) at ($(b)!1/3!(c)$);
      \coordinate (j) at (-1/3,3^0.5);
      \coordinate (k) at ($(c)!1/3!(d)$);
      \coordinate (l) at ($(c)!2/3!(d)$);
      \coordinate (m) at ($(d)!1/3!(e)$);
      \coordinate (n) at ($(d)!2/3!(e)$);
      \coordinate (o) at (-1/3,0);
      \coordinate (p) at ($(e)!2/3!(f)$);
      \coordinate (q) at ($(m)!1/4!(h)$);
      \coordinate (t) at ($(m)!1/2!(h)$);
      \coordinate (r) at ($(j)!1/3!(o)$);
      \coordinate (s) at ($(j)!2/3!(o)$);
      \coordinate (u) at ($(a)!1/3!(f)$);
      \coordinate (v) at ($(a)!2/3!(f)$);
      \draw[blue,name path=B_1] (j)--(o);
      \draw[red,name path=A_1](m)--(h);
      \draw[black,name path=C_1](l)--(v);
      \draw[black](k)--(u);
      \node[above,black]at($(b)!1/2!(c)$){$C_3$};
      \node[below,black]at($(e)!1/2!(f)$){$C_0$};
      \node[right,blue]at($(a)!1/2!(b)$){$B_0$};
      \node[left,blue]at($(d)!1/2!(e)$){$B_3$};
      \node[left,red]at($(c)!1/2!(d)$){$A_0$};
      \node[right,red]at($(a)!1/2!(f)$){$A_3$};
      \draw[red]($(d)!1/30!(e)$)--($(c)!1/60!(f)$);
      \draw[red]($(c)!1/60!(f)$)--($(b)!1/30!(a)$);
      \draw[blue]($(e)!2/30!(f)$)--($(d)!2/60!(a)$);
      \draw[blue]($(d)!2/60!(a)$)--($(c)!2/30!(b)$);
      \path [name intersections={of=A_1 and B_1}]; 
      \coordinate [label=left:$P$] (P) at (intersection-1);
      \path[fill=green](P)circle[radius=0.03];
\node at (1.7,3^0.5/2){$\rightsquigarrow$};
      \coordinate (a') at (1+4,3^0.5/2);
      \coordinate (b') at (1/2+4,3^0.5);
      \coordinate (c') at (-1/2+4,3^0.5);
      \coordinate (d') at (-1+4,3^0.5/2);
      \coordinate (e') at (-1/2+4,0);
      \coordinate (f') at (1/2+4,0);
      \draw[blue](a')--(b');
      \draw[black](b')--(c');
      \draw[blue,ultra thick](c')--(d');
      \draw[blue](d')--(e');
      \draw[black](e')--(f');
      \draw[red](f')--(a');
      \coordinate (g') at ($(a')!1/3!(b')$);
      \coordinate (h') at ($(a')!2/3!(b')$);
      \coordinate (i') at ($(b')!1/3!(c')$);
      \coordinate (j') at (-1/3+4,3^0.5);
      \coordinate (k') at ($(c')!1/3!(d')$);
      \coordinate (l') at ($(c')!2/3!(d')$);
      \coordinate (m') at ($(d')!1/3!(e')$);
      \coordinate (n') at ($(d')!2/3!(e')$);
      \coordinate (o') at (-1/3+4,0);
      \coordinate (p') at ($(e')!2/3!(f')$);
      \coordinate (q') at ($(m')!1/4!(h')$);
      \coordinate (t') at ($(m')!1/2!(h')$);
      \coordinate (r') at ($(j')!1/3!(o')$);
      \coordinate (s') at ($(j')!2/3!(o')$);
      \coordinate (u') at ($(a')!1/3!(f')$);
      \coordinate (v') at ($(a')!2/3!(f')$);
      \draw[blue] (j')--(r');
      \draw[blue](s')--(o');
      \draw[red](m')--(q');
      \draw[red](t')--(h');
      \draw[green](q')--(r');
      \draw[black](k')--(u');
      \draw[black](l')--(v');
      \draw[red]($(c')!1/30!(d')$)--($(b')!1/30!(a')$);
      \draw[blue]($(e')!1/30!(f')$)--($(d')!1/30!(c')$);
      \node[above,black]at($(b')!1/2!(c')$){$C_3$};
      \node[left,blue]at($(d')!1/2!(e')$){$B_3$};
      \node[below,black]at($(l')!2/3!(v')$){$C_1$};

      \coordinate (w') at ($(a')!1+3/4!(c')$);
      \coordinate (x') at ($(f')!1+3/4!(d')$);
      \draw[black](c')--(w');
      \draw[blue](d')--(x');
      \draw[red](w')--(x');
      \draw[black](k')--($(u')!1+3/4!(k')$);
      \draw[black](l')--($(v')!1+3/4!(l')$);
      \draw [red]($(c')!1/30!(d')$)--($(d')!2/3!(x')$);
      \draw[blue]($(d')!1/30!(c')$)--($(c')!2/3!(w')$);

\end{tikzpicture}     \\

\par\end{center}

We get the lc model and call it Case 4, by contracting $C_{3}$ and
$B_{3}$ in the component $Bl_{4}\mathbb{P}^{2}$. In the component
$\mathbb{F}_{1}$, the curve $A_{0}$ is the $(-1)$-curve $s$, curves
$B_{3},C_{1},C_{2},C_{3}$ are fibers $f$, and curves $A_{2},B_{2}$
are sections of the numerical type $s+f$. In the component $\mathbb{F}_{0}$,
the double locus is the diagonal $s+f$ and all of the other curves
are fibers.

\begin{center}
\begin{tikzpicture}[xscale=2,yscale=2][font=\tiny]
      \coordinate (a) at (0,0);
      \coordinate (b) at (0,1);
      \coordinate (c) at (-1,1);
      \coordinate (d) at (-1,0);
      \draw [name path=triangle,blue,ultra thick] (a)--(b);
      \draw [name path=C_3,black] (b)--(c);
      \draw [name path=A_0,red] (c)--(d); 
      \draw [name path=B_3,blue] (d)--(a);
      \draw [name path=C_2,black] ($(a)!1/3!(b)$)--($(d)!1/3!(c)$);
      \draw [name path=C_1,black] ($(a)!2/3!(b)$)--($(d)!2/3!(c)$);
      \draw [name path=B_2,blue] (a)--($(b)!2/3!(c)$);
      \draw [name path=A_2,red] (b)--($(a)!2/3!(d)$);
      \node[right] at (a){1};
      \node[right] at (b){4};
      \node[right] at ($(a)!1/3!(b)$){2};
      \node[right] at ($(a)!2/3!(b)$){3};
      \node[red,left] at ($(c)!1/2!(d)$){$A_0$};
      \node[blue,below] at ($(a)!5/6!(d)$){$B_3$};     
      \node[black,above] at ($(b)!5/6!(c)$){$C_3$};
      \node[black] at ($(a)!1/2!(d)+(0,-1/4)$){$\mathbb{F}_{1}$};
      \node[black,above] at ($(a)!2/3!(c)$){$C_2$};
      
      \coordinate (e) at (2,0);\node [below] at (e){1};
      \coordinate (f) at (2,1);\node [above] at (f){4};
      \coordinate (g) at (2-3^0.5/2,1/2);
      \coordinate (h) at (2+3^0.5/2,1/2);
      \draw [name path=triangle,ultra thick,blue] (e)--(f);
      \draw [name path=B_1,blue] (g)--(f); 
      \draw [name path=B_0,blue] (f)--(h);    
      \draw [name path=C_0,black] (g)--(e);
      \draw [name path=A_1,red] (e)--(h);
      \draw [name path=E,green] ($(g)!1/3!(f)$)--($(e)!1/3!(h)$);
      \draw [name path=A_3,red] ($(e)!1/3!(g)$)--($(h)!1/3!(f)$);
      \draw [name path=E,green]($(f)!1/3!(g)$)--($(h)!1/3!(e)$);
      \draw [name path=C_2,black]($(g)!2/3!(f)$)--($(e)!2/3!(h)$);
      \node [red,right] at ($(e)!1/2!(h)$){$A_1$};
      \node [blue,right] at ($(h)!1/3!(f)$){$B_0$};
      \node [blue,left] at ($(f)!1/2!(g)$){$B_1$};
      \node [black,left] at ($(g)!2/3!(e)$){$C_0$};
      \node[black] at ($(e)+(0,-1/4)$){$\mathbb{P}^{1}\times\mathbb{P}^{1}$};
\node[black] at (3/4,1/2){$\cup$};
\node at ($(3/4,1/2)+(0,-1)$){$Case$ $4$};
\end{tikzpicture}
\par\end{center}

There exists some further degenerations of Case 4. Take a one-parameter
family with general fibers $\mathbb{F}_{1}\cup\mathbb{F}_{0}$ we
obtained above. In the central fiber, the curve $C_{2}$ coincides
with $C_{3}$, $A_{2}$ degenerates to $A_{0}+C_{3}$ in the component
$\mathbb{F}_{1}$ and the curve $A_{3}$ coincide with $B_{1}$ in
the component $\mathbb{P}^{1}\times\mathbb{P}^{1}$. The total space
of the one-parameter family is a union of two nonsingular three dimensional
spaces $\mathbb{A}^{1}\times\mathbb{F}_{1}$ and $\mathbb{A}^{1}\times\mathbb{F}_{0}$.
Blowing up the total space along a line in the central fiber is the
same as blowing up the line in each three dimensional space first
and then gluing the two resulting surfaces together in the central
fiber. Now let's see the degenerations in each component. \\

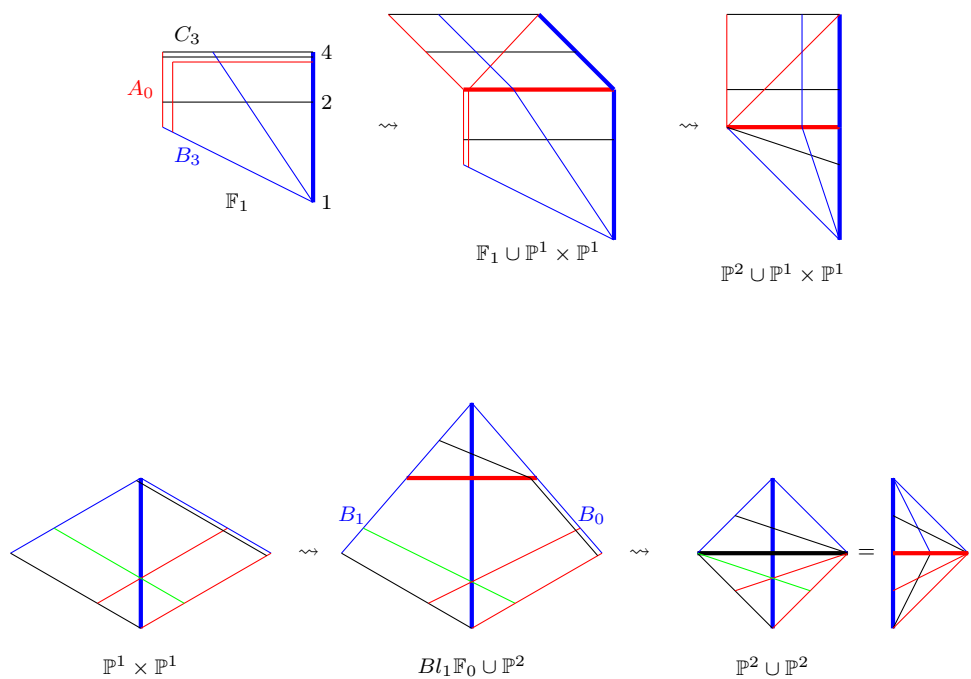
\begin{center}
\begin{tikzpicture}[xscale=2,yscale=2][font=\tiny]
      \coordinate (a) at (0,0);
      \coordinate (b) at (0,1);
      \coordinate (c) at (-1,1);
      \coordinate (d) at (-1,1/2);
      \draw [name path=triangle,blue,ultra thick] (a)--(b);
      \draw [name path=C_3,black] (b)--(c);
      \draw [name path=A_0,red] (c)--(d); 
      \draw [name path=B_3,blue] (d)--(a);
      \draw [name path=C_2,black] ($(b)!2/60!(a)$)--($(c)!2/30!(d)$);
      \draw [name path=C_1,black] ($(b)!1/3!(a)$)--($(c)!2/3!(d)$);
      \draw [name path=B_2,blue] (a)--($(b)!2/3!(c)$);
      \draw [name path=A_2,red] ($(b)!4/60!(a)$)--($(c)!4/60!(a)$);
      \draw [name path=A_2',red] ($(c)!4/60!(a)$)--($(d)!4/60!(a)$);
      \node[right] at (a){1};
      \node[right] at (b){4};
      \node[right] at ($(a)!2/3!(b)$){2};
      \node[red,left] at ($(c)!1/2!(d)$){$A_0$};
      \node[blue,below] at ($(a)!5/6!(d)$){$B_3$};     
      \node[black,above] at ($(b)!5/6!(c)$){$C_3$};
      \node[black] at ($(a)!1/2!(d)+(0,-1/4)$){$\mathbb{F}_{1}$};
\node at (0.5,0.5){$\rightsquigarrow$};
      \coordinate (a) at (0+2,0-1/4);
      \coordinate (b) at (0+2,1-1/4);
      \coordinate (c) at (-1+2,1-1/4);
      \coordinate (d) at (-1+2,1/2-1/4);
      \coordinate (e) at (-1/2+2,3/2-1/4);
      \coordinate (f) at (-3/2+2,3/2-1/4); 
      \draw [name path=triangle,blue,ultra thick] (a)--(b);
      \draw [name path=C_3,red,ultra thick] (b)--(c);
      \draw [name path=A_0,red] (c)--(d); 
      \draw [name path=B_3,blue] (d)--(a);
      \draw [name path=C_1,black] ($(b)!1/3!(a)$)--($(c)!2/3!(d)$);
      \draw [name path=B_2,blue] (a)--($(b)!2/3!(c)$);
      \draw [name path=A_2',red] ($(c)!2/60!(b)$)--($(d)!2/60!(a)$);
      \node[black] at ($(a)!1/2!(d)+(0,-1/3)$){$\mathbb{F}_{1}\cup\mathbb{P}^{1}\times\mathbb{P}^{1}$};
      \draw [name path=triangle,blue,ultra thick] (b)--(e);
      \draw [name path=C_3,black] (e)--(f);
      \draw [name path=C_2,black] ($(b)!1/2!(e)$)--($(c)!1/2!(f)$);
      \draw [name path=A_0,red] (c)--(f);
      \draw [name path=A_2,red] (e)--($(c)!2/60!(b)$);
      \draw [name path=B_2,blue] ($(b)!2/3!(c)$)--($(e)!2/3!(f)$);
\node at (0.5+2,0.5){$\rightsquigarrow$};
      \coordinate (a) at (0+3.5,0-1/4);
      \coordinate (b) at (0+3.5,3/4-1/4);
      \coordinate (c) at (-3/4+3.5,3/4-1/4);
      \coordinate (e) at (0+3.5,3/2-1/4);
      \coordinate (f) at (-3/4+3.5,3/2-1/4); 
      \draw [name path=triangle,blue,ultra thick] (a)--(b);
      \draw [name path=C_3,red,ultra thick] (b)--(c);
      \draw [name path=B_3,blue] (c)--(a);
      \draw [name path=C_1,black] ($(b)!1/3!(a)$)--(c);
      \draw [name path=B_2,blue] (a)--($(b)!1/3!(c)$);
      \node[black] at ($(a)!1/2!(c)+(0,-3/5)$){$\mathbb{P}^{2}\cup\mathbb{P}^{1}\times\mathbb{P}^{1}$};
      \draw [name path=triangle,blue,ultra thick] (b)--(e);
      \draw [name path=C_3,black] (e)--(f);
      \draw [name path=C_2,black] ($(b)!1/3!(e)$)--($(c)!1/3!(f)$);
      \draw [name path=A_0,red] (c)--(f);
      \draw [name path=A_2,red] (e)--(c);
      \draw [name path=B_2,blue] ($(b)!1/3!(c)$)--($(e)!1/3!(f)$);  
\end{tikzpicture}\\

\par\end{center}

\begin{center}
\begin{tikzpicture}[xscale=2,yscale=2][font=\tiny]
      \coordinate (a) at (0,0);
      \coordinate (b) at (3^0.5/2,1/2);
      \coordinate (c) at (0,1);
      \coordinate (d) at (-3^0.5/2,1/2);
      \draw [name path=triangle,ultra thick,blue] (a)--(c);
      \draw [name path=B_1,blue] (d)--(c); 
      \draw [name path=B_0,blue] (c)--(b);    
      \draw [name path=C_0,black] (d)--(a);
      \draw [name path=A_1,red] (a)--(b);
      \draw [name path=A_3,red] ($(b)!1/3!(c)$)--($(a)!1/3!(d)$);
      \draw [name path=E,green] ($(d)!1/3!(c)$)--($(a)!1/3!(b)$);
      \draw [name path=C_2,black]($(c)!1/30!(d)$)--($(b)!1/30!(a)$);
        \node[black] at ($(a)+(0,-1/4)$){$\mathbb{P}^{1}\times\mathbb{P}^{1}$};
      \node at ($(b)+(1/4,0)$){$\rightsquigarrow$};
      \coordinate (a) at (0+2.2,0);
      \coordinate (b) at (3^0.5/2+2.2,1/2);
      \coordinate (c) at (0+2.2,1);
      \coordinate (d) at (-3^0.5/2+2.2,1/2);
      \coordinate (e) at (3^0.5/2/2+2.2,1);
      \coordinate (f) at (-3^0.5/2/2+2.2,1);
      \coordinate (g) at (0+2.2,3/2);
      \draw [name path=triangle,ultra thick,blue] (a)--(g);
      \draw [name path=B_1,blue] (d)--(f); 
      \draw [name path=B_0,blue] (e)--(b);    
      \draw [name path=C_0,black] (d)--(a);
      \draw [name path=E,green] ($(d)!1/3!(f)$)--($(a)!1/3!(b)$);
      \draw [name path=triangle',red, ultra thick] (e)--(f);
      \draw [name path=A_1,red] (a)--(b);
      \draw [name path=A_3,red] ($(a)!1/3!(d)$)--($(b)!1/3!(e)$);
      \draw [name path=C_2,black]($(b)!1/30!(a)$)--($(e)!1/20!(f)$);
       \node[black] at ($(a)+(0,-1/4)$){$Bl_{1}\mathbb{F}_{0}\cup\mathbb{P}^{2}$};
      \draw [blue] (e)--(g);
      \draw [blue] (f)--(g);
      \draw [black,name path=C_2] ($(e)!1/20!(f)$)--($(f)!1/2!(g)$);
      \node[right, blue] at ($(b)!1/2!(e)$){$B_0$};
      \node[left, blue] at ($(d)!1/2!(f)$){$B_1$};
    \node at ($(b)+(1/4,0)$){$\rightsquigarrow$};
      \coordinate (a) at (0+4.2,0);
      \coordinate (c) at (0+4.2,1/2);
      \coordinate (e) at (1/2+4.2,1/2);
      \coordinate (f) at (-1/2+4.2,1/2);
      \coordinate (g) at (0+4.2,1);
      \draw [name path=triangle,ultra thick,blue] (a)--(g);
      \draw [name path=C_0,black] (f)--(a);
      \draw [name path=triangle',black, ultra thick] (e)--(f);
      \draw [name path=A_1,red] (a)--(e);
      \draw [name path=A_3,red] ($(f)!1/2!(a)$)--(e);
      \draw [name path=E,green] ($(e)!1/2!(a)$)--(f); 
             \node[black] at ($(a)+(0,-1/4)$){$\mathbb{P}^{2}\cup\mathbb{P}^{2}$};
      \draw [blue] (e)--(g);
      \draw [blue] (f)--(g);
      \draw [black] (e)--($(f)!1/2!(g)$);
  \node at ($(e)+(1/8,0)$){$=$};
      \coordinate (a) at (0+5,0);
      \coordinate (b) at (1/2+5,1/2);
      \coordinate (c) at (0+5,1);
      \coordinate (d) at ($(a)!1/2!(c)$);
      \draw [blue,ultra thick] (a)--(c);
      \draw [blue] (b)--(c);
      \draw [blue] (c)--($(b)!1/2!(d)$);
      \draw [black] (b)--($(c)!1/2!(d)$);
      \draw [black] (a)--($(b)!1/2!(d)$);
      \draw [red,ultra thick] (b)--(d);
      \draw [red] (a)--(b);
      \draw [red] (b)--($(a)!1/2!(d)$);
     
\end{tikzpicture}\\

\par\end{center}

In the central fiber, the surface $\mathbb{F}_{1}$ splits into $\mathbb{P}^{2}\cup\mathbb{P}^{1}\times\mathbb{P}^{1}$
and $\mathbb{P}^{1}\times\mathbb{P}^{1}$ becomes $\mathbb{P}^{2}\cup\mathbb{P}^{2}$.
Gluing the two resulting surfaces together, we obtain a further degeneration
as in the first figure below. Another possible degeneration is that
$C_{2}$ moves to $B_{3}$, and $B_{2}$ degenerates to $A_{0}+B_{3}$
(the central fiber of lc model is the second figure below). Since
both the first and second figures contain a component $\mathbb{P}^{1}\times\mathbb{P}^{1}$,
the diagonals can again be degenerated to the section $s+f$, and
the limits are the two remaining figures below, each of which are
just the previous cases 2 and 3. \\

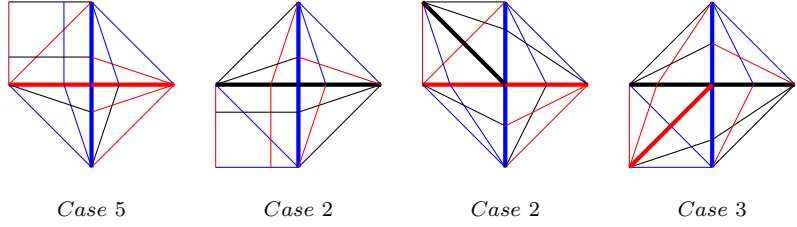
\begin{center}
\begin{tikzpicture}[xscale=1.1,yscale=1.1][font=\tiny]
      \coordinate (a) at (1,0);
      \coordinate (b) at (0,1);
      \coordinate (c) at (-1,1);
      \coordinate (d) at (-1,0);
      \coordinate (e) at (0,-1);
      \coordinate (O) at (0,0);
      \draw [red] (c)--(d); 
      \draw [red, ultra thick] (O)--(d);
      \draw [blue] (a)--(b);
      \draw [blue] (d)--(e);
      \draw [blue, ultra thick] (b)--(e);
      \draw [black] (b)--(c);
      \draw [red] (e)--(a);
      \draw [red, ultra thick] (O)--(a);
      \draw [red] ($(O)!1/3!(b)$)--(a);
      \draw [black] ($(O)!1/3!(a)$)--(e);
      \draw [red] (b)--(d);
      \draw [blue] ($(O)!1/3!(a)$)--(b);
      \draw [blue] ($(O)!1/3!(d)$)--(e);
      \draw [blue] ($(O)!1/3!(d)$)--($(b)!1/3!(c)$);
      \draw [red] ($(O)!1/3!(e)$)--(a);
      \draw [black] ($(O)!1/3!(e)$)--(d);
      \draw [black] ($(O)!1/3!(b)$)--($(d)!1/3!(c)$);
      \node at ($(e)+(0,-0.5)$){$Case$ $5$};
      \coordinate (a) at (1+2.5,0);
      \coordinate (e) at (0+2.5,1);
      \coordinate (c) at (-1+2.5,-1);
      \coordinate (d) at (-1+2.5,0);
      \coordinate (b) at (0+2.5,-1);
      \coordinate (O) at (0+2.5,0);
      \draw [red] (c)--(d); 
      \draw [black, ultra thick] (O)--(d);
      \draw [black] (a)--(b);
      \draw [black] (d)--(e);
      \draw [blue, ultra thick] (b)--(e);
      \draw [blue] (b)--(c);
      \draw [blue] (e)--(a);
      \draw [black, ultra thick] (O)--(a);
      \draw [black] ($(O)!1/3!(b)$)--(a);
      \draw [blue] ($(O)!1/3!(a)$)--(e);
      \draw [blue] (b)--(d);
      \draw [red] ($(O)!1/3!(a)$)--(b);
      \draw [red] ($(O)!1/3!(d)$)--(e);
      \draw [red] ($(O)!1/3!(d)$)--($(b)!1/3!(c)$);
      \draw [red] ($(O)!1/3!(e)$)--(a);
      \draw [black] ($(O)!1/3!(e)$)--(d);
      \draw [black] ($(O)!1/3!(b)$)--($(d)!1/3!(c)$);
     \node at ($(b)+(0,-0.5)$){$Case$ $2$};

      \coordinate (a) at (1+5,0);
      \coordinate (b) at (0+5,1);
      \coordinate (c) at (-1+5,1);
      \coordinate (d) at (-1+5,0);
      \coordinate (e) at (0+5,-1);
      \coordinate (O) at (0+5,0);
      \draw [red] (c)--(d); 
      \draw [red, ultra thick] (O)--(d);
      \draw [blue] (a)--(b);
      \draw [blue] (d)--(e);
      \draw [blue, ultra thick] (b)--(e);
      \draw [black] (b)--(c);
      \draw [red] (e)--(a);
      \draw [red, ultra thick] (O)--(a);
      \draw [black] ($(O)!2/3!(b)$)--(a);
      \draw [black] ($(O)!1/2!(a)$)--(e);
      \draw [red] (b)--(d);
      \draw [blue] ($(O)!1/2!(a)$)--(b);
      \draw [blue] ($(O)!2/3!(d)$)--(e);
      \draw [blue] ($(O)!2/3!(d)$)--(c);
      \draw [red] ($(O)!1/2!(e)$)--(a);
      \draw [black] ($(O)!1/2!(e)$)--(d);
      \draw [black] ($(O)!2/3!(b)$)--(c);
      \draw [black,ultra thick] (O)--(c);  
      \node at ($(e)+(0,-0.5)$){$Case$ $2$};

      \coordinate (a) at (1+7.5,0);
      \coordinate (e) at (0+7.5,1);
      \coordinate (c) at (-1+7.5,-1);
      \coordinate (d) at (-1+7.5,0);
      \coordinate (b) at (0+7.5,-1);
      \coordinate (O) at (0+7.5,0);
      \draw [red] (c)--(d); 
      \draw [black, ultra thick] (O)--(d);
      \draw [black] (a)--(b);
      \draw [black] (d)--(e);
      \draw [blue, ultra thick] (b)--(e);
      \draw [blue] (b)--(c);
      \draw [blue] (e)--(a);
      \draw [black, ultra thick] (O)--(a);
      \draw [black] ($(O)!2/3!(b)$)--(a);
      \draw [blue] ($(O)!1/2!(a)$)--(e);
      \draw [blue] (b)--(d);
      \draw [red] ($(O)!1/2!(a)$)--(b);
      \draw [red] ($(O)!2/3!(d)$)--(e);
      \draw [red] ($(O)!2/3!(d)$)--(c);
      \draw [red] ($(O)!1/2!(e)$)--(a);
      \draw [black] ($(O)!1/2!(e)$)--(d);
      \draw [black] ($(O)!2/3!(b)$)--(c);
      \draw [red, ultra thick] (O)--(c);
      \node at ($(b)+(0,-0.5)$){$Case$ $3$};
\end{tikzpicture}  
\par\end{center}

\textbf{Case 6. }When the five lines $A_{1},A_{2},B_{1},B_{2},C_{1},C_{2}$
meet at the point $P$. We first blow up the total space along the
curve $L_{P}$, then blow up the point $P$ in the central fiber,
which is the intersection of $A_{1},B_{1},C_{1}$ in the exceptional
divisor $\mathbb{P}^{2}$ of the blowup. The resulting central fiber
contains two components $Bl_{4}\mathbb{P}^{2}\cup\mathbb{F}_{1}$,
which is the central fiber of the minimal model.\\

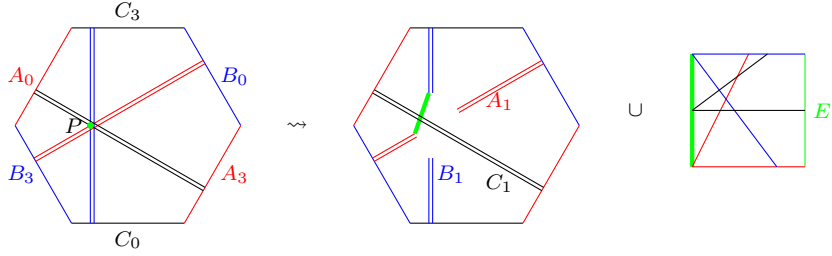
\begin{center}
\begin{tikzpicture}[xscale=1.5,yscale=1.5][font=\tiny]
      \coordinate (f) at (1/2,0);
      \coordinate (a) at (1,3^0.5/2);
      \coordinate (b) at (1/2,3^0.5);
      \coordinate (c) at (-1/2,3^0.5);
      \coordinate (d) at (-1,3^0.5/2);
      \coordinate (e) at (-1/2,0);
      \draw[blue](a)--(b);
      \draw[black](b)--(c);
      \draw[red](c)--(d);
      \draw[blue](d)--(e);
      \draw[black](e)--(f);
      \draw[red](f)--(a);
      \coordinate (h) at ($(a)!2/3!(b)$);
      \coordinate (j) at (-1/3,3^0.5);
      \coordinate (l) at ($(c)!2/3!(d)$);
      \coordinate (m) at ($(d)!1/3!(e)$);
      \coordinate (o) at (-1/3,0);
      \coordinate (v) at ($(a)!2/3!(f)$);
      \draw[blue,name path=B_1] (j)--(o);
      \draw[blue,name path=B_2] ($(j)+(1/30,0)$)--($(o)+(1/30,0)$);
      \draw[red,name path=A_1] (h)--(m);
      \draw[red,name path=A_2] ($(e)!2/3-1/30!(d)$)--($(a)!2/3-1/30!(b)$);
      \draw[black,name path=C_1](l)--(v);
      \draw[black,name path=C_2]($(c)!2/3-1/30!(d)$)--($(a)!2/3-1/30!(f)$);
      \node[above,black]at($(b)!1/2!(c)$){$C_3$};
      \node[below,black]at($(e)!1/2!(f)$){$C_0$};
      \node[right,blue]at($(a)!1/2!(b)$){$B_0$};
      \node[left,blue]at($(d)!1/2!(e)$){$B_3$};
      \node[left,red]at($(c)!1/2!(d)$){$A_0$};
      \node[right,red]at($(a)!1/2!(f)$){$A_3$};
      \path [name intersections={of=A_1 and B_1}]; 
      \coordinate [label=left:$P$] (P) at (intersection-1);
      \path[fill=green](P)circle[radius=0.03];
\node at ($(a)+(0.5,0)$){$\rightsquigarrow$};
      \coordinate (f) at (1/2+3,0);
      \coordinate (a) at (1+3,3^0.5/2);
      \coordinate (b) at (1/2+3,3^0.5);
      \coordinate (c) at (-1/2+3,3^0.5);
      \coordinate (d) at (-1+3,3^0.5/2);
      \coordinate (e) at (-1/2+3,0);
      \draw[blue](a)--(b);
      \draw[black](b)--(c);
      \draw[red](c)--(d);
      \draw[blue](d)--(e);
      \draw[black](e)--(f);
      \draw[red](f)--(a);
      \coordinate (h) at ($(a)!2/3!(b)$);
      \coordinate (j) at (-1/3+3,3^0.5);
      \coordinate (l) at ($(c)!2/3!(d)$);
      \coordinate (m) at ($(d)!1/3!(e)$);
      \coordinate (o) at (-1/3+3,0);
      \coordinate (q) at ($(m)!1/4!(h)$);
      \coordinate (t) at ($(m)!1/2!(h)$);
      \coordinate (r) at ($(j)!1/3!(o)$);
      \coordinate (s) at ($(j)!2/3!(o)$);
      \coordinate (v) at ($(a)!2/3!(f)$);
      \draw[blue] (j)--(r);
      \draw[blue](s)--(o);
      \draw[blue]($(s)+(1/30,0)$)--($(o)+(1/30,0)$);
      \draw[blue]($(j)+(1/30,0)$)--($(r)+(1/30,0)$);
      \draw[red](m)--(q);
      \draw[red](t)--(h);
      \draw[red]($(t)+(1/60,-3^0.5/60)$)--($(a)!2/3-1/30!(b)$);
      \draw[red]($(q)+(1/60,-3^0.5/60)$)--($(e)!2/3-1/30!(d)$);
      \draw[green,ultra thick](q)--(r);
      \draw[black](l)--(v);
      \draw[black]($(c)!2/3-1/30!(d)$)--($(a)!2/3-1/30!(f)$);
      \node[below,black]at($(v)!1/4!(l)$){$C_1$};
      \node[right,blue]at($(j)!3/4!(o)$){$B_1$};
      \node[below,red]at($(h)!1/4!(m)$){$A_1$};
      \coordinate (a) at (5,0+1/2);
      \coordinate (b) at (6,0+1/2);
      \coordinate (c) at (6,1+1/2);
      \coordinate (d) at (5,1+1/2);
      \draw[green,ultra thick] (a)--(d);
      \draw[green] (b)--(c);
      \draw[red] (a)--(b);
      \draw[red] (a)--($(d)!1/2!(c)$);
      \draw[blue] (c)--(d);
      \draw[blue] (d)--($(a)!3/4!(b)$);
      \draw[black] ($(a)!1/2!(d)$)--($(b)!1/2!(c)$);
      \draw[black] ($(a)!1/2!(d)$)--($(d)!2/3!(c)$);
      \node[green,right] at ($(b)!1/2!(c)$){$E$};
      \node [black] at (4.5,1){$\cup$};
\end{tikzpicture}    \\

\par\end{center}

Running the minimal model program, we contract $A_{1},B_{1},C_{1}$
and get the log canonical model with the central fiber $\mathbb{P}^{1}\times\mathbb{P}^{1}\cup\mathbb{F}_{1}$,
where $E$ is the (-1)-curve in $\mathbb{F}_{1}$. There is no further
degeneration for this case.

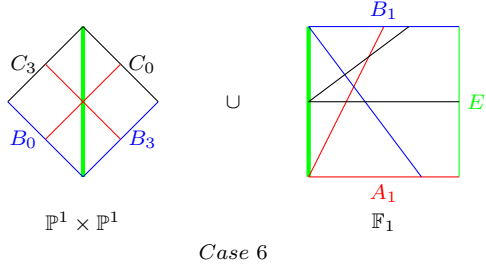
\begin{center}
\begin{tikzpicture}[xscale=2,yscale=2][font=\tiny]
      \coordinate (O) at (0,0);
      \coordinate (a) at (1/2,1/2);
      \coordinate (b) at (0,1);
      \coordinate (c) at (-1/2,1/2);
      \draw [name path=triangle,green, ultra thick] (O)--(b);
      \draw [name path=B_3,blue] (O)--(a);
      \draw [name path=B_0,blue] (O)--(c);
      \draw [name path=C_0,black] (a)--(b);
      \draw [name path=C_3,black] (b)--(c);
      \draw [name path=A_3,red] ($(O)!1/2!(a)$)--($(b)!1/2!(c)$);
      \draw [name path=A_0,red] ($(O)!1/2!(c)$)--($(a)!1/2!(b)$);
      \node[blue,right] at ($(O)!1/2!(a)$){$B_3$};
      \node[blue,left] at ($(O)!1/2!(c)$){$B_0$};
      \node[black,right] at ($(a)!1/2!(b)$){$C_0$};
      \node[black,left] at ($(b)!1/2!(c)$){$C_3$};
      \node at ($(O)+(0,-0.3)$){$\mathbb{P}^{1}\times\mathbb{P}^{1}$};
      \coordinate (a) at (3/2,0);
      \coordinate (b) at (3/2+1,0);
      \coordinate (c) at (3/2+1,1);
      \coordinate (d) at (3/2,1);
      \draw[name path=triangle,green,ultra thick] (a)--(d);
      \draw[name path=E,green] (b)--(c);
      \draw[name path=A_1,red] (a)--(b);
      \draw[name path=A_2,red] (a)--($(d)!1/2!(c)$);
      \draw[name path=B_1,blue] (c)--(d);
      \draw[name path=B_2,blue] (d)--($(a)!3/4!(b)$);
      \draw[name path=C_1,black] ($(a)!1/2!(d)$)--($(b)!1/2!(c)$);
      \draw[name path=C_2,black] ($(a)!1/2!(d)$)--($(d)!2/3!(c)$);
     \node[green,right] at ($(b)!1/2!(c)$){$E$};
     \node[blue,above] at ($(c)!1/2!(d)$){$B_1$};
     \node[red,below] at ($(a)!1/2!(b)$){$A_1$};
      \coordinate (e) at ($(a)!1/2!(b)$);
      \node at ($(e)+(0,-0.3)$){$\mathbb{F}_{1}$};
      \node [black] at (1,1/2){$\cup$};
      \node [black] at ($(1,1/2)+(0,-1)$){$Case$ $6$};
\end{tikzpicture}
\par\end{center}

Case 6 can be obtained from Case 9 with $K^{2}=6$, by taking $A_{1},B_{1},C_{1}$
to have a common intersection. For the above surface, if $B_{2},C_{2}$
in the second component $\mathbb{F}_{1}$ degenerate to $B_{1}+E$,
$C_{1}+E$, then it is the central fiber of lc model of the degeneration
comes from Case 10 with $K^{2}=6$ .

\subsection{Log canonical degenerations.}

Case 1,8 and 5 with $K^{2}=6$ are special. Case 5 with $K^{2}=6$
does not produce any degenerations with $K^{2}=5$. 

For Case 5 with $K^{2}=6$ (the left figure below), there is no corresponding
degeneration with $K^{2}=5$. Since $A_{1}$, $B_{1}$, $C_{1}$ must
intersect, the resulting degeneration has an infinite automorphism
group, and therefore does not correspond to an irreducible component
of a stable pair.

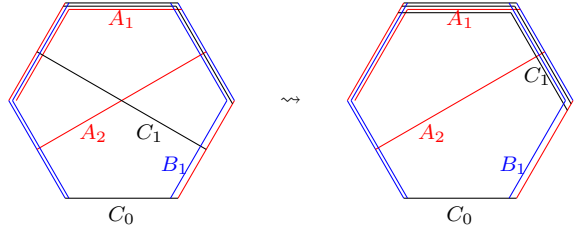
\begin{center}
\begin{tikzpicture}[xscale=1.5,yscale=1.5][font=\tiny]
      \coordinate (z) at (0,0);
      \coordinate (a) at ($(1,3^0.5/2)+(z)$);
      \coordinate (b) at ($(1/2,3^0.5)+(z)$);
      \coordinate (c) at ($(-1/2,3^0.5)+(z)$);
      \coordinate (d) at ($(-1,3^0.5/2)+(z)$);
      \coordinate (e) at ($(-1/2,0)+(z)$);
      \coordinate (f) at ($(1/2,0)+(z)$);
      \draw[blue](a)--(b);
      \draw[black](b)--(c);
      \draw[red](c)--(d);
      \draw[blue](d)--(e);
      \draw[black](e)--(f);
      \draw[red](f)--(a);
      \coordinate (g) at ($(a)!1/3!(b)$);
      \coordinate (h) at ($(a)!2/3!(b)$);
      \coordinate (i) at ($(b)!1/3!(c)$);
      \coordinate (j) at ($(b)!2/3!(c)$);
      \coordinate (k) at ($(c)!1/3!(d)$);
      \coordinate (l) at ($(c)!2/3!(d)$);
      \coordinate (m) at ($(a)!1/2!(b)$);
      \coordinate (n) at ($(d)!1/2!(e)$);
      \coordinate (q) at ($(a)!1/2!(f)$);
      \coordinate (r) at ($(c)!1/2!(d)$);
      \node[below,black]at($(e)!1/2!(f)$){$C_0$};
      \draw[red,name path=A_1]($(d)!2/60!(a)$)--($(c)!2/60!(f)$);
      \draw[red,name path=A_1]($(c)!2/60!(f)$)--($(b)!2/30!(a)$);     
      \draw[red,name path=A_2]($(d)!1/2!(e)$)--($(b)!1/2!(a)$);
      \draw[black,name path=C_1]($(d)!1/2!(c)$)--($(a)!1/2!(f)$);
      \draw[black,name path=C_2]($(c)!1/30!(d)$)--($(b)!1/60!(e)$);
      \draw[black,name path=C_2]($(b)!1/60!(e)$)--($(a)!1/30!(f)$);     
      \draw[blue,name path=B_2]($(e)!1/30!(f)$)--($(d)!1/60!(a)$);
      \draw[blue,name path=B_2]($(d)!1/60!(a)$)--($(c)!1/30!(b)$);
      \draw[blue,name path=B_1]($(b)!2/30!(c)$)--($(a)!2/60!(d)$);
      \draw[blue,name path=B_1]($(a)!2/60!(d)$)--($(f)!2/30!(e)$);      
      \node[black,below]at($(q)!1/3!(r)$){$C_1$};
      \node[red,below]at($(n)!1/3!(m)$){$A_2$};
      \node[blue,left]at($(f)!1/3!(a)$){$B_1$};
      \node[red,below]at($(b)!1/2!(c)$){$A_1$};
\node at ($(a)+(0.5,0)$){$\rightsquigarrow$};
      \coordinate (z) at (3,0);
      \coordinate (a) at ($(1,3^0.5/2)+(z)$);
      \coordinate (b) at ($(1/2,3^0.5)+(z)$);
      \coordinate (c) at ($(-1/2,3^0.5)+(z)$);
      \coordinate (d) at ($(-1,3^0.5/2)+(z)$);
      \coordinate (e) at ($(-1/2,0)+(z)$);
      \coordinate (f) at ($(1/2,0)+(z)$);
      \draw[blue](a)--(b);
      \draw[black](b)--(c);
      \draw[red](c)--(d);
      \draw[blue](d)--(e);
      \draw[black](e)--(f);
      \draw[red](f)--(a);
      \coordinate (g) at ($(a)!1/3!(b)$);
      \coordinate (h) at ($(a)!2/3!(b)$);
      \coordinate (i) at ($(b)!1/3!(c)$);
      \coordinate (j) at ($(b)!2/3!(c)$);
      \coordinate (k) at ($(c)!1/3!(d)$);
      \coordinate (l) at ($(c)!2/3!(d)$);
      \coordinate (m) at ($(a)!1/2!(b)$);
      \coordinate (n) at ($(d)!1/2!(e)$);
      \coordinate (q) at ($(a)!1/2!(f)$);
      \coordinate (r) at ($(c)!1/2!(d)$);
      \node[below,black]at($(e)!1/2!(f)$){$C_0$};
      \draw[red,name path=A_1]($(d)!2/60!(a)$)--($(c)!2/60!(f)$);
      \draw[red,name path=A_1]($(c)!2/60!(f)$)--($(b)!2/30!(a)$);     
      \draw[red,name path=A_2]($(d)!1/2!(e)$)--($(b)!1/2!(a)$);
      \draw[black,name path=C_1]($(c)!3/30!(d)$)--($(b)!3/60!(e)$);
      \draw[black,name path=C_1]($(b)!3/60!(e)$)--($(a)!3/30!(f)$);
      \draw[black,name path=C_2]($(c)!1/30!(d)$)--($(b)!1/60!(e)$);
      \draw[black,name path=C_2]($(b)!1/60!(e)$)--($(a)!1/30!(f)$);     
      \draw[blue,name path=B_2]($(e)!1/30!(f)$)--($(d)!1/60!(a)$);
      \draw[blue,name path=B_2]($(d)!1/60!(a)$)--($(c)!1/30!(b)$);
      \draw[blue,name path=B_1]($(b)!2/30!(c)$)--($(a)!2/60!(d)$);
      \draw[blue,name path=B_1]($(a)!2/60!(d)$)--($(f)!2/30!(e)$);      
      \node[black,left]at($(a)!1/4!(b)$){$C_1$};
      \node[red,below]at($(n)!1/3!(m)$){$A_2$};
      \node[blue,left]at($(f)!1/3!(a)$){$B_1$};
      \node[red,below]at($(b)!1/2!(c)$){$A_1$};
     
\end{tikzpicture}     
\par\end{center}

Case 1 and 8 with $K^{2}=6$ produce degenerations with $K^{2}=5$.
But it is surprising that the lc models of the degenerations are irreducible
and are the same as some lc degenerations. We elaborate on the special
cases 1 and 8 as following. 

We first look at Case 8 with $K^{2}=6$ which is also a degeneration
with $K^{2}=5$. When all of the five lines $A_{1},A_{2},B_{1},B_{2},C_{1}$
meet at a point $P$, we blow up the total space along the curve $L_{P}$,
and the resulting central fiber contains two components $Bl_{4}\mathbb{P}^{2}\cup\mathbb{F}_{1}$.
Running the minimal model program, the whole component $Bl_{1}\mathbb{P}^{2}$
is contracted and the central fiber of the lc model is irreducible,
which is $Bl_{4}\mathbb{P}^{2}$.\\

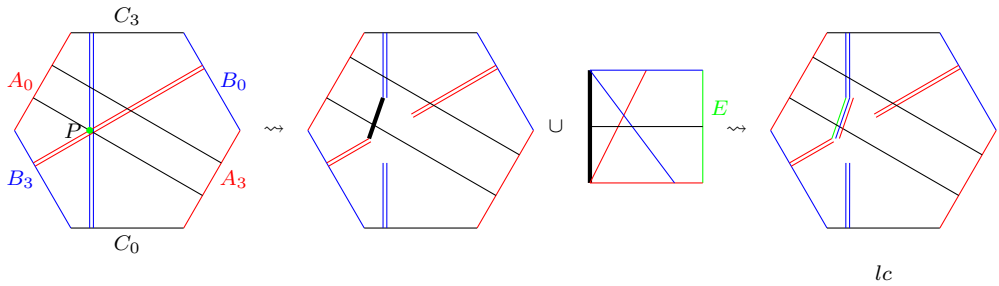
\begin{center}
\begin{tikzpicture}[xscale=1.5,yscale=1.5][font=\tiny]
      \coordinate (f) at (1/2,0);
      \coordinate (a) at (1,3^0.5/2);
      \coordinate (b) at (1/2,3^0.5);
      \coordinate (c) at (-1/2,3^0.5);
      \coordinate (d) at (-1,3^0.5/2);
      \coordinate (e) at (-1/2,0);
      \draw[blue](a)--(b);
      \draw[black](b)--(c);
      \draw[red](c)--(d);
      \draw[blue](d)--(e);
      \draw[black](e)--(f);
      \draw[red](f)--(a);
      \coordinate (h) at ($(a)!2/3!(b)$);
      \coordinate (j) at (-1/3,3^0.5);
      \coordinate (l) at ($(c)!2/3!(d)$);
      \coordinate (m) at ($(d)!1/3!(e)$);
      \coordinate (o) at (-1/3,0);
      \coordinate (v) at ($(a)!2/3!(f)$);
      \draw[blue,name path=B_1] (j)--(o);
      \draw[blue,name path=B_2] ($(j)+(1/30,0)$)--($(o)+(1/30,0)$);
      \draw[red,name path=A_1] (h)--(m);
      \draw[red,name path=A_2] ($(e)!2/3-1/30!(d)$)--($(a)!2/3-1/30!(b)$);
      \draw[black,name path=C_1](l)--(v);
      \draw[black,name path=C_2]($(c)!1/3!(d)$)--($(a)!1/3!(f)$);
      \node[above,black]at($(b)!1/2!(c)$){$C_3$};
      \node[below,black]at($(e)!1/2!(f)$){$C_0$};
      \node[right,blue]at($(a)!1/2!(b)$){$B_0$};
      \node[left,blue]at($(d)!1/2!(e)$){$B_3$};
      \node[left,red]at($(c)!1/2!(d)$){$A_0$};
      \node[right,red]at($(a)!1/2!(f)$){$A_3$};
      \path [name intersections={of=A_1 and B_1}]; 
      \coordinate [label=left:$P$] (P) at (intersection-1);
      \path[fill=green](P)circle[radius=0.03];
\node at ($(a)+(0.3,0)$){$\rightsquigarrow$};
      \coordinate (f) at (1/2+2.6,0);
      \coordinate (a) at (1+2.6,3^0.5/2);
      \coordinate (b) at (1/2+2.6,3^0.5);
      \coordinate (c) at (-1/2+2.6,3^0.5);
      \coordinate (d) at (-1+2.6,3^0.5/2);
      \coordinate (e) at (-1/2+2.6,0);
      \draw[blue](a)--(b);
      \draw[black](b)--(c);
      \draw[red](c)--(d);
      \draw[blue](d)--(e);
      \draw[black](e)--(f);
      \draw[red](f)--(a);
      \coordinate (h) at ($(a)!2/3!(b)$);
      \coordinate (j) at (-1/3+2.6,3^0.5);
      \coordinate (l) at ($(c)!2/3!(d)$);
      \coordinate (m) at ($(d)!1/3!(e)$);
      \coordinate (o) at (-1/3+2.6,0);
      \coordinate (q) at ($(m)!1/4!(h)$);
      \coordinate (t) at ($(m)!1/2!(h)$);
      \coordinate (r) at ($(j)!1/3!(o)$);
      \coordinate (s) at ($(j)!2/3!(o)$);
      \coordinate (v) at ($(a)!2/3!(f)$);
      \draw[blue] (j)--(r);
      \draw[blue](s)--(o);
      \draw[blue]($(s)+(1/30,0)$)--($(o)+(1/30,0)$);
      \draw[blue]($(j)+(1/30,0)$)--($(r)+(1/30,0)$);
      \draw[red](m)--(q);
      \draw[red](t)--(h);
      \draw[red]($(t)+(1/60,-3^0.5/60)$)--($(a)!2/3-1/30!(b)$);
      \draw[red]($(q)+(1/60,-3^0.5/60)$)--($(e)!2/3-1/30!(d)$);
      \draw[black,ultra thick](q)--(r);
      \draw[black](l)--(v);
      \draw[black,name path=C_2]($(c)!1/4!(d)$)--($(a)!1/4!(f)$);
\node [black] at (3.8,1-0.1){$\cup$};
      \coordinate (a) at (4.1,0+1/2-0.1);
      \coordinate (b) at (5.1,0+1/2-0.1);
      \coordinate (c) at (5.1,1+1/2-0.1);
      \coordinate (d) at (4.1,1+1/2-0.1);
      \draw[black,ultra thick] (a)--(d);
      \draw[green] (b)--(c);
      \draw[red] (a)--(b);
      \draw[red] (a)--($(d)!1/2!(c)$);
      \draw[blue] (c)--(d);
      \draw[blue] (d)--($(a)!3/4!(b)$);
      \draw[black] ($(a)!1/2!(d)$)--($(b)!1/2!(c)$);
      \node[green,right] at ($(b)!2/3!(c)$){$E$};

      \coordinate (x) at (6.7,0);
      \coordinate (f) at ($(1/2,0)+(x)$);
      \coordinate (a) at ($(1,3^0.5/2)+(x)$);
      \coordinate (b) at ($(1/2,3^0.5)+(x)$);
      \coordinate (c) at ($(-1/2,3^0.5)+(x)$);
      \coordinate (d) at ($(-1,3^0.5/2)+(x)$);
      \coordinate (e) at ($(-1/2,0)+(x)$);
      \draw[blue](a)--(b);
      \draw[black](b)--(c);
      \draw[red](c)--(d);
      \draw[blue](d)--(e);
      \draw[black](e)--(f);
      \draw[red](f)--(a);
      \coordinate (h) at ($(a)!2/3!(b)$);
      \coordinate (j) at ($(-1/3,3^0.5)+(x)$);
      \coordinate (l) at ($(c)!2/3!(d)$);
      \coordinate (m) at ($(d)!1/3!(e)$);
      \coordinate (o) at ($(-1/3,0)+(x)$);
      \coordinate (q) at ($(m)!1/4!(h)$);
      \coordinate (t) at ($(m)!1/2!(h)$);
      \coordinate (r) at ($(j)!1/3!(o)$);
      \coordinate (s) at ($(j)!2/3!(o)$);
      \coordinate (v) at ($(a)!2/3!(f)$);
      \draw[blue] (j)--(r);
      \draw[blue](s)--(o);
      \draw[blue]($(s)+(1/30,0)$)--($(o)+(1/30,0)$);
      \draw[blue]($(j)+(1/30,0)$)--($(r)+(1/30,0)$);
      \draw[red](m)--(q);
      \draw[red](t)--(h);
      \draw[red]($(t)+(1/60,-3^0.5/60)$)--($(a)!2/3-1/30!(b)$);
      \draw[red]($(q)+(1/60,-3^0.5/60)$)--($(e)!2/3-1/30!(d)$);
      \draw[green,name path=E](q)--(r);
      \draw[black](l)--(v);
      \draw[blue,name path=B_2]($(r)+(1/30,0)$)--($(q)+(1/30,0)$);
      \draw[red,name path=A_2]($(r)+(2/30,0)$)--($(q)+(2/30,0)$);
      \draw[black,name path=C_2]($(c)!1/4!(d)$)--($(a)!1/4!(f)$);
\node at ($(d)-(0.3,0)$){$\rightsquigarrow$};
\coordinate (z) at ($(e)!1/2!(f)$);
\node at ($(z)+((0,-0.4)$){$lc$};
\end{tikzpicture}    
\par\end{center}

For Case 1 with $K^{2}=6$, we can degenerate $B_{1}$ to $A_{0}+B_{3}$
to produce the degeneration with $K^{2}=5$, which is the first figure
below. We first blow up the total space along the curve $A_{0}$ in
the central fiber, then blow up the total space along the strict preimage
of $C_{3}$ in the central fiber. The resulting central fiber is $Bl_{3}\mathbb{P}^{2}\cup Bl_{1}\mathbb{F}_{1}\cup\mathbb{F}_{0}$,
which is the second figure below. 

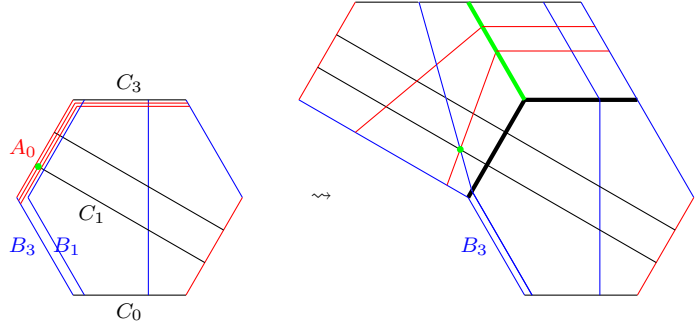
\begin{center}
\begin{tikzpicture}[xscale=1.5,yscale=1.5][font=\tiny]
      \coordinate (z) at (0,0);
      \coordinate (a) at ($(1,3^0.5/2)+(z)$);
      \coordinate (b) at ($(1/2,3^0.5)+(z)$);
      \coordinate (c) at ($(-1/2,3^0.5)+(z)$);
      \coordinate (d) at ($(-1,3^0.5/2)+(z)$);
      \coordinate (e) at ($(-1/2,0)+(z)$);
      \coordinate (f) at ($(1/2,0)+(z)$);
      \draw[blue](a)--(b);
      \draw[black](b)--(c);
      \draw[red](c)--(d);
      \draw[blue](d)--(e);
      \draw[black](e)--(f);
      \draw[red](f)--(a);
      \coordinate (g) at ($(a)!1/3!(b)$);
      \coordinate (h) at ($(a)!2/3!(b)$);
      \coordinate (i) at ($(b)!1/3!(c)$);
      \coordinate (j) at ($(b)!2/3!(c)$);
      \coordinate (k) at ($(c)!1/3!(d)$);
      \coordinate (l) at ($(c)!2/3!(d)$);
      \coordinate (m) at ($(d)!1/3!(e)$);
      \coordinate (n) at ($(d)!2/3!(e)$);
      \coordinate (o) at ($(e)!1/3!(f)$);
      \coordinate (p) at ($(e)!2/3!(f)$);
      \coordinate (q) at ($(f)!1/3!(a)$);
      \coordinate (r) at ($(f)!2/3!(a)$);
      \draw[blue,name path=B_2] (i)--(p);
      \draw[black,name path=C_1](k)--(r);
      \draw[black,name path=C_2](l)--(q);
      \node[above,black]at($(b)!1/2!(c)$){$C_3$};
      \node[below,black]at($(e)!1/2!(f)$){$C_0$};
      \node[left,blue]at($(d)!1/2!(e)$){$B_3$};
      \node[left,red]at($(c)!1/2!(d)$){$A_0$};
      \draw[red,name path=A_1]($(d)!1/30!(e)$)--($(c)!1/60!(f)$);
      \draw[red,name path=A_1]($(c)!1/60!(f)$)--($(b)!1/30!(a)$);
      \draw[red,name path=A_2]($(d)!2/30!(e)$)--($(c)!2/60!(f)$);
      \draw[red,name path=A_2]($(c)!2/60!(f)$)--($(b)!2/30!(a)$);
      \draw[blue,name path=B_1]($(e)!3/30!(f)$)--($(d)!3/60!(a)$);
      \draw[blue,name path=B_1]($(d)!3/60!(a)$)--($(c)!3/30!(b)$);
      \node[black,below]at($(l)!1/3!(q)$){$C_1$};
      \node[blue,right]at($(d)!1/2!(e)$){$B_1$};
      \coordinate (P) at ($(l)!1/60!(q)$);
      \path[fill=green](P)circle[radius=0.03];
      \node at ($(a)+(0.7,0)$){$\rightsquigarrow$};
      \coordinate (z) at (4,0);
      \coordinate (a) at ($(1,3^0.5/2)+(z)$);
      \coordinate (b) at ($(1/2,3^0.5)+(z)$);
      \coordinate (c) at ($(-1/2,3^0.5)+(z)$);
      \coordinate (d) at ($(-1,3^0.5/2)+(z)$);
      \coordinate (e) at ($(-1/2,0)+(z)$);
      \coordinate (f) at ($(1/2,0)+(z)$);
      \draw[blue](a)--(b);
      \draw[black,ultra thick](b)--(c);
      \draw[red](c)--(d);
      \draw[blue](d)--(e);
      \draw[black](e)--(f);
      \draw[red](f)--(a);
      \coordinate (g) at ($(a)!1/3!(b)$);
      \coordinate (h) at ($(a)!2/3!(b)$);
      \coordinate (i) at ($(b)!1/3!(c)$);
      \coordinate (j) at ($(b)!2/3!(c)$);
      \coordinate (k) at ($(c)!1/3!(d)$);
      \coordinate (l) at ($(c)!2/3!(d)$);
      \coordinate (m) at ($(d)!1/3!(e)$);
      \coordinate (n) at ($(d)!2/3!(e)$);
      \coordinate (o) at ($(e)!1/3!(f)$);
      \coordinate (p) at ($(e)!2/3!(f)$);
      \coordinate (q) at ($(f)!1/3!(a)$);
      \coordinate (r) at ($(f)!2/3!(a)$);
      \draw[blue,name path=B_2] (i)--(p);
      \draw[blue,name path=B_1] ($(e)!2/30!(f)$)--($(d)!2/30!(c)$);
      \node[left,blue]at($(d)!1/2!(e)$){$B_3$};
      \draw[black,ultra thick](c)--(d);
      \draw[blue,name path=B_1]($(e)!2/30!(f)$)--($(d)!2/30!(c)$);
      \coordinate (w) at ($(a)!2!(c)$);
      \coordinate (x) at ($(f)!2!(d)$);
      \coordinate (s) at ($(f)!1+1/2!(c)$);
      \coordinate (t) at ($(a)!2!(b)$);
      \draw[black,name path=C_3](w)--(t);
      \draw[blue,name path=B_3](d)--(x);
      \draw[red,name path=A_0](w)--(x);
      \draw[black,name path=C_2](r)--($(x)!2/3!(w)$);
      \draw[black,name path=C_1](q)--($(x)!1/3!(w)$);
      \draw[red,name path=A_2]($(c)!3/4!(s)$)--($(d)!2/3!(x)$);
      \draw[red,name path=A_1]($(c)!1/2!(s)$)--($(d)!1/8!(x)$);
      \draw[green,ultra thick](c)--(s);
      \path [name intersections={of=C_1 and A_1}]; 
      \coordinate (P) at (intersection-1);  
      \coordinate (v) at ($(d)!2/30!(c)$);
      \draw[white,name path=B_1](v)--($(v)!5!(P)$);
      \path [name intersections={of=B_1 and C_3}]; 
      \coordinate (B) at (intersection-1);
      \draw[blue,name path=B_1](v)--(B);
     \draw[blue,name path=B_0](b)--(t);
      \draw[blue,name path=B_2](i)--($(t)!1/3!(s)$);
      \draw[red,name path=A_2]($(c)!3/4!(s)$)--($(b)!3/4!(t)$);
      \draw[red,name path=A_1]($(c)!1/2!(s)$)--($(b)!1/2!(t)$);
      \path[fill=green](P)circle[radius=0.03];
      
\end{tikzpicture}     
\par\end{center}

Consider the curve $B_{3}$ in the component $Bl_{3}\mathbb{P}^{2}$
of the central fiber, we have $\left(K_{\mathcal{Y}}+\dfrac{1}{2}\mathcal{D}\right).B_{3}=-\dfrac{1}{2}<0$
and $K_{\Sigma_{0}}|_{Bl_{3}\mathbb{P}^{2}}.B_{3}=0$. When run the
minimal model program, there will be a flip for $\left(\mathcal{Y},\dfrac{1}{2}\mathcal{D}\right)$.
The normal bundle of $B_{3}$ in the total space is $\mathcal{O}(-1)\oplus\mathcal{O}(-1)$.
The flip for $\left(\mathcal{Y},\dfrac{1}{2}\mathcal{D}\right)$ is
the Atiyah flop for $\mathcal{Y}$. The process is as follows.

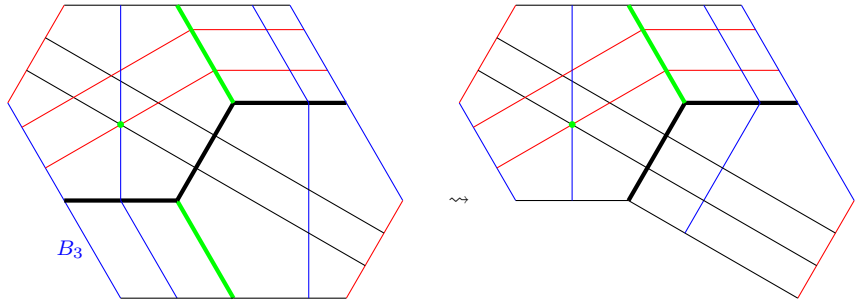
\begin{center}
\begin{tikzpicture}[xscale=1.5,yscale=1.5][font=\tiny]
      \coordinate (z) at (0,0);
      \coordinate (a) at ($(1,3^0.5/2)+(z)$);
      \coordinate (b) at ($(1/2,3^0.5)+(z)$);
      \coordinate (c) at ($(-1/2,3^0.5)+(z)$);
      \coordinate (d) at ($(-1,3^0.5/2)+(z)$);
      \coordinate (e) at ($(-1/2,0)+(z)$);
      \coordinate (f) at ($(1/2,0)+(z)$);
      \draw[blue](a)--(b);
      \draw[black,ultra thick](b)--(c);
      \draw[red](c)--(d);
      \draw[green,ultra thick](d)--(e);
      \draw[black](e)--(f);
      \draw[red](f)--(a);
      \coordinate (g) at ($(a)!1/3!(b)$);
      \coordinate (h) at ($(a)!2/3!(b)$);
      \coordinate (i) at ($(b)!1/3!(c)$);
      \coordinate (j) at ($(b)!2/3!(c)$);
      \coordinate (k) at ($(c)!1/3!(d)$);
      \coordinate (l) at ($(c)!2/3!(d)$);
      \coordinate (m) at ($(d)!1/3!(e)$);
      \coordinate (n) at ($(d)!2/3!(e)$);
      \coordinate (o) at ($(e)!1/3!(f)$);
      \coordinate (p) at ($(e)!2/3!(f)$);
      \coordinate (q) at ($(f)!1/3!(a)$);
      \coordinate (r) at ($(f)!2/3!(a)$);
      \draw[blue,name path=B_2] (i)--(p);
      \draw[black,ultra thick](c)--(d);
      \coordinate (w) at ($(a)!2!(c)$);
      \coordinate (x) at ($(f)!2!(d)$);
      \coordinate (s) at ($(f)!1+1/2!(c)$);
      \coordinate (t) at ($(a)!2!(b)$);
      \coordinate (u) at ($(a)!1+1/2!(d)$);
      \coordinate (v) at ($(f)!2!(e)$);
      \draw[black,name path=C_3](w)--(t);
      \draw[blue,name path=B_3](x)--(v);
      \draw[red,name path=A_0](w)--(x);
      \draw[black,name path=C_2](r)--($(w)!1/3!(x)$);
      \draw[black,name path=C_1](q)--($(x)!1/3!(w)$);
      \draw[red,name path=A_2]($(x)!1/4!(u)$)--($(s)!1/4!(c)$);
      \draw[red,name path=A_2']($(s)!1/4!(c)$)--($(t)!1/4!(b)$);
      \draw[red,name path=A_1]($(u)!1/3!(x)$)--($(c)!1/3!(s)$);
      \draw[red,name path=A_1']($(c)!1/3!(s)$)--($(b)!1/3!(t)$);
      \draw[green,ultra thick](c)--(s);
      \draw[blue,name path=B_0](b)--(t);
      \draw[blue,name path=B_2](i)--($(t)!1/3!(s)$);
      \draw[black,ultra thick](d)--(u);
      \draw[blue,name path=B_1]($(s)!1/2!(w)$)--($(d)!1/2!(u)$);
      \draw[blue,name path=B_1']($(d)!1/2!(u)$)--($(e)!1/2!(v)$);
      \draw[black,name path=C_0'](e)--(v);
      \node[left,blue]at($(u)!1/2!(v)$){$B_3$};
      \path [name intersections={of=B_1 and A_1}]; 
      \coordinate (P) at (intersection-1);  
      \path[fill=green](P)circle[radius=0.03];
\node at ($(a)+(0.5,0)$){$\rightsquigarrow$};
      \coordinate (z) at (4,0);
      \coordinate (a) at ($(1,3^0.5/2)+(z)$);
      \coordinate (b) at ($(1/2,3^0.5)+(z)$);
      \coordinate (c) at ($(-1/2,3^0.5)+(z)$);
      \coordinate (d) at ($(-1,3^0.5/2)+(z)$);
      \coordinate (e) at ($(-1/2,0)+(z)$);
      \coordinate (f) at ($(1/2,0)+(z)$);
      \draw[blue](a)--(b);
      \draw[black,ultra thick](b)--(c);
      \draw[red](c)--(d);
      \draw[red](f)--(a);
      \draw[black,name path=C_0](d)--(f);
      \coordinate (g) at ($(a)!1/3!(b)$);
      \coordinate (h) at ($(a)!2/3!(b)$);
      \coordinate (i) at ($(b)!1/3!(c)$);
      \coordinate (j) at ($(b)!2/3!(c)$);
      \coordinate (k) at ($(c)!1/3!(d)$);
      \coordinate (l) at ($(c)!2/3!(d)$);
      \coordinate (m) at ($(d)!1/3!(e)$);
      \coordinate (n) at ($(d)!2/3!(e)$);
      \coordinate (o) at ($(e)!1/3!(f)$);
      \coordinate (p) at ($(e)!2/3!(f)$);
      \coordinate (q) at ($(f)!1/3!(a)$);
      \coordinate (r) at ($(f)!2/3!(a)$);
      \draw[blue,name path=B_2](i)--($(d)!1/3!(f)$);
      \draw[black,ultra thick](c)--(d);
      \coordinate (w) at ($(a)!2!(c)$);
      \coordinate (x) at ($(f)!2!(d)$);
      \coordinate (s) at ($(f)!1+1/2!(c)$);
      \coordinate (t) at ($(a)!2!(b)$);
      \coordinate (u) at ($(a)!1+1/2!(d)$);
      \coordinate (v) at ($(f)!2!(e)$);
      \draw[black,name path=C_3](w)--(t);
      \draw[blue,name path=B_3](x)--(u);
      \draw[red,name path=A_0](w)--(x);
      \draw[black,name path=C_2](r)--($(w)!1/3!(x)$);
      \draw[black,name path=C_1](q)--($(x)!1/3!(w)$);
      \draw[red,name path=A_2]($(x)!1/4!(u)$)--($(s)!1/4!(c)$);
      \draw[red,name path=A_2']($(s)!1/4!(c)$)--($(t)!1/4!(b)$);
      \draw[red,name path=A_1]($(u)!1/3!(x)$)--($(c)!1/3!(s)$);
      \draw[red,name path=A_1']($(c)!1/3!(s)$)--($(b)!1/3!(t)$);
      \draw[green,ultra thick](c)--(s);
      \draw[blue,name path=B_0](b)--(t);
      \draw[blue,name path=B_2](i)--($(t)!1/3!(s)$);
      \draw[black,name path=C_0](d)--(u);
      \draw[blue,name path=B_1]($(s)!1/2!(w)$)--($(d)!1/2!(u)$);
      \path [name intersections={of=B_1 and A_1}]; 
      \coordinate (P) at (intersection-1);  
      \path[fill=green](P)circle[radius=0.03];

\end{tikzpicture}     
\par\end{center}

After applying the flip to $\left(\mathcal{Y},\dfrac{1}{2}\mathcal{D}\right)$,
the central fiber of the resulting 3-fold space is $Bl_{4}\mathbb{P}^{2}\cup\mathbb{F}_{0}\cup Bl_{2}\mathbb{P}^{2}$.
Finally we blow up the total space along the strict preimage of $L_{P}$.
The general fibers and the central fiber are all blown up at one point.
$ $Now we have $K_{\Sigma_{0}}.C\geq0$ for all the curve $C$ in
$\Sigma_{0}$. Running the minimal model program, both the components
$\mathbb{F}_{0}$ and $Bl_{2}\mathbb{P}^{2}$ in the central fiber
are contracted. The central fiber becomes $Bl_{4}\mathbb{P}^{2}$,
which is a lc degeneration of the general fibers $Bl_{4}\mathbb{P}^{2}$.

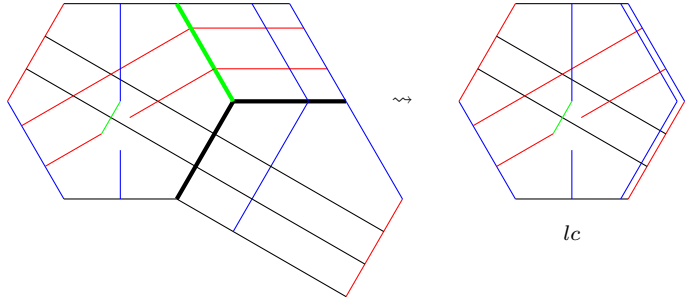
\begin{center}
\begin{tikzpicture}[xscale=1.5,yscale=1.5][font=\tiny]
      \coordinate (z) at (0,0);
      \coordinate (a) at ($(1,3^0.5/2)+(z)$);
      \coordinate (b) at ($(1/2,3^0.5)+(z)$);
      \coordinate (c) at ($(-1/2,3^0.5)+(z)$);
      \coordinate (d) at ($(-1,3^0.5/2)+(z)$);
      \coordinate (e) at ($(-1/2,0)+(z)$);
      \coordinate (f) at ($(1/2,0)+(z)$);
      \draw[blue](a)--(b);
      \draw[black,ultra thick](b)--(c);
      \draw[red](c)--(d);
      \draw[red](f)--(a);
      \draw[black,name path=C_0](d)--(f);
      \coordinate (g) at ($(a)!1/3!(b)$);
      \coordinate (h) at ($(a)!2/3!(b)$);
      \coordinate (i) at ($(b)!1/3!(c)$);
      \coordinate (j) at ($(b)!2/3!(c)$);
      \coordinate (k) at ($(c)!1/3!(d)$);
      \coordinate (l) at ($(c)!2/3!(d)$);
      \coordinate (m) at ($(d)!1/3!(e)$);
      \coordinate (n) at ($(d)!2/3!(e)$);
      \coordinate (o) at ($(e)!1/3!(f)$);
      \coordinate (p) at ($(e)!2/3!(f)$);
      \coordinate (q) at ($(f)!1/3!(a)$);
      \coordinate (r) at ($(f)!2/3!(a)$);
      \draw[blue,name path=B_2](i)--($(d)!1/3!(f)$);
      \draw[black,ultra thick](c)--(d);
      \coordinate (w) at ($(a)!2!(c)$);
      \coordinate (x) at ($(f)!2!(d)$);
      \coordinate (s) at ($(f)!1+1/2!(c)$);
      \coordinate (t) at ($(a)!2!(b)$);
      \coordinate (u) at ($(a)!1+1/2!(d)$);
      \coordinate (v) at ($(f)!2!(e)$);
      \draw[black,name path=C_3](w)--(t);
      \draw[blue,name path=B_3](x)--(u);
      \draw[red,name path=A_0](w)--(x);
      \draw[black,name path=C_2](r)--($(w)!1/3!(x)$);
      \draw[black,name path=C_1](q)--($(x)!1/3!(w)$);
      \draw[red,name path=A_2]($(x)!1/4!(u)$)--($(s)!1/4!(c)$);
      \draw[red,name path=A_2']($(s)!1/4!(c)$)--($(t)!1/4!(b)$);
      \draw[red,name path=A_1']($(c)!1/3!(s)$)--($(b)!1/3!(t)$);
      \draw[green,ultra thick](c)--(s);
      \draw[blue,name path=B_0](b)--(t);
      \draw[blue,name path=B_2](i)--($(t)!1/3!(s)$);
      \draw[black,name path=C_0](d)--(u);
      \coordinate (i) at ($(s)!1/2!(w)$);
      \coordinate (j) at ($(u)!1/3!(x)$);
      \coordinate (k) at ($(d)!1/2!(u)$);
      \coordinate (l) at ($(c)!1/3!(s)$);
      \draw[red,name path=A_1](j)--($(j)!1/3!(l)$);
      \draw[red,name path=A_1'](l)--($(l)!1/2!(j)$);
      \draw[blue,name path=B_1](i)--($(i)!1/2!(k)$); 
      \draw[blue,name path=B_1'](k)--($(k)!1/4!(i)$);
      \draw[green,name path=E]($(j)!1/3!(l)$)--($(i)!1/2!(k)$);
\node at ($(b)+(0.5,0)$){$\rightsquigarrow$};

      \coordinate (z) at (4,0);
      \coordinate (a) at ($(1,3^0.5/2)+(z)$);
      \coordinate (b) at ($(1/2,3^0.5)+(z)$);
      \coordinate (c) at ($(-1/2,3^0.5)+(z)$);
      \coordinate (d) at ($(-1,3^0.5/2)+(z)$);
      \coordinate (e) at ($(-1/2,0)+(z)$);
      \coordinate (f) at ($(1/2,0)+(z)$);
      \coordinate (g) at ($(a)!1/3!(b)$);
      \coordinate (h) at ($(a)!2/3!(b)$);
      \coordinate (i) at ($(b)!1/3!(c)$);
      \coordinate (j) at ($(b)!2/3!(c)$);
      \coordinate (k) at ($(c)!1/3!(d)$);
      \coordinate (l) at ($(c)!2/3!(d)$);
      \coordinate (m) at ($(d)!1/3!(e)$);
      \coordinate (n) at ($(d)!2/3!(e)$);
      \coordinate (o) at ($(e)!1/3!(f)$);
      \coordinate (p) at ($(e)!2/3!(f)$);
      \coordinate (q) at ($(f)!1/3!(a)$);
      \coordinate (r) at ($(f)!2/3!(a)$);
      \draw[red](c)--(d);
      \coordinate (w) at ($(a)!2!(c)$);
      \coordinate (x) at ($(f)!2!(d)$);
      \coordinate (s) at ($(f)!1+1/2!(c)$);
      \coordinate (t) at ($(a)!2!(b)$);
      \coordinate (u) at ($(a)!1+1/2!(d)$);
      \coordinate (v) at ($(f)!2!(e)$);
      \draw[black,name path=C_3](w)--(s);
      \draw[blue,name path=B_3](x)--(u);
      \draw[red,name path=A_0](w)--(x);
      \draw[black,name path=C_2]($(c)!1/3!(d)$)--($(w)!1/3!(x)$);
      \draw[black,name path=C_1]($(c)!2/3!(d)$)--($(x)!1/3!(w)$);
      \draw[red,name path=A_2]($(x)!1/4!(u)$)--($(s)!1/4!(c)$);
      \draw[blue](c)--(s);
      \draw[black,name path=C_0](d)--(u);
      \draw[blue]($(s)!2/30!(w)$)--($(c)!2/60!(x)$);
      \draw[blue]($(c)!2/60!(x)$)--($(d)!2/30!(u)$);

      \coordinate (i) at ($(s)!1/2!(w)$);
      \coordinate (j) at ($(u)!1/3!(x)$);
      \coordinate (k) at ($(d)!1/2!(u)$);
      \coordinate (l) at ($(c)!1/3!(s)$);
      \draw[red,name path=A_1](j)--($(j)!1/3!(l)$);
      \draw[red,name path=A_1'](l)--($(l)!1/2!(j)$);
      \draw[blue,name path=B_1](i)--($(i)!1/2!(k)$); 
      \draw[blue,name path=B_1'](k)--($(k)!1/4!(i)$);
      \draw[green,name path=E]($(j)!1/3!(l)$)--($(i)!1/2!(k)$);
  \node at ($(k)+(0,-0.3)$){$lc$};

\end{tikzpicture}     
\par\end{center}

From the cases discussed above, we conclude that there are 6 types
of degenerate configurations with reduced log canonical models for
the moduli space of Burniat surfaces with $K^{2}=5$, up to the symmetry
group $\mathbb{Z}_{6}$. All of the 6 cases could be obtained from
the degenerating cases for $K^{2}=6$ listed in \cite{AP09}, with
the additional condition that $A_{1},B_{1},C_{1}$ meet at a point
$P$.\emph{ }All the lc models of the degenerate configurations come
from Case 1 and 8 for $K^{2}=6$ are irreducible. Case 5 for $K^{2}=6$
does not produce any degenerations for $K^{2}=5$. We give a table
with the relations between cases for $K^{2}=5$ and cases for $K^{2}=6$.
This table describes how to get cases with $K^{2}=5$ possibly from
cases with $K^{2}=6$ with $A_{1},B_{1},C_{1}$ meeting at one point
$P$\emph{.}

\begin{center}
\begin{tabular}{|c|c|c|}
\hline 
$K^{2}=6$ & $K^{2}=5$ & $K^{2}=5$, further degenerations \tabularnewline
\hline 
\hline 
Case 1 & lc & \tabularnewline
\hline 
Case 2 & Case 1 & Case 2,3\tabularnewline
\hline 
Case 3 & Case 2 & Case 3\tabularnewline
\hline 
Case 4 & Case 3 & \tabularnewline
\hline 
Case 6 & Case 4 & Case 5,3\tabularnewline
\hline 
Case 5  & none & \tabularnewline
\hline 
Case 7 & Case 5 & Case 3\tabularnewline
\hline 
Case 8 & lc & \tabularnewline
\hline 
Case 9,10 & Case 6 & \tabularnewline
\hline 
\end{tabular}
\par\end{center}

\section{Burniat surfaces with $K^{2}=4$\label{degree4}}

We consider $\mathbb{P}^{2}$ with $9$ lines. There are two cases
with two distinct points $P_{1},P_{2}$ which are the intersections
of three lines inside the triangle. We denote these two cases as a
``nodal case'' and a ``non-nodal case''. \\

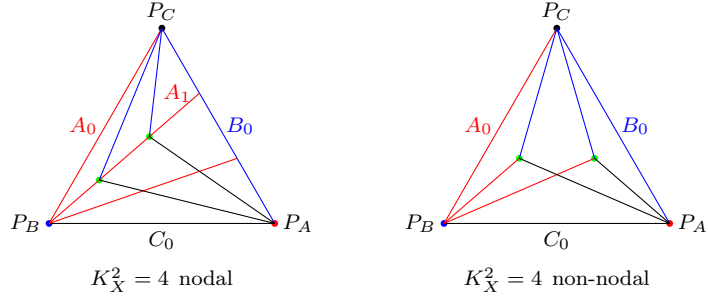
\begin{center}
\begin{tikzpicture}[xscale=1.5,yscale=1.5][font=\tiny]
      \coordinate (P_A) at (1,0);
      \node[right] at (P_A){$P_A$};
      \path[fill=red](P_A)circle[radius=0.03];
      \coordinate (P_B) at (-1,0);
      \node[left] at (P_B){$P_B$};
      \path[fill=blue](P_B)circle[radius=0.03];
      \coordinate (P_C) at (0,3^0.5);
      \node[above] at (P_C){$P_C$};
      \path[fill=black](P_C)circle[radius=0.03];
      \draw[black](P_A)--(P_B);
      \node[below,black]at($(P_A)!1/2!(P_B)$){$C_0$};
      \draw[red](P_B)--(P_C);
      \node[left,red]at($(P_B)!1/2!(P_C)$){$A_0$};
      \draw[blue](P_C)--(P_A);
      \node[right,blue]at($(P_C)!1/2!(P_A)$){$B_0$};
      \coordinate (a) at ($(P_A)!1/3!(P_C)$);
      \coordinate (b) at ($(P_A)!2/3!(P_C)$);
      \coordinate (c) at ($(b)!1/3!(P_B)$);
      \coordinate (d) at ($(b)!2/3!(P_B)$);
      \path[fill=green](c)circle[radius=0.03];
      \path[fill=green](d)circle[radius=0.03];
      \draw [name path=A_1,red](b)--(P_B);
      \draw [name path=A_2,red](a)--(P_B);
      \node [left,red] at (b){$A_1$};
      \draw [name path=B_1,blue](P_C)--(d);
      \draw [name path=B_2,blue](P_C)--(c);
      \draw [name path=C_1,black](P_A)--(d);
      \draw [name path=C_2,black](P_A)--(c);
      \coordinate (e) at ($(P_A)!1/2!(P_B)$);
      \node at ($(e)-(0,1/2)$) {$K_{X}^{2}=4\:\:\mbox{nodal}$};
      \coordinate (P_A) at (1+3.5,0);
      \node[right] at (P_A){$P_A$};
      \path[fill=red](P_A)circle[radius=0.03];
      \coordinate (P_B) at (-1+3.5,0);
      \node[left] at (P_B){$P_B$};
      \path[fill=blue](P_B)circle[radius=0.03];
      \coordinate (P_C) at (0+3.5,3^0.5);
      \node[above] at (P_C){$P_C$};
      \path[fill=black](P_C)circle[radius=0.03];
      \draw[black](P_A)--(P_B);
      \node[below,black]at($(P_A)!1/2!(P_B)$){$C_0$};
      \draw[red](P_B)--(P_C);
      \node[left,red]at($(P_B)!1/2!(P_C)$){$A_0$};
      \draw[blue](P_C)--(P_A);
      \node[right,blue]at($(P_C)!1/2!(P_A)$){$B_0$};
      \coordinate (a) at ($(P_A)!1/3!(P_C)$);
      \coordinate (b) at ($(P_B)!1/3!(P_C)$);
      \coordinate (c) at ($(a)!1/4!(b)$);
      \coordinate (d) at ($(a)!3/4!(b)$);
      \path[fill=green](c)circle[radius=0.03];
      \path[fill=green](d)circle[radius=0.03];
      \draw[red](P_B)--(c);
      \draw[red](P_B)--(d);
      \draw[blue](P_C)--(c);
      \draw[blue](P_C)--(d);
      \draw[black](P_A)--(c);
      \draw[black](P_A)--(d);
      \coordinate (e) at ($(P_A)!1/2!(P_B)$);
      \node at ($(e)-(0,1/2)$) {$K_{X}^{2}=4\:\:\mbox{non-nodal}$};
\end{tikzpicture}    
\par\end{center}

Let $\Sigma=Bl_{5}\mathbb{P}^{2}$ be the blowup of $\mathbb{P}^{2}$
at 5 points $P_{A},P_{B},P_{C}$ and $P_{1},P_{2}$.\\

\begin{defn}
A Burniat surface $X$ in $M_{Bur}^{3}$ is the canonical model of
a $\mathbb{Z}_{2}^{2}$-cover of $\Sigma=\mbox{Bl}_{5}\mathbb{P}^{2}$
for the building data $D_{a}=\sum_{i=0}^{3}A_{i}$, $D_{b}=\sum_{i=0}^{3}B_{i}$,
$D_{c}=\sum_{i=0}^{3}C_{i}$, where $a,b,c$ are the $3$ nonzero
elements of $\mathbb{Z}_{2}^{2}$.
\end{defn}
Let $D=\dfrac{1}{2}\left(D_{a}+D_{b}+D_{c}\right)$, then we have
\[
K_{X}^{2}=\left(\pi^{*}(K_{\Sigma}+D)\right)^{2}=\left(\pi^{*}(-\dfrac{1}{2}K_{\Sigma})\right)^{2}=4\left(\dfrac{1}{4}K_{\Sigma}^{2}\right)=4.
\]

For the nodal case, the curve $A_{1}$ in $\Sigma$ is a $\left(-2\right)$-curve
and $K_{\Sigma}.A_{1}=0$. The anti-canonical divisor $-K_{\Sigma}$
is nef but not ample, so $K_{\Sigma}+D=-\dfrac{1}{2}K_{\Sigma}$ is
not ample which implies $X$ is not ample. Stable Burniat surfaces
$X$ with $K_{X}^{2}=4$ are $\mathbb{Z}_{2}^{2}$-covers of the canonical
models $\Sigma^{c}$ of $\Sigma$ with the building data $\frac{1}{2}D$.

For the non-nodal case, we have that $X$ is stable as that $-K_{\Sigma}$
is ample, and stable Burniat surfaces $X$ with $K_{X}^{2}=4$ are
$\mathbb{Z}_{2}^{2}$-covers of $\Sigma$ with the building data $\frac{1}{2}D$.

To compactify the moduli space of stable pairs $(Y,\frac{1}{2}D)$,
we will study one-parameter families of configurations in the moduli
space. For the nodal case, the general fiber $\Sigma^{c}$ is the
blown down of a (-2)-curve $A_{1}$ of $\Sigma=Bl_{5}\mathbb{P}^{2}$;
for the non-nodal case, the general fiber $\Sigma$ is $Bl_{5}\mathbb{P}^{2}$.

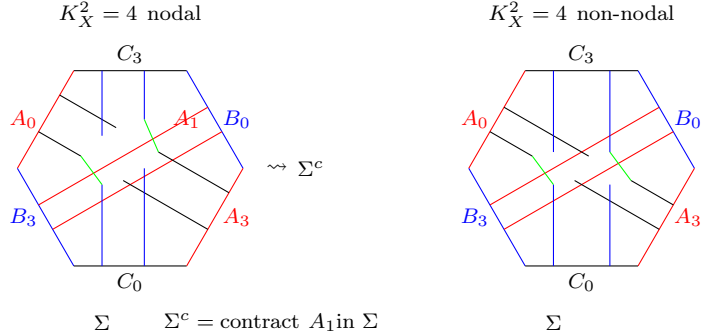
\begin{center}
\begin{tikzpicture}[xscale=1.5,yscale=1.5][font=\tiny]
      \coordinate (a) at (1,0);
      \coordinate (b) at (1/2,3^0.5/2);
      \coordinate (c) at (-1/2,3^0.5/2);
      \coordinate (d) at (-1,0);
      \coordinate (e) at (-1/2,-3^0.5/2);
      \coordinate (f) at (1/2,-3^0.5/2);
      \draw[blue](a)--(b);
      \draw[black](b)--(c);
      \draw[red](c)--(d);
      \draw[blue](d)--(e);
      \draw[black](e)--(f);
      \draw[red](f)--(a);
      \coordinate (g) at ($(a)!5/8!(b)$);
      \coordinate (h) at ($(b)!3/8!(c)$);
      \coordinate (i) at ($(b)!3/4!(c)$);
      \coordinate (j) at ($(c)!1/4!(d)$);
      \coordinate (k) at ($(c)!5/8!(d)$);
      \coordinate (l) at ($(d)!3/8!(e)$);
      \coordinate (m) at ($(e)!1/4!(f)$);
      \coordinate (n) at ($(e)!5/8!(f)$);
      \coordinate (o) at ($(f)!3/8!(a)$);
      \coordinate (p) at ($(f)!3/4!(a)$);
      \draw [name path=A_1,red](g)--(l);
      \draw [name path=A_2,red]($(a)!3/8!(b)$)--($(e)!3/8!(d)$);
      \draw [name path=C_1,black](k)--($(k)!1/4!(o)$);
      \draw [name path=C_1',black] (o)--($(o)!1/2!(k)$);
      \draw [name path=C_2,black](j)--($(j)!1/3!(p)$);
      \draw [name path=C_2',black](p)--($(p)!5/12!(j)$);
      \draw [name path=B_1,blue](i)--($(i)!1/3!(m)$);
      \draw [name path=B_1',blue](m)--($(m)!5/12!(i)$);
      \draw [name path=B_2,blue](h)--($(h)!1/4!(n)$);
      \draw [name path=B_2',blue](n)--($(n)!1/2!(h)$);
      \draw [name path=E_1,green]($(k)!1/4!(o)$)--($(m)!5/12!(i)$);
      \draw [name path=E_2,green]($(h)!1/4!(n)$)--($(p)!5/12!(j)$);
      \node[above,black]at($(b)!1/2!(c)$){$C_3$};
      \node[below,black]at($(e)!1/2!(f)$){$C_0$};
      \node[right,blue]at($(a)!1/2!(b)$){$B_0$};
      \node[left,blue]at($(d)!1/2!(e)$){$B_3$};
      \node[left,red]at($(c)!1/2!(d)$){$A_0$};
      \node[right,red]at($(a)!1/2!(f)$){$A_3$};
      \node[red]at($(b)!1/4!(f)$){$A_1$};
\node at ($(a)+(0.3,0)$){$\rightsquigarrow$};
\node at ($(a)+(0.6,0)$){$\Sigma^{c}$};
\node at ($(m)+(0,-1/2)$) {$\Sigma$};
\node at ($(m)+(1.5,-1/2)$) {$\Sigma^{c}=\mbox{contract \ensuremath{A_{1}}in \ensuremath{\Sigma}}$};
      \node at ($(i)+(0+1/4,1/2)$) {$K_{X}^{2}=4\:\:\mbox{nodal}$};
      \coordinate (a) at (1+4,0);
      \coordinate (b) at (1/2+4,3^0.5/2);
      \coordinate (c) at (-1/2+4,3^0.5/2);
      \coordinate (d) at (-1+4,0);
      \coordinate (e) at (-1/2+4,-3^0.5/2);
      \coordinate (f) at (1/2+4,-3^0.5/2);
      \draw[blue](a)--(b);
      \draw[black](b)--(c);
      \draw[red](c)--(d);
      \draw[blue](d)--(e);
      \draw[black](e)--(f);
      \draw[red](f)--(a);
      \coordinate (g) at ($(a)!5/8!(b)$);
      \coordinate (h) at ($(b)!1/4!(c)$);
      \coordinate (i) at ($(b)!3/4!(c)$);
      \coordinate (j) at ($(c)!3/8!(d)$);
      \coordinate (k) at ($(c)!5/8!(d)$);
      \coordinate (l) at ($(d)!3/8!(e)$);
      \coordinate (m) at ($(e)!1/4!(f)$);
      \coordinate (n) at ($(e)!3/4!(f)$);
      \coordinate (o) at ($(f)!3/8!(a)$);
      \coordinate (p) at ($(f)!5/8!(a)$);
      \draw [name path=A_1,red](g)--(l);
      \draw [name path=A_2,red]($(a)!3/8!(b)$)--($(e)!3/8!(d)$);
      \draw [name path=C_1,black](k)--($(k)!1/4!(o)$);
      \draw [name path=C_1',black](o)--($(o)!1/2!(k)$);
      \draw [name path=C_2,black](j)--($(j)!1/2!(p)$);
      \draw [name path=C_2',black](p)--($(p)!1/4!(j)$);

      \draw [name path=B_1,blue](i)--($(i)!5/12!(m)$);
      \draw [name path=B_1',blue](m)--($(m)!5/12!(i)$);
      \draw [name path=B_2,blue](h)--($(h)!5/12!(n)$);
      \draw [name path=B_2',blue](n)--($(n)!5/12!(h)$);

      \draw [name path=E_1,green]($(k)!1/4!(o)$)--($(m)!5/12!(i)$);
      \draw [name path=E_2,green]($(h)!5/12!(n)$)--($(p)!1/4!(j)$);
      \node[above,black]at($(b)!1/2!(c)$){$C_3$};
      \node[below,black]at($(e)!1/2!(f)$){$C_0$};
      \node[right,blue]at($(a)!1/2!(b)$){$B_0$};
      \node[left,blue]at($(d)!1/2!(e)$){$B_3$};
      \node[left,red]at($(c)!1/2!(d)$){$A_0$};
      \node[right,red]at($(a)!1/2!(f)$){$A_3$};
\node at ($(m)+(0,-1/2)$) {$\Sigma$};
      \node at ($(i)+(0+1/4,1/2)$) {$K_{X}^{2}=4\:\:\mbox{non-nodal}$};
\end{tikzpicture}
\par\end{center}

The general fiber $\Sigma^{c}$ is a singular surface with an $A_{1}$-singularity,
which is obtained from $Bl_{5}\mathbb{P}^{2}$ by contracting the
(-2)-curve. To see the degenerating arrangements with $K^{2}=4$,
we will start with surfaces $Bl_{3}\mathbb{P}^{2}$ which are shown
in the following figures.\\

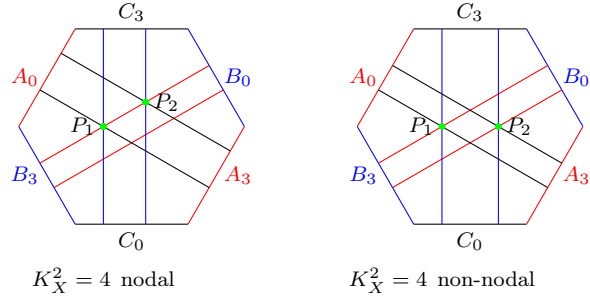
\begin{center}
\begin{tikzpicture}[xscale=1.5,yscale=1.5][font=\tiny]
      \coordinate (a) at (1,0);
      \coordinate (b) at (1/2,3^0.5/2);
      \coordinate (c) at (-1/2,3^0.5/2);
      \coordinate (d) at (-1,0);
      \coordinate (e) at (-1/2,-3^0.5/2);
      \coordinate (f) at (1/2,-3^0.5/2);
      \draw[blue](a)--(b);
      \draw[black](b)--(c);
      \draw[red](c)--(d);
      \draw[blue](d)--(e);
      \draw[black](e)--(f);
      \draw[red](f)--(a);
      \coordinate (g) at ($(a)!5/8!(b)$);
      \coordinate (h) at ($(b)!3/8!(c)$);
      \coordinate (i) at ($(b)!3/4!(c)$);
      \coordinate (j) at ($(c)!1/4!(d)$);
      \coordinate (k) at ($(c)!5/8!(d)$);
      \coordinate (l) at ($(d)!3/8!(e)$);
      \coordinate (m) at ($(e)!1/4!(f)$);
      \coordinate (n) at ($(e)!5/8!(f)$);
      \coordinate (o) at ($(f)!3/8!(a)$);
      \coordinate (p) at ($(f)!3/4!(a)$);
      \draw [name path=A_1,red](g)--(l);
      \draw [name path=A_2,red]($(a)!3/8!(b)$)--($(e)!3/8!(d)$);
      \draw [name path=C_1,black](k)--(o);
      \draw [name path=C_2,black](j)--(p);
      \draw [name path=B_1,blue](i)--(m);
      \draw [name path=B_2,blue](h)--(n);
      \node[above,black]at($(b)!1/2!(c)$){$C_3$};
      \node[below,black]at($(e)!1/2!(f)$){$C_0$};
      \node[right,blue]at($(a)!1/2!(b)$){$B_0$};
      \node[left,blue]at($(d)!1/2!(e)$){$B_3$};
      \node[left,red]at($(c)!1/2!(d)$){$A_0$};
      \node[right,red]at($(a)!1/2!(f)$){$A_3$};
      \path [name intersections={of=A_1 and B_1}]; 
      \coordinate [label=left:$P_1$] (P_1) at (intersection-1);
      \path[fill=green](P_1)circle[radius=0.03];
      \path [name intersections={of=A_1 and B_2}]; 
      \coordinate [label=right:$P_2$] (P_2) at (intersection-1);
      \path[fill=green](P_2)circle[radius=0.03];

      \path[fill=green]($(b)!3/8!(e)$)circle[radius=0.03];
      \path[fill=green]($(a)!5/8!(d)$)circle[radius=0.03];
      \node at ($(m)-(0,1/2)$) {$K_{X}^{2}=4\:\:\mbox{nodal}$};
      \coordinate (a) at (1+3,0);
      \coordinate (b) at (1/2+3,3^0.5/2);
      \coordinate (c) at (-1/2+3,3^0.5/2);
      \coordinate (d) at (-1+3,0);
      \coordinate (e) at (-1/2+3,-3^0.5/2);
      \coordinate (f) at (1/2+3,-3^0.5/2);
      \draw[blue](a)--(b);
      \draw[black](b)--(c);
      \draw[red](c)--(d);
      \draw[blue](d)--(e);
      \draw[black](e)--(f);
      \draw[red](f)--(a);
      \coordinate (g) at ($(a)!5/8!(b)$);
      \coordinate (h) at ($(b)!1/4!(c)$);
      \coordinate (i) at ($(b)!3/4!(c)$);
      \coordinate (j) at ($(c)!3/8!(d)$);
      \coordinate (k) at ($(c)!5/8!(d)$);
      \coordinate (l) at ($(d)!3/8!(e)$);
      \coordinate (m) at ($(e)!1/4!(f)$);
      \coordinate (n) at ($(e)!3/4!(f)$);
      \coordinate (o) at ($(f)!3/8!(a)$);
      \coordinate (p) at ($(f)!5/8!(a)$);
      \draw [name path=A_1,red](g)--(l);
      \draw [name path=A_2,red]($(a)!3/8!(b)$)--($(e)!3/8!(d)$);
      \draw [name path=C_1,black](k)--(o);
      \draw [name path=C_2,black](j)--(p);
      \draw [name path=B_1,blue](i)--(m);
      \draw [name path=B_2,blue](h)--(n);
      \node[above,black]at($(b)!1/2!(c)$){$C_3$};
      \node[below,black]at($(e)!1/2!(f)$){$C_0$};
      \node[right,blue]at($(a)!1/2!(b)$){$B_0$};
      \node[left,blue]at($(d)!1/2!(e)$){$B_3$};
      \node[left,red]at($(c)!1/2!(d)$){$A_0$};
      \node[right,red]at($(a)!1/2!(f)$){$A_3$};
      \path [name intersections={of=A_1 and B_1}]; 
      \coordinate [label=left:$P_1$] (P_1) at (intersection-1);
      \path[fill=green](P_1)circle[radius=0.03];
      \path [name intersections={of=A_2 and B_2}]; 
      \coordinate [label=right:$P_2$] (P_2) at (intersection-1);
      \path[fill=green](P_2)circle[radius=0.03];
      \node at ($(m)-(0,1/2)$) {$K_{X}^{2}=4\:\:\mbox{non-nodal}$};
\end{tikzpicture}
\par\end{center}

We first consider the nodal case with $K^{2}=4$.

\textbf{Case 1.} When the curve $A_{2}$ degenerates to $A_{0}+C_{3}$,
$B_{2}$ degenerates to $A_{0}+B_{3}$, and $C_{2}$ degenerates to
$B_{3}+C_{0}$. Blowing up the total space $\mathcal{Y}$ along the
curve $L_{P_{1}}$ and the curve $A_{0}$ in $\Sigma_{0}$, we see
the general fibers are $Bl_{4}\mathbb{P}^{2}$ and the central fiber
is $\Sigma_{0}=Bl_{4}\mathbb{P}^{2}\cup\mathbb{F}_{1}$. Next we blow
up the total space along the strict preimage of $B_{3}$ in the component
$Bl_{4}\mathbb{P}^{2}$ of $\Sigma_{0}$ and along the strict transform
$\tilde{L}_{P_{2}}$ , which results in the central fiber becoming
a union of three components $Bl_{4}\mathbb{P}^{2}\cup Bl_{1}\mathbb{F}_{1}\cup\mathbb{F}_{0}$
(see the first figure of the second row below). Running the minimal
model program, we get the lc model with central fiber $\mathbb{F}_{0}\cup\mathbb{P}^{2}\cup\mathbb{P}^{2}$
and we call it Case 1. The further degeneration is 4 copies of $\mathbb{P}^{2}$
and we call it Case 2. 

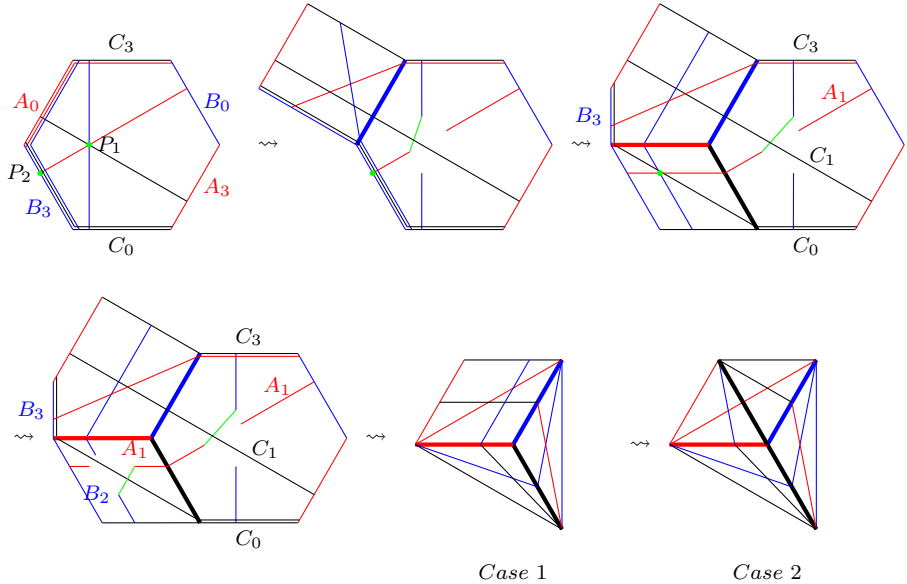
\begin{center}
\begin{tikzpicture}[xscale=1.3,yscale=1.3][font=\tiny]
      \coordinate (a) at (1,3^0.5/2);
      \coordinate (b) at (1/2,3^0.5);
      \coordinate (c) at (-1/2,3^0.5);
      \coordinate (d) at (-1,3^0.5/2);
      \coordinate (e) at (-1/2,0);
      \coordinate (f) at (1/2,0);
      \draw[blue](a)--(b);
      \draw[black](b)--(c);
      \draw[red](c)--(d);
      \draw[blue](d)--(e);
      \draw[black](e)--(f);
      \draw[red](f)--(a);
      \coordinate (g) at ($(a)!1/3!(b)$);
      \coordinate (h) at ($(a)!2/3!(b)$);
      \coordinate (i) at ($(b)!1/3!(c)$);
      \coordinate (j) at (-1/3,3^0.5);
      \coordinate (k) at ($(c)!1/3!(d)$);
      \coordinate (l) at ($(c)!2/3!(d)$);
      \coordinate (m) at ($(d)!1/3!(e)$);
      \coordinate (n) at ($(d)!2/3!(e)$);
      \coordinate (o) at (-1/3,0);
      \coordinate (p) at ($(e)!2/3!(f)$);
      \coordinate (q) at ($(m)!1/4!(h)$);
      \coordinate (t) at ($(m)!1/2!(h)$);
      \coordinate (r) at ($(j)!1/3!(o)$);
      \coordinate (s) at ($(j)!2/3!(o)$);
      \coordinate (u) at ($(a)!1/3!(f)$);
      \coordinate (v) at ($(a)!2/3!(f)$);
      \draw[blue,name path=B_1] (j)--(o);
      \draw[red,name path=A_1](m)--(h);
      \draw[black,name path=C_1](l)--(v);
      \node[above,black]at($(b)!1/2!(c)$){$C_3$};
      \node[below,black]at($(e)!1/2!(f)$){$C_0$};
      \node[right,blue]at($(a)!1/2!(b)$){$B_0$};
      \node[left,blue]at($(d)!3/4!(e)$){$B_3$};
      \node[left,red]at($(c)!1/2!(d)$){$A_0$};
      \node[right,red]at($(a)!1/2!(f)$){$A_3$};
      \path [name intersections={of=A_1 and B_1}]; 
      \coordinate (P) at (intersection-1);
      \path[fill=green](P)circle[radius=0.03];
      \node[right,black] at(P){$P_1$};
      \path[fill=green](m)circle[radius=0.03];
      \node[left,black] at (m){$P_2$};
      \draw[red]($(d)!1/30!(e)$)--($(c)!1/60!(f)$);
      \draw[red]($(c)!1/60!(f)$)--($(b)!1/30!(a)$);
      \draw[blue]($(e)!2/30!(f)$)--($(d)!2/60!(a)$);
      \draw[blue]($(d)!2/60!(a)$)--($(c)!2/30!(b)$);
      \draw[black]($(d)!1/30!(c)$)--($(e)!1/60!(b)$);
      \draw[black]($(e)!1/60!(b)$)--($(f)!1/30!(a)$);
  \node at ($(a)+(0.5,0)$){$\rightsquigarrow$};
      \coordinate (a) at (1+3.4,3^0.5/2);
      \coordinate (b) at (1/2+3.4,3^0.5);
      \coordinate (c) at (-1/2+3.4,3^0.5);
      \coordinate (d) at (-1+3.4,3^0.5/2);
      \coordinate (e) at (-1/2+3.4,0);
      \coordinate (f) at (1/2+3.4,0);
      \draw[blue](a)--(b);
      \draw[black](b)--(c);
      \draw[ultra thick,blue](c)--(d);
      \draw[blue](d)--(e);
      \draw[black](e)--(f);
      \draw[red](f)--(a);
      \coordinate (g) at ($(a)!1/3!(b)$);
      \coordinate (h) at ($(a)!2/3!(b)$);
      \coordinate (i) at ($(b)!1/3!(c)$);
      \coordinate (j) at (-1/3+3.4,3^0.5);
      \coordinate (k) at ($(c)!1/3!(d)$);
      \coordinate (l) at ($(c)!2/3!(d)$);
      \coordinate (m) at ($(d)!1/3!(e)$);
      \coordinate (n) at ($(d)!2/3!(e)$);
      \coordinate (o) at (-1/3+3.4,0);
      \coordinate (p) at ($(e)!2/3!(f)$);
      \coordinate (q) at ($(m)!1/4!(h)$);
      \coordinate (t) at ($(m)!1/2!(h)$);
      \coordinate (r) at ($(j)!1/3!(o)$);
      \coordinate (s) at ($(j)!2/3!(o)$);
      \coordinate (u) at ($(a)!1/3!(f)$);
      \coordinate (v) at ($(a)!2/3!(f)$);
      \draw[blue] (j)--(r);
      \draw[blue](s)--(o);
      \draw[red](m)--(q);
      \draw[red](t)--(h);
      \draw[green](q)--(r);
      \draw[black](l)--(v);
      \draw[red]($(c)!1/30!(d)$)--($(b)!1/30!(a)$);
      \draw[blue]($(e)!2/30!(f)$)--($(d)!2/30!(c)$);
      \draw[black]($(d)!1/30!(c)$)--($(e)!1/60!(b)$);
      \draw[black]($(e)!1/60!(b)$)--($(f)!1/30!(a)$);
   \path[fill=green](m)circle[radius=0.03];
      \coordinate (w) at ($(a)!1+2/3!(c)$);
      \coordinate (x) at ($(f)!1+2/3!(d)$);
      \draw[black](c)--(w);
      \draw[black](l)--($(v)!1+2/3!(l)$);
      \draw[black]($(d)!1/30!(c)$)--($(x)!1/30!(w)$);
      \draw[red](w)--(x);
      \draw[red]($(c)!1/30!(d)$)--($(d)!2/3!(x)$);
      \draw[blue]($(d)!2/30!(c)$)--($(c)!2/3!(w)$);
      \draw[blue](d)--(x);
 \node at ($(a)+(0.3,0)$){$\rightsquigarrow$};
      \coordinate (a) at (1+3+4,3^0.5/2);
      \coordinate (b) at (1/2+3+4,3^0.5);
      \coordinate (c) at (-1/2+3+4,3^0.5);
      \coordinate (d) at (-1+3+4,3^0.5/2);
      \coordinate (e) at (-1/2+3+4,0);
      \coordinate (f) at (1/2+3+4,0);
      \draw[blue](a)--(b);
      \draw[black](b)--(c);
      \draw[ultra thick,blue](c)--(d);
      \draw[ultra thick,black](d)--(e);
      \draw[black](e)--(f);
      \draw[red](f)--(a);
      \coordinate (g) at ($(a)!1/3!(b)$);
      \coordinate (h) at ($(a)!2/3!(b)$);
      \coordinate (i) at ($(b)!1/3!(c)$);
      \coordinate (j) at (-1/3+3+4.2,3^0.5);
      \coordinate (k) at ($(c)!1/3!(d)$);
      \coordinate (l) at ($(c)!2/3!(d)$);
      \coordinate (m) at ($(d)!1/3!(e)$);
      \coordinate (n) at ($(d)!2/3!(e)$);
      \coordinate (o) at (-1/3+3+4.2,0);
      \coordinate (p) at ($(e)!2/3!(f)$);
      \coordinate (q) at ($(m)!1/4!(h)$);
      \coordinate (t) at ($(m)!1/2!(h)$);
      \coordinate (r) at ($(j)!1/3!(o)$);
      \coordinate (s) at ($(j)!2/3!(o)$);
      \coordinate (u) at ($(a)!1/3!(f)$);
      \coordinate (v) at ($(a)!2/3!(f)$);
      \draw[blue,name path=B_1] (j)--(r);
      \draw[blue,name path=B_1'](s)--(o);
      \draw[red,name path=A_1](m)--(q);
      \draw[red,name path=A_1'](t)--(h);
      \draw[green,name path=E](q)--(r);
      \draw[black,name path=C_1](l)--(v);
      \draw[red]($(c)!1/30!(d)$)--($(b)!1/30!(a)$);
      \draw[black]($(e)!1/30!(d)$)--($(f)!1/30!(a)$);
      \node[black,above] at ($(b)!1/2!(c)$){$C_3$};
      \node[black,below] at ($(e)!1/2!(f)$){$C_0$};
      \node[black,above] at ($(v)!1/3!(l)$){$C_1$};
      \node[red,above] at ($(h)!1/4!(m)$){$A_1$};
      \coordinate (w) at ($(a)!1+2/3!(c)$);
      \coordinate (x) at ($(f)!1+2/3!(d)$);
      \coordinate (y) at ($(a)!1+1/2!(d)$);
      \draw[black](c)--(w);
      \draw[black](l)--($(v)!1+2/3!(l)$);
      \draw[black]($(x)!2/30!(w)$)--($(y)!1/30!(d)$);
      \draw[blue]($(c)!1/2!(w)$)--($(d)!2/3!(y)$);
      \draw[red](w)--(x);
      \draw[red]($(c)!1/30!(d)$)--($(x)!2/3!(y)$);
      \draw[red,ultra thick](d)--(y);
      \draw[blue,name path=B_3](x)--(y);
      \node[blue,left] at ($(x)!1/2!(y)$){$B_3$};
      \coordinate (z) at ($(f)!2!(e)$);
      \draw[blue](y)--(z);
      \draw[black](e)--(z);
      \draw[black,name path=C_2]($(y)!1/30!(d)$)--($(e)!1/30!(d)$);
      \draw[blue,name path=B_2]($(d)!2/3!(y)$)--($(z)!1/3!(e)$);
      \draw[red,name path=A_1]($(m)$)--($(y)!1/3!(z)$);
      \path [name intersections={of=A_1 and B_2}]; 
      \coordinate (P) at (intersection-1);
      \path[fill=green](P)circle[radius=0.03];  
      \coordinate (x) at (1.3,-3);
      \coordinate (a) at ($(1,3^0.5/2)+(x)$);
      \coordinate (b) at ($(1/2,3^0.5)+(x)$);
      \coordinate (c) at ($(-1/2,3^0.5)+(x)$);
      \coordinate (d) at ($(-1,3^0.5/2)+(x)$);
      \coordinate (e) at ($(-1/2,0)+(x)$);
      \coordinate (f) at ($(1/2,0)+(x)$);
      \node at ($(d)+(-1.3,0)$){$\rightsquigarrow$};
      \draw[blue](a)--(b);
      \draw[black](b)--(c);
      \draw[ultra thick,blue](c)--(d);
      \draw[ultra thick,black](d)--(e);
      \draw[black](e)--(f);
      \draw[red](f)--(a);
      \coordinate (g) at ($(a)!1/3!(b)$);
      \coordinate (h) at ($(a)!2/3!(b)$);
      \coordinate (i) at ($(b)!1/3!(c)$);
      \coordinate (j) at ($(-1/3+0.2,3^0.5)+(x)$);
      \coordinate (k) at ($(c)!1/3!(d)$);
      \coordinate (l) at ($(c)!2/3!(d)$);
      \coordinate (m) at ($(d)!1/3!(e)$);
      \coordinate (n) at ($(d)!2/3!(e)$);
      \coordinate (o) at ($(-1/3+0.2,0)+(x)$);
      \coordinate (p) at ($(e)!2/3!(f)$);
      \coordinate (q) at ($(m)!1/4!(h)$);
      \coordinate (t) at ($(m)!1/2!(h)$);
      \coordinate (r) at ($(j)!1/3!(o)$);
      \coordinate (s) at ($(j)!2/3!(o)$);
      \coordinate (u) at ($(a)!1/3!(f)$);
      \coordinate (v) at ($(a)!2/3!(f)$);
      \draw[blue,name path=B_1] (j)--(r);
      \draw[blue,name path=B_1'](s)--(o);
      \draw[red,name path=A_1](m)--(q);
      \draw[red,name path=A_1'](t)--(h);
      \draw[green,name path=E](q)--(r);
      \draw[black,name path=C_1](l)--(v);
      \draw[red]($(c)!1/30!(d)$)--($(b)!1/30!(a)$);
      \draw[black]($(e)!1/30!(d)$)--($(f)!1/30!(a)$);
      \node[black,above] at ($(b)!1/2!(c)$){$C_3$};
      \node[black,below] at ($(e)!1/2!(f)$){$C_0$};
      \node[black,above] at ($(v)!1/3!(l)$){$C_1$};
      \node[red,above] at ($(h)!1/4!(m)$){$A_1$};
      \coordinate (w) at ($(a)!1+2/3!(c)$);
      \coordinate (x) at ($(f)!1+2/3!(d)$);
      \coordinate (y) at ($(a)!1+1/2!(d)$);
      \draw[black](c)--(w);
      \draw[black](l)--($(v)!1+2/3!(l)$);
      \draw[black]($(x)!2/30!(w)$)--($(y)!1/30!(d)$);
      \draw[blue]($(c)!1/2!(w)$)--($(d)!2/3!(y)$);
      \draw[red](w)--(x);
      \draw[red]($(c)!1/30!(d)$)--($(x)!2/3!(y)$);
      \draw[red,ultra thick](d)--(y);
      \draw[blue,name path=B_3](x)--(y);
      \node[blue,left] at ($(x)!1/2!(y)$){$B_3$};
      \coordinate (z) at ($(f)!2!(e)$);
      \coordinate (y') at ($(y)!1/3!(z)$);
      \coordinate (z') at ($(z)!1/3!(e)$);
      \coordinate (d') at ($(y)!1/3!(d)$);
      \draw[blue](y)--(z);
      \draw[black](e)--(z);
      \draw[black,name path=C_2]($(y)!1/30!(d)$)--($(e)!1/30!(d)$);
      \draw[blue,name path=B_2]($(d')!2/3!(z')$)--(z');
      \draw[blue,name path=B_2']($(d')!1/5!(z')$)--(d');
      \draw[red,name path=A_1]($(y')!2/3!(m)$)--(m);
      \draw[red,name path=A_1']($(y')!1/5!(m)$)--(y');
      \draw[green,name path=E_2]($(d')!2/3!(z')$)--($(y')!2/3!(m)$);
      \node [red,above] at ($(y')!2/3!(m)$){$A_1$}; 
      \node [blue,left] at ($(d')!2/3!(z')$){$B_2$}; 
 \node at ($(a)+(0.3,0)$){$\rightsquigarrow$};
      \coordinate (x) at (4,-2.2);
      \coordinate (O) at ($(0,0)+(x)$);
      \coordinate (a) at ($(1/2,3^0.5/2)+(x)$);
      \coordinate (b) at ($(-1/2,3^0.5/2)+(x)$);
      \coordinate (c) at ($(-1,0)+(x)$);
      \coordinate (d) at ($(-1/2,-3^0.5/2)+(x)$);
      \coordinate (e) at ($(1/2,-3^0.5/2)+(x)$);
      \draw[blue,ultra thick] (O)--(a);
      \draw[blue] ($(O)!1/2!(e)$)--(a);
      \draw[blue] (a)--(e);
      \draw[black,ultra thick] (O)--(e);
      \draw[red] ($(O)!1/2!(a)$)--(e);
      \draw[red,ultra thick] (O)--(c);
      \draw[red] (a)--(c);
      \draw[red] (c)--(b);
      \draw[black] ($(O)!1/2!(a)$)--($(c)!1/2!(b)$);
      \draw[black] (a)--(b);
      \draw[blue] ($(O)!1/3!(c)$)--($(a)!1/3!(b)$); 
      \draw[blue] (c)--($(O)!1/2!(e)$);
      \draw[black] ($(O)!1/3!(c)$)--(e);
      \draw[black] (c)--(e);
      \node at ($(O)+(0,-1.3)$){$Case$ $1$};
      \coordinate (x) at (6.6,-2.2);
      \coordinate (O) at ($(0,0)+(x)$);
      \coordinate (a) at ($(1/2,3^0.5/2)+(x)$);
      \coordinate (b) at ($(-1/2,3^0.5/2)+(x)$);
      \coordinate (c) at ($(-1,0)+(x)$);
      \coordinate (d) at ($(-1/2,-3^0.5/2)+(x)$);
      \coordinate (e) at ($(1/2,-3^0.5/2)+(x)$);
     \node at ($(c)+(-0.3,0)$){$\rightsquigarrow$};
      \draw[blue,ultra thick] (O)--(a);
      \draw[blue] ($(O)!1/2!(e)$)--(a);
      \draw[blue] (a)--(e);
      \draw[black,ultra thick] (O)--(e);
      \draw[red] ($(O)!1/2!(a)$)--(e);
      \draw[red,ultra thick] (O)--(c);
      \draw[red] (a)--(c);
      \draw[red] (c)--(b);
      \draw[black] ($(O)!1/2!(a)$)--(b);
      \draw[black] (a)--(b);
      \draw[blue] ($(O)!1/3!(c)$)--(b); 
      \draw[black,ultra thick] (O)--(b);
      \draw[blue] (c)--($(O)!1/2!(e)$);
      \draw[black] ($(O)!1/3!(c)$)--(e);
      \draw[black] (c)--(e);
      \node at ($(O)+(0,-1.3)$){$Case$ $2$};

\end{tikzpicture}     
\par\end{center}

\textbf{Case 3.} When one of the two points $P_{1},P_{2}$ is on $B_{0}$
or $B_{3}$. Each degenerating arrangement is the same up to rotation.
WLOG, we can assume that $P_{2}$ is on $B_{0}$. To get the minimal
model, we first blow up the total space $\mathcal{Y}$ along the curve
$B_{0}$ in the central fiber. Let curves $\widetilde{L_{P_{1}}}$
and $\widetilde{L_{P_{2}}}$ be the proper transform of $L_{P_{1}}$
and $L_{P_{2}}$. Then blow up $\mathcal{Y}_{1}=Bl_{1}\mathcal{Y}$
along $\widetilde{L_{P_{1}}}$ and $\widetilde{L_{P_{2}}}$. The central
fiber $\Sigma_{0}$ becomes $Bl_{4}\mathbb{P}^{2}\cup Bl_{1}\mathbb{F}_{1}$. 

\begin{center}
\begin{tikzpicture}[xscale=1.5,yscale=1.5][font=\tiny]
      \coordinate (a) at (1,0);
      \coordinate (b) at (1/2,3^0.5/2);
      \coordinate (c) at (-1/2,3^0.5/2);
      \coordinate (d) at (-1,0);
      \coordinate (e) at (-1/2,-3^0.5/2);
      \coordinate (f) at (1/2,-3^0.5/2);
      \draw[blue](a)--(b);
      \draw[black](b)--(c);
      \draw[red](c)--(d);
      \draw[blue](d)--(e);
      \draw[black](e)--(f);
      \draw[red](f)--(a);
      \coordinate (g) at ($(a)!5/8!(b)$);
      \coordinate (h) at ($(b)!3/8!(c)$);
      \coordinate (i) at ($(b)!3/4!(c)$);
      \coordinate (j) at ($(c)!1/4!(d)$);
      \coordinate (k) at ($(c)!5/8!(d)$);
      \coordinate (l) at ($(d)!3/8!(e)$);
      \coordinate (m) at ($(e)!1/4!(f)$);
      \coordinate (n) at ($(e)!5/8!(f)$);
      \coordinate (o) at ($(f)!3/8!(a)$);
      \coordinate (p) at ($(f)!3/4!(a)$);
      \draw [name path=A_1,red](g)--(l);
      \draw [name path=A_2,red]($(a)!3/8!(b)$)--($(e)!3/8!(d)$);
      \draw [name path=C_1,black](k)--(o);
      \draw [name path=C_2,black]($(c)!1/30!(d)$)--($(b)!1/60!(e)$);
      \draw [name path=C_2',black]($(b)!1/60!(e)$)--($(a)!1/30!(f)$);
      \draw [name path=B_1,blue](i)--(m);
      \draw [name path=B_2,blue]($(b)!2/30!(c)$)--($(a)!2/60!(d)$);
      \draw [name path=B_2',blue]($(a)!2/60!(d)$)--($(f)!2/30!(e)$);
      \node[above,black]at($(b)!1/2!(c)$){$C_3$};
      \node[below,black]at($(e)!1/2!(f)$){$C_0$};
      \node[right,blue]at($(a)!1/2!(b)$){$B_0$};
      \node[left,blue]at($(d)!1/2!(e)$){$B_3$};
      \node[left,red]at($(c)!1/2!(d)$){$A_0$};
      \node[right,red]at($(a)!1/2!(f)$){$A_3$};
      \path[fill=green](g)circle[radius=0.03];
      \path[fill=green]($(a)!5/8!(d)$)circle[radius=0.03];
      \node [left] at (g){$P_1$};
      \node [left] at ($(a)!5/8!(d)$){$P_2$};

\node at ($(a)+(0.3,0)$){$\rightsquigarrow$};
      \coordinate (a) at (1+2.7,0);
      \coordinate (b) at (1/2+2.7,3^0.5/2);
      \coordinate (c) at (-1/2+2.7,3^0.5/2);
      \coordinate (d) at (-1+2.7,0);
      \coordinate (e) at (-1/2+2.7,-3^0.5/2);
      \coordinate (f) at (1/2+2.7,-3^0.5/2);
      \draw[black,ultra thick](a)--(b);
      \draw[black](b)--(c);
      \draw[red](c)--(d);
      \draw[blue](d)--(e);
      \draw[black](e)--(f);
      \draw[red](f)--(a);
      \coordinate (g) at ($(a)!5/8!(b)$);
      \coordinate (h) at ($(b)!3/8!(c)$);
      \coordinate (i) at ($(b)!3/4!(c)$);
      \coordinate (j) at ($(c)!1/4!(d)$);
      \coordinate (k) at ($(c)!5/8!(d)$);
      \coordinate (l) at ($(d)!3/8!(e)$);
      \coordinate (m) at ($(e)!1/4!(f)$);
      \coordinate (n) at ($(e)!5/8!(f)$);
      \coordinate (o) at ($(f)!3/8!(a)$);
      \coordinate (p) at ($(f)!3/4!(a)$);
      \draw [name path=A_1,red](g)--(l);
      \draw [name path=A_2,red]($(a)!3/8!(b)$)--($(e)!3/8!(d)$);
      \draw [name path=C_1,black](k)--(o);
      \draw [name path=C_2,black]($(c)!1/30!(d)$)--($(b)!1/30!(a)$);
      \draw [name path=B_1,blue](i)--(m);
      \draw [name path=B_2,blue]($(a)!1/30!(b)$)--($(f)!1/30!(e)$);
      \path[fill=green]($(a)!5/8!(d)$)circle[radius=0.03];
      \coordinate (r) at ($(d)!1+2/3!(b)$);
      \coordinate (s) at ($(e)!1+2/3!(a)$);
      \coordinate (t) at ($(r)!3/8!(s)$);
      \draw [name path=C_3,black](b)--(r);
      \draw [name path=B_0,blue](r)--(s);
      \draw [name path=A_3,red](a)--(s);
      \draw [name path=A_1,red](g)--($(l)!1+2/3!(g)$);
      \draw [name path=A_2,red]($(a)!3/8!(b)$)--($(s)!3/8!(r)$);
      \draw [name path=C_2,black]($(b)!1/30!(a)$)--($(a)!2/3!(s)$);
      \draw [name path=B_2,blue]($(a)!1/30!(b)$)--($(b)!2/5!(r)$);
      \path[fill=green]($(g)!1/4!(t)$)circle[radius=0.03];
\node at ($(a)+(0.5,0)$){$\rightsquigarrow$};
      \coordinate (a) at (1+5.7,0);
      \coordinate (b) at (1/2+5.7,3^0.5/2);
      \coordinate (c) at (-1/2+5.7,3^0.5/2);
      \coordinate (d) at (-1+5.7,0);
      \coordinate (e) at (-1/2+5.7,-3^0.5/2);
      \coordinate (f) at (1/2+5.7,-3^0.5/2);
      \draw[black,ultra thick](a)--(b);
      \draw[black](b)--(c);
      \draw[red](c)--(d);
      \draw[blue](d)--(e);
      \draw[black](e)--(f);
      \draw[red](f)--(a);
      \coordinate (g) at ($(a)!5/8!(b)$);
      \coordinate (h) at ($(b)!3/8!(c)$);
      \coordinate (i) at ($(b)!3/4!(c)$);
      \coordinate (j) at ($(c)!1/4!(d)$);
      \coordinate (k) at ($(c)!5/8!(d)$);
      \coordinate (l) at ($(d)!3/8!(e)$);
      \coordinate (m) at ($(e)!1/4!(f)$);
      \coordinate (n) at ($(e)!5/8!(f)$);
      \coordinate (o) at ($(f)!3/8!(a)$);
      \coordinate (p) at ($(f)!3/4!(a)$);
      \draw [name path=A_1,red](g)--(l);
      \draw [name path=A_2,red]($(a)!3/8!(b)$)--($(e)!3/8!(d)$);
      \draw [name path=C_1,black](k)--($(k)!1/4!(o)$);
      \draw [name path=C_1',black](o)--($(o)!1/2!(k)$);
      \draw [name path=C_2,black]($(c)!1/30!(d)$)--($(b)!1/30!(a)$);
      \draw [name path=B_1,blue](i)--($(i)!1/3!(m)$);
      \draw [name path=B_1',blue](m)--($(m)!5/12!(i)$);
      \draw [name path=B_2,blue]($(a)!1/30!(b)$)--($(f)!1/30!(e)$);
      \draw [name path=E_1,green]($(k)!1/4!(o)$)--($(m)!5/12!(i)$);
      \node[above,black]at($(b)!1/2!(c)$){$C_3$};
      \node[right,red]at($(a)!1/2!(f)$){$A_3$};
      \node[above,red]at($(g)!1/3!(l)$){$A_1$};
      \coordinate (r) at ($(d)!1+2/3!(b)$);
      \coordinate (s) at ($(e)!1+2/3!(a)$);
      \coordinate (t) at ($(r)!3/8!(s)$);
      \coordinate (u) at ($(a)!2/3!(s)$);
      \coordinate (v) at ($(b)!2/5!(r)$);
      \draw [name path=C_3,black](b)--(r);
      \draw [name path=B_0,blue](r)--(s);
      \draw [name path=A_3,red](a)--(s);
      \draw [name path=A_1,red](g)--(t);
      \draw [name path=A_2,red]($(a)!3/8!(b)$)--($(s)!3/8!(r)$);
      \draw [name path=C_2,black]($(b)!1/30!(a)$)--($(b)!3/10!(u)$);
      \draw [name path=C_2',black](u)--($(b)!1/2!(u)$);
      \draw [name path=B_2,blue]($(a)!1/30!(b)$)--($(a)!1/2!(v)$);
      \draw [name path=B_2',blue](v)--($(v)!1/4!(a)$);
      \draw [name path=E_2,green]($(b)!3/10!(u)$)--($(a)!1/2!(v)$);
      \node [above,red]at($(g)!1/2!(t)$){$A_1$};
\node at($(m)+(0.3,-0.3)$){$Bl_{4}\mathbb{P}^{2}\cup Bl_{1}\mathbb{F}_{1}$};
\end{tikzpicture}
\par\end{center}

The canonical model of $\Sigma_{0}$ is $\Sigma_{0}^{c}=\mathbb{F}_{0}\cup\mathbb{F}_{0}$
which we denote as Case 3. The first figure below is obtained from
the component $Bl_{4}\mathbb{P}^{2}$ by contracting 3 curves $A_{1},A_{3},C_{3}$.
The second figure below is obtained from the component $Bl_{1}\mathbb{F}_{1}$
by contracting the curve $A_{1}$. This case could be obtained from
Case 4 for $K^{2}=5$. 

\begin{center}
\begin{tikzpicture}[xscale=2,yscale=2][font=\tiny]
      \coordinate (e) at (0,0);\node [below] at (e){1};
      \coordinate (f) at (0,1);\node [above] at (f){4};
      \coordinate (g) at (-3^0.5/2,1/2);
      \coordinate (h) at (3^0.5/2,1/2);
      \draw [name path=triangle,ultra thick,black] (e)--(f);
      \draw [name path=B_1,blue] (g)--(f); 
      \draw [name path=A_0,red] (f)--(h);    
      \draw [name path=C_0,black] (g)--(e);
      \draw [name path=C_1,black] (e)--(h);
      \draw [name path=A_2,red] ($(g)!1/3!(f)$)--($(e)!1/3!(h)$);
      \draw [name path=B_3,blue] ($(e)!2/3!(g)$)--($(h)!2/3!(f)$);
      \draw [name path=E,green]($(f)!1/3!(g)$)--($(h)!1/3!(e)$);
      \node [black,right] at ($(e)!1/3!(f)$){2};
      \node [black,right] at ($(e)!2/3!(f)$){3};
      \node[black] at ($(e)+(0,-1/4)$){$\mathbb{P}^{1}\times\mathbb{P}^{1}$};
\node[black] at ($(h)+(0.4,0)$){$\cup$};
\node[black] at ($(h)+(0.4,-1)$){$Case$ $3$};

      \coordinate (e) at (0+2.5,0);\node [below] at (e){1};
      \coordinate (f) at (0+2.5,1);\node [above] at (f){4};
      \coordinate (g) at (-3^0.5/2+2.5,1/2);
      \coordinate (h) at (3^0.5/2+2.5,1/2);
      \draw [name path=triangle,ultra thick,black] (e)--(f);
      \draw [name path=C_2,black] (g)--(f); 
      \draw [name path=C_3,black] (f)--(h);    
      \draw [name path=A_3,red] (g)--(e);
      \draw [name path=B_2,blue] (e)--(h);
      \draw [name path=A_2,red] ($(g)!1/3!(f)$)--($(e)!1/3!(h)$);
      \draw [name path=B_0,blue] ($(e)!2/3!(g)$)--($(h)!2/3!(f)$);
      \draw [name path=E_2,green] ($(f)!1/3!(g)$)--($(h)!1/3!(e)$);
      \node [black,right] at ($(e)!1/3!(f)$){2};
      \node [black,right] at ($(e)!2/3!(f)$){3};           
      \node[black] at ($(e)+(0,-1/4)$){$\mathbb{P}^{1}\times\mathbb{P}^{1}$};
\end{tikzpicture}
\par\end{center}

There are further degenerations; however, the further degenerations
do not produce any new cases. For example, when the point 3 on the
double locus goes to the point 4 in the above figures, the lc model
of the further degeneration is the same as Case 2. \\

\textbf{Case 4.} When $P_{1},P_{2}$ coincide. Blow up the total space
$\mathcal{Y}$ at the point $P$ in the central fiber, the central
fiber is as follows

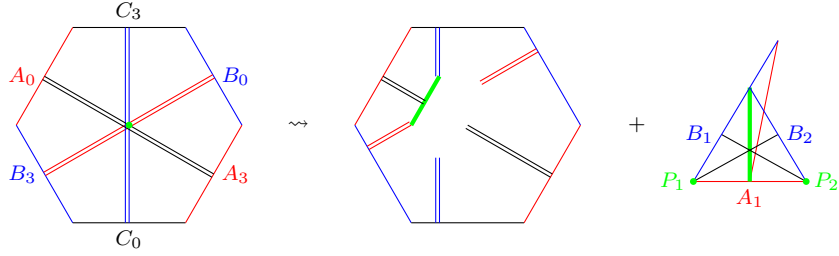
\begin{center}
\begin{tikzpicture}[xscale=1.5,yscale=1.5][font=\tiny]
      \coordinate (a) at (1,0);
      \coordinate (b) at (1/2,3^0.5/2);
      \coordinate (c) at (-1/2,3^0.5/2);
      \coordinate (d) at (-1,0);
      \coordinate (e) at (-1/2,-3^0.5/2);
      \coordinate (f) at (1/2,-3^0.5/2);
      \draw[blue](a)--(b);
      \draw[black](b)--(c);
      \draw[red](c)--(d);
      \draw[blue](d)--(e);
      \draw[black](e)--(f);
      \draw[red](f)--(a);
      \coordinate (g) at ($(a)!1/2-1/60!(b)$);
      \coordinate (h) at ($(a)!1/2+1/60!(b)$);
      \coordinate (i) at ($(b)!1/2!(c)$);
      \coordinate (i') at ($(c)!1/2-1/30!(b)$);
      \coordinate (j) at ($(c)!1/2!(d)$);
      \coordinate (j') at ($(c)!1/2+1/30!(d)$);
      \coordinate (k) at ($(d)!1/2-1/60!(e)$);
      \coordinate (l) at ($(d)!1/2+1/60!(e)$);
      \coordinate (m) at ($(e)!1/2!(f)$);
      \coordinate (m') at ($(e)!1/2-1/30!(f)$);
      \coordinate (n) at ($(f)!1/2!(a)$);
      \coordinate (n') at ($(f)!1/2-1/30!(a)$);
      \draw [name path=A_1,red](h)--(k);
      \draw [name path=A_2,red](g)--(l);
      \draw [name path=B_1,blue](i)--(m);
      \draw [name path=B_2,blue](i')--(m');
      \draw [name path=C_1,black](j)--(n);
      \draw [name path=C_2,black](j')--(n');
      \node[above,black]at($(b)!1/2!(c)$){$C_3$};
      \node[below,black]at($(e)!1/2!(f)$){$C_0$};
      \node[right,blue]at($(a)!1/2!(b)$){$B_0$};
      \node[left,blue]at($(d)!1/2!(e)$){$B_3$};
      \node[left,red]at($(c)!1/2!(d)$){$A_0$};
      \node[right,red]at($(a)!1/2!(f)$){$A_3$};
      \path[fill=green]($(a)!1/2!(d)$)circle[radius=0.03];
\node at ($(a)+(0.5,0)$) {$\rightsquigarrow$};

      \coordinate (a) at (1+3,0);
      \coordinate (b) at (1/2+3,3^0.5/2);
      \coordinate (c) at (-1/2+3,3^0.5/2);
      \coordinate (d) at (-1+3,0);
      \coordinate (e) at (-1/2+3,-3^0.5/2);
      \coordinate (f) at (1/2+3,-3^0.5/2);
      \draw[blue](a)--(b);
      \draw[black](b)--(c);
      \draw[red](c)--(d);
      \draw[blue](d)--(e);
      \draw[black](e)--(f);
      \draw[red](f)--(a);
      \coordinate (g) at ($(a)!3/4+1/30!(b)$);
      \coordinate (h) at ($(a)!3/4!(b)$);
      \coordinate (i) at ($(b)!3/4!(c)$);
      \coordinate (i') at ($(c)!1/4-1/30!(b)$);
      \coordinate (j) at ($(c)!1/2!(d)$);
      \coordinate (j') at ($(c)!1/2+1/30!(d)$);
      \coordinate (k) at ($(d)!1/4!(e)$);
      \coordinate (l) at ($(d)!1/4-1/30!(e)$);
      \coordinate (m) at ($(e)!1/4!(f)$);
      \coordinate (m') at ($(e)!1/4-1/30!(f)$);
      \coordinate (n) at ($(f)!1/2!(a)$);
      \coordinate (n') at ($(f)!1/2-1/30!(a)$);
      \coordinate (o) at ($(c)!1/2!(f)$);
      \draw [name path=A_1,red](h)--($(h)!1/3!(k)$);
      \draw [name path=A_1',red](k)--($(o)!1/2!(d)$);
      \draw [name path=A_2,red](g)--($(g)!1/3!(l)$);
      \draw [name path=A_2',red](l)--($(l)!1/4!(g)$);
      \draw [name path=E,green,ultra thick]($(o)!1/2!(c)$)--($(o)!1/2!(d)$);
      \draw [name path=B_1,blue](i)--($(i)!1/4!(m)$);
      \draw [name path=B_1',blue](m)--($(m)!1/3!(i)$);
      \draw [name path=B_2,blue](i')--($(i')!1/4!(m')$);
      \draw [name path=B_2',blue](m')--($(m')!1/3!(i')$);
      \draw [name path=C_1,black](j)--($(j)!1/4!(n)$);
      \draw [name path=C_1',black](o)--(n);
      \draw [name path=C_2,black](j')--($(j')!1/4!(n')$);
      \draw [name path=C_2',black]($(e)!1/2-1/60!(b)$)--(n');
\node at ($(a)+(0.5,0)$) {$+$};
      \coordinate (a) at (2+3,-1/2);
      \coordinate (b) at (3+3,-1/2);
      \coordinate (c) at (2.5+3,1/3);
      \coordinate (d) at (2+3,-1/2);
      \coordinate (e) at (3+3,-1/2);
      \draw[red](a)--(b);
      \draw[red]($(a)!1/2!(b)$)--($(d)!1+1/2!(c)$);
      \draw[green,ultra thick](c)--($(a)!1/2!(b)$);
      \draw[black]($(a)!1/2!(d)$)--($(e)!1/2!(c)$);
      \draw[black]($(b)!1/2!(e)$)--($(d)!1/2!(c)$);
      \draw[blue](d)--($(d)!1+1/2!(c)$);
      \draw[blue](e)--(c);
      \draw[green](a)--(d);
      \draw[green](b)--(e);
      \node[left,blue]at($(d)!1/2!(c)$){$B_1$};
      \node[right,blue]at($(e)!1/2!(c)$){$B_2$};
      \node[below,red]at($(a)!1/2!(b)$){$A_1$};
      \node[left,green]at($(a)!1/2!(d)$){$P_1$};
      \node[right,green]at($(b)!1/2!(e)$){$P_2$};
      \path[fill=green]($(a)!1/2!(d)$)circle[radius=0.03];
      \path[fill=green]($(b)!1/2!(e)$)circle[radius=0.03];
\end{tikzpicture}
\par\end{center}

Then we blow up the total space along the proper transform $\tilde{L}_{P_{1}}$
and $\tilde{L}_{P_{2}}$. The central fiber of the minimal model is
$Bl_{4}\mathbb{P}^{2}\cup Bl_{2}\mathbb{P}^{2}$. Running the minimal
model program, we contract $A_{1},A_{2},B_{1},B_{2},C_{1},C_{2}$
in the component $Bl_{4}\mathbb{P}^{2}$ and $ $$A_{1}$ in the component
$Bl_{2}\mathbb{P}^{2}$ in the central fiber. The central fiber of
the canonical model is $ $$\mathbb{F}_{0}\cup\mathbb{F}_{0}$ as
follows 

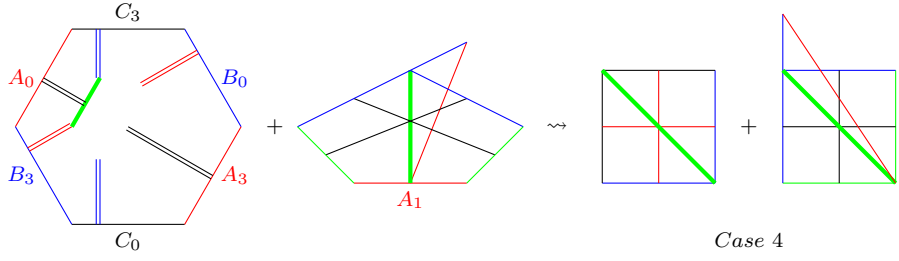
\begin{center}
\begin{tikzpicture}[xscale=1.5,yscale=1.5][font=\tiny]
      \coordinate (x) at (-3,0);
      \coordinate (a) at ($(1,0)+(x)$);
      \coordinate (b) at ($(1/2,3^0.5/2)+(x)$);
      \coordinate (c) at ($(-1/2,3^0.5/2)+(x)$);
      \coordinate (d) at ($(-1,0)+(x)$);
      \coordinate (e) at ($(-1/2,-3^0.5/2)+(x)$);
      \coordinate (f) at ($(1/2,-3^0.5/2)+(x)$);
      \draw[blue](a)--(b);
      \draw[black](b)--(c);
      \draw[red](c)--(d);
      \draw[blue](d)--(e);
      \draw[black](e)--(f);
      \draw[red](f)--(a);
      \coordinate (g) at ($(a)!3/4+1/30!(b)$);
      \coordinate (h) at ($(a)!3/4!(b)$);
      \coordinate (i) at ($(b)!3/4!(c)$);
      \coordinate (i') at ($(c)!1/4-1/30!(b)$);
      \coordinate (j) at ($(c)!1/2!(d)$);
      \coordinate (j') at ($(c)!1/2+1/30!(d)$);
      \coordinate (k) at ($(d)!1/4!(e)$);
      \coordinate (l) at ($(d)!1/4-1/30!(e)$);
      \coordinate (m) at ($(e)!1/4!(f)$);
      \coordinate (m') at ($(e)!1/4-1/30!(f)$);
      \coordinate (n) at ($(f)!1/2!(a)$);
      \coordinate (n') at ($(f)!1/2-1/30!(a)$);
      \coordinate (o) at ($(c)!1/2!(f)$);
      \draw [name path=A_1,red](h)--($(h)!1/3!(k)$);
      \draw [name path=A_1',red](k)--($(o)!1/2!(d)$);
      \draw [name path=A_2,red](g)--($(g)!1/3!(l)$);
      \draw [name path=A_2',red](l)--($(l)!1/4!(g)$);
      \draw [name path=E,green,ultra thick]($(o)!1/2!(c)$)--($(o)!1/2!(d)$);
      \draw [name path=B_1,blue](i)--($(i)!1/4!(m)$);
      \draw [name path=B_1',blue](m)--($(m)!1/3!(i)$);
      \draw [name path=B_2,blue](i')--($(i')!1/4!(m')$);
      \draw [name path=B_2',blue](m')--($(m')!1/3!(i')$);
      \draw [name path=C_1,black](j)--($(j)!1/4!(n)$);
      \draw [name path=C_1',black](o)--(n);
      \draw [name path=C_2,black](j')--($(j')!1/4!(n')$);
      \draw [name path=C_2',black]($(e)!1/2-1/60!(b)$)--(n');
      \node[above,black]at($(b)!1/2!(c)$){$C_3$};
      \node[below,black]at($(e)!1/2!(f)$){$C_0$};
      \node[right,blue]at($(a)!1/2!(b)$){$B_0$};
      \node[left,blue]at($(d)!1/2!(e)$){$B_3$};
      \node[left,red]at($(c)!1/2!(d)$){$A_0$};
      \node[right,red]at($(a)!1/2!(f)$){$A_3$};
\node at ($(a)+(0.3,0)$) {$+$};
      \coordinate (x) at (-3,0);
      \coordinate (a) at ($(2,-1/2)+(x)$);
      \coordinate (b) at ($(3,-1/2)+(x)$);
      \coordinate (c) at ($(2.5,1/2)+(x)$);
      \coordinate (d) at ($(1.5,0)+(x)$);
      \coordinate (e) at ($(3.5,0)+(x)$);
      \draw[red](a)--(b);
      \draw[red]($(a)!1/2!(b)$)--($(d)!1+1/2!(c)$);
      \draw[green,ultra thick](c)--($(a)!1/2!(b)$);
      \draw[black]($(a)!1/2!(d)$)--($(e)!1/2!(c)$);
      \draw[black]($(b)!1/2!(e)$)--($(d)!1/2!(c)$);
      \draw[blue](d)--($(d)!1+1/2!(c)$);
      \draw[blue](e)--(c);
      \draw[green](a)--(d);
      \draw[green](b)--(e);

      \node[below,red]at($(a)!1/2!(b)$){$A_1$};
  \coordinate (x) at (1.2,-0.5);
  \coordinate (a) at ($(0,0)+(x)$);
  \coordinate (b) at ($(1,0)+(x)$);
  \coordinate (c) at ($(1,1)+(x)$);
  \coordinate (d) at ($(0,1)+(x)$);
  \draw [name path=B_3,blue](a)--(b);
  \draw [name path=B_0,blue](b)--(c);
  \draw [name path=C_3,black](c)--(d);
  \draw [name path=C_0,black](d)--(a);
  \draw [name path=A_0,red]($(a)!1/2!(b)$)--($(c)!1/2!(d)$);
  \draw [name path=A_3,red]($(b)!1/2!(c)$)--($(d)!1/2!(a)$);
  \draw [name path=triangle,green,ultra thick](b)--(d);
  \node at ($(e)+(0.3,0)$){$\rightsquigarrow$};
\node at ($(1.3,0.5)+(x)$) {$+$};
\node at ($(1.3,0.5)+(0,-1)+(x)$){$Case$ $4$};
  \coordinate (x) at (0.8,-0.5);
  \coordinate (a) at ($(0+2,0)+(x)$);
  \coordinate (b) at ($(1+2,0)+(x)$);
  \coordinate (c) at ($(1+2,1)+(x)$);
  \coordinate (d) at ($(0+2,1)+(x)$); 
  \draw [name path=E_1,green](a)--(b);
  \draw [name path=E_2,green](b)--(c);
  \draw [name path=B_2,blue](c)--(d);
  \draw [name path=B_1,blue](a)--($(a)!1+1/2!(d)$);
  \draw [name path=C_1,black]($(a)!1/2!(b)$)--($(c)!1/2!(d)$);
  \draw [name path=C_2,black]($(b)!1/2!(c)$)--($(d)!1/2!(a)$);
  \draw [name path=triangle,green,ultra thick](b)--(d);
  \draw [name path=A_2,red](b)--($(a)!1+1/2!(d)$);

\end{tikzpicture}
\par\end{center}

\textbf{Case 5. }The further degeneration of Case 4 above. After blowing
up the total space at $P$, the points $P_{1},P_{2}$ could still
coincide in the exceptional divisor $\mathbb{P}^{2}$ of the blowup.
We need to blow up the total space at the point $P$ first, then $P_{1},P_{2}$
will be distinct. Now we can blow up the total space along the lines
$\tilde{L}_{P_{1}}$ and $\tilde{L}_{P_{2}}$ . The following figures
are only the second component of the lc model, with the first component
$Bl_{4}\mathbb{P}^{2}$, which is the same as the first figure above.

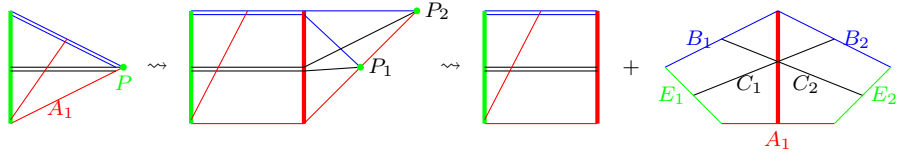
\begin{center}
\begin{tikzpicture}[xscale=1.5,yscale=1.5][font=\tiny]
      \coordinate (x) at (0,0);
      \coordinate (a) at ($(0,-1/2)+(x)$);
      \coordinate (b) at ($(1,0)+(x)$);
      \coordinate (c) at ($(0,1/2)+(x)$);
      \coordinate (o) at ($(a)!1/2!(c)$);
      \draw[red,name path=A_1](a)--(b);
      \draw[blue,name path=B_1](b)--(c);            
      \draw[blue,name path=B_2]($(c)+(0,-1/30)$)--($(b)+(-2/30,0)$);
      \draw[green,ultra thick](a)--(c);
      \draw[black,name path=C_1](o)--(b);
      \draw[black,name path=C_2]($(o)+(0,-1/30)$)--($(b)+(-2/30,-1/30)$);
      \draw[red,name path=A_2]($(a)+(0,1/30)$)--($(b)!1/2!(c)$);
      \node[right,red]at($(a)!1/4!(b)$){$A_1$};
      \node[below,green]at(b){$P$};
      \path[fill=green](b)circle[radius=0.03];
      \node at ($(b)+(0.3,0)$){$\rightsquigarrow$};

      \coordinate (x) at (1.6,0);
      \coordinate (a) at ($(0,-1/2)+(x)$);
      \coordinate (b) at ($(1,-1/2)+(x)$);
      \coordinate (c) at ($(0,1/2)+(x)$);
      \coordinate (d) at ($(1,1/2)+(x)$);
      \coordinate (o) at ($(a)!1/2!(c)$);
      \coordinate (f) at ($(b)!1/2!(d)$);
      \draw[red,name path=A_1](a)--(b);
      \draw[blue,name path=B_1](c)--(d);            
      \draw[blue,name path=B_2]($(c)+(0,-1/30)$)--($(d)+(0,-1/30)$);
      \draw[green,ultra thick](a)--(c);
      \draw[red,ultra thick](b)--(d);
      \draw[black,name path=C_1](o)--(f);
      \draw[black,name path=C_2]($(o)+(0,-1/30)$)--($(f)+(0,-1/30)$);
      \draw[red,name path=A_2]($(a)+(0,1/30)$)--($(d)!1/2!(c)$);
      \coordinate (e) at ($(c)!2!(d)$);
      \draw[blue,name path=B_2](d)--(e);
      \draw[blue,name path=B_1]($(d)+(0,-1/30)$)--($(b)!1/2!(e)$);
      \draw[black,name path=C_1](f)--(e);
      \draw[black,name path=C_2]($(f)+(0,-1/30)$)--($(b)!1/2!(e)$);
      \draw[red,name path=A_1](b)--(e);
      \path[fill=green](e)circle[radius=0.03];
      \path[fill=green]($(b)!1/2!(e)$)circle[radius=0.03];
      \node[right] at (e){$P_2$};
      \node[right] at ($(b)!1/2!(e)$){$P_1$};
   \node at ($(f)+(1.3,0)$){$\rightsquigarrow$}; 
      \coordinate (x) at (4.2,0);
      \coordinate (a) at ($(0,-1/2)+(x)$);
      \coordinate (b) at ($(1,-1/2)+(x)$);
      \coordinate (c) at ($(0,1/2)+(x)$);
      \coordinate (d) at ($(1,1/2)+(x)$);
      \coordinate (o) at ($(a)!1/2!(c)$);
      \coordinate (f) at ($(b)!1/2!(d)$);
      \draw[red,name path=A_1](a)--(b);
      \draw[blue,name path=B_1](c)--(d);            
      \draw[blue,name path=B_2]($(c)+(0,-1/30)$)--($(d)+(0,-1/30)$);
      \draw[green,ultra thick](a)--(c);
      \draw[red,ultra thick](b)--(d);
      \draw[black,name path=C_1](o)--(f);
      \draw[black,name path=C_2]($(o)+(0,-1/30)$)--($(f)+(0,-1/30)$);
      \draw[red,name path=A_2]($(a)+(0,1/30)$)--($(d)!1/2!(c)$);
   \node at ($(f)+(0.3,0)$){$+$};
      \coordinate (x) at (4.3,0);
      \coordinate (a) at ($(2,-1/2)+(x)$);
      \coordinate (b) at ($(3,-1/2)+(x)$);
      \coordinate (c) at ($(2.5,1/2)+(x)$);
      \coordinate (d) at ($(1.5,0)+(x)$);
      \coordinate (e) at ($(3.5,0)+(x)$);
      \draw[red](a)--(b);
      \draw[red,ultra thick](c)--($(a)!1/2!(b)$);
      \draw[black]($(a)!1/2!(d)$)--($(e)!1/2!(c)$);
      \draw[black]($(b)!1/2!(e)$)--($(d)!1/2!(c)$);
      \draw[blue](d)--(c);
      \draw[blue](e)--(c);
      \draw[green](a)--(d);
      \draw[green](b)--(e);
      \node[below,red]at($(a)!1/2!(b)$){$A_1$};
      \node[left,blue]at ($(d)!1/2!(c)$){$B_1$};
      \node[right,blue]at ($(e)!1/2!(c)$){$B_2$};
      \node[below,black]at ($(a)!1/2!(c)$){$C_1$};
      \node[below,black]at ($(b)!1/2!(c)$){$C_2$};
      \node[left,green]at ($(a)!1/2!(d)$){$E_1$};
      \node[right,green]at ($(b)!1/2!(e)$){$E_2$};
\end{tikzpicture}
\par\end{center}

Consider the line $A_{1}$ in the component $Bl_{2}\mathbb{P}^{2}$
of the central fiber, we have

\begin{eqnarray*}
(K_{\mathcal{Y}}+\mathcal{D}).A_{1} & = & -\dfrac{1}{2}<0
\end{eqnarray*}
and

\begin{eqnarray*}
K_{\Sigma_{0}}|_{Bl_{2}\mathbb{P}^{2}}.A_{1} & = & \left(K_{Bl_{2}\mathbb{P}^{2}}+\triangle\right).A_{1}=0
\end{eqnarray*}
According to the minimal model program, there will be a flip for $\left(\mathcal{Y},\dfrac{1}{2}\mathcal{D}\right)$.
The normal bundle of $A_{1}$ in the total space $\mathcal{Y}_{2}$
is $\mathcal{O}(-1)\oplus\mathcal{O}(-1)$. The flip for $\left(\mathcal{Y},\dfrac{1}{2}\mathcal{D}\right)$
is the Atiyah flop for $\mathcal{Y}$.

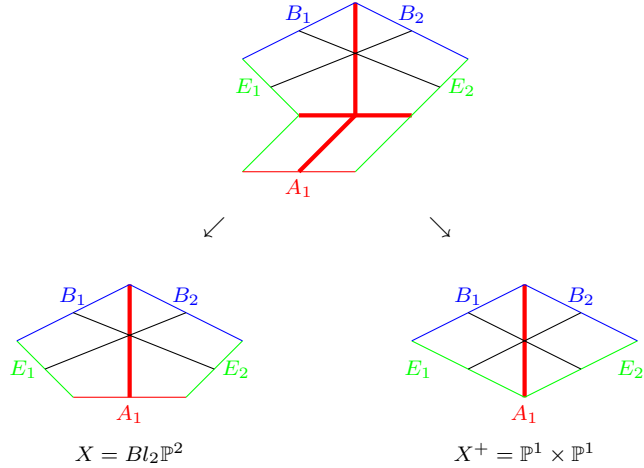
\begin{center}
\begin{tikzpicture}[xscale=1.5,yscale=1.5][font=\tiny]
      \coordinate (a) at (2.5-3,0);
      \coordinate (b) at (3.5-3,0);
      \coordinate (c) at (3-3,1);
      \coordinate (d) at (2-3,1/2);
      \coordinate (e) at (4-3,1/2);
      \draw[name path=A_1,red,ultra thick](a)--(b);
      \draw[red,ultra thick](c)--($(a)!1/2!(b)$);
      \draw[name path=E_1,green](b)--(e);
      \draw[name path=E_2,green](a)--(d);
      \draw[name path=B_1,blue](e)--(c);
      \draw[name path=B_2,blue](d)--(c);
      \draw[name path=C_1,black]($(a)!1/2!(d)$)--($(e)!1/2!(c)$);
      \draw[name path=C_2,black]($(b)!1/2!(e)$)--($(c)!1/2!(d)$);
      \node[above,blue]at($(d)!1/2!(c)$){$B_1$};
      \node[above,blue]at($(e)!1/2!(c)$){$B_2$};
      \node[left,green]at($(a)!1/2!(d)$){$E_1$};
      \node[right,green]at($(b)!1/2!(e)$){$E_2$};
      \coordinate (f) at (2-3,-1/2);
      \coordinate (g) at (3-3,-1/2);
      \draw[green](a)--(f);
      \draw[green](b)--(g);
      \draw[red,ultra thick]($(a)!1/2!(b)$)--($(f)!1/2!(g)$);
      \draw[red](f)--(g);
      \node[below,red]at($(f)!1/2!(g)$){$A_1$};
      \coordinate (a) at (2.5-5,-1/2-2);
      \coordinate (b) at (3.5-5,-1/2-2);
      \coordinate (c) at (3-5,1/2-2);
      \coordinate (d) at (2-5,0-2);
      \coordinate (e) at (4-5,0-2);
      \coordinate (f) at ($(a)!1/2!(b)$);
      \draw[name path=A_1,red](a)--(b);
      \draw[red,ultra thick](c)--($(a)!1/2!(b)$);
      \draw[name path=E_1,green](b)--(e);
      \draw[name path=E_2,green](a)--(d);
      \draw[name path=B_1,blue](e)--(c);
      \draw[name path=B_2,blue](d)--(c);
      \draw[name path=C_1,black]($(a)!1/2!(d)$)--($(e)!1/2!(c)$);
      \draw[name path=C_2,black]($(b)!1/2!(e)$)--($(c)!1/2!(d)$);
      \node[above,blue]at($(d)!1/2!(c)$){$B_1$};
      \node[above,blue]at($(e)!1/2!(c)$){$B_2$};
      \node[below,red]at($(a)!1/2!(b)$){$A_1$};
      \node[left,green]at($(a)!1/2!(d)$){$E_1$};
      \node[right,green]at($(b)!1/2!(e)$){$E_2$};
  \node[below,black] at ($(c)!1+1/3!(f)$){$X=Bl_{2}\mathbb{P}^{2}$};
      \coordinate (a) at (2.5-1.5,-1/2-2);
      \coordinate (b) at (3.5-1.5,-1/2-2);
      \coordinate (c) at (3-1.5,1/2-2);
      \coordinate (d) at (2-1.5,0-2);
      \coordinate (e) at (4-1.5,0-2);
      \coordinate (f) at ($(a)!1/2!(b)$);
      \draw[red,ultra thick](c)--($(a)!1/2!(b)$);
      \draw[name path=E_1,green](d)--(f);
      \draw[name path=E_2,green](e)--(f);
      \draw[name path=B_1,blue](e)--(c);
      \draw[name path=B_2,blue](d)--(c);
      \draw[name path=C_1,black]($(d)!1/2!(f)$)--($(e)!1/2!(c)$);
      \draw[name path=C_2,black]($(f)!1/2!(e)$)--($(c)!1/2!(d)$);
      \node[above,blue]at($(d)!1/2!(c)$){$B_1$};
      \node[above,blue]at($(e)!1/2!(c)$){$B_2$};
      \node[below,red]at($(a)!1/2!(b)$){$A_1$};
      \node[left,green]at($(a)!1/2!(d)$){$E_1$};
      \node[right,green]at($(b)!1/2!(e)$){$E_2$};
  \node[below,black] at ($(c)!1+1/3!(f)$){$X^{+}=\mathbb{P}^{1}\times\mathbb{P}^{1}$};
      \node at ($(-2.5,-2)!1/2!(0,0)$){$\swarrow$};
      \node at ($(1.5,-2)!1/2!(0,0)$){$\searrow$};
\end{tikzpicture}\\

\par\end{center}

When we apply the flip, the component $\mathbb{F}_{1}$ is blown up
at one point on the double locus. The central fiber becomes a union
of three components $Bl_{4}\mathbb{P}^{2}\cup Bl_{1}\mathbb{F}_{1}\cup\mathbb{F}_{0}$. 

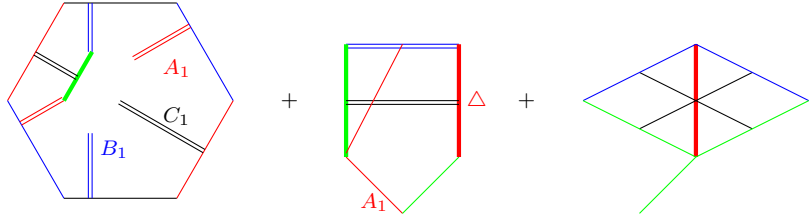
\begin{center}
\begin{tikzpicture}[xscale=1.5,yscale=1.5][font=\tiny]
      \coordinate (x) at (-3,0);
      \coordinate (a) at ($(1,0)+(x)$);
      \coordinate (b) at ($(1/2,3^0.5/2)+(x)$);
      \coordinate (c) at ($(-1/2,3^0.5/2)+(x)$);
      \coordinate (d) at ($(-1,0)+(x)$);
      \coordinate (e) at ($(-1/2,-3^0.5/2)+(x)$);
      \coordinate (f) at ($(1/2,-3^0.5/2)+(x)$);
      \draw[blue](a)--(b);
      \draw[black](b)--(c);
      \draw[red](c)--(d);
      \draw[blue](d)--(e);
      \draw[black](e)--(f);
      \draw[red](f)--(a);
      \coordinate (g) at ($(a)!3/4+1/30!(b)$);
      \coordinate (h) at ($(a)!3/4!(b)$);
      \coordinate (i) at ($(b)!3/4!(c)$);
      \coordinate (i') at ($(c)!1/4-1/30!(b)$);
      \coordinate (j) at ($(c)!1/2!(d)$);
      \coordinate (j') at ($(c)!1/2+1/30!(d)$);
      \coordinate (k) at ($(d)!1/4!(e)$);
      \coordinate (l) at ($(d)!1/4-1/30!(e)$);
      \coordinate (m) at ($(e)!1/4!(f)$);
      \coordinate (m') at ($(e)!1/4-1/30!(f)$);
      \coordinate (n) at ($(f)!1/2!(a)$);
      \coordinate (n') at ($(f)!1/2-1/30!(a)$);
      \coordinate (o) at ($(c)!1/2!(f)$);
      \draw [name path=A_1,red](h)--($(h)!1/3!(k)$);
      \draw [name path=A_1',red](k)--($(o)!1/2!(d)$);
      \draw [name path=A_2,red](g)--($(g)!1/3!(l)$);
      \draw [name path=A_2',red](l)--($(l)!1/4!(g)$);
      \draw [name path=E,green,ultra thick]($(o)!1/2!(c)$)--($(o)!1/2!(d)$);
      \draw [name path=B_1,blue](i)--($(i)!1/4!(m)$);
      \draw [name path=B_1',blue](m)--($(m)!1/3!(i)$);
      \draw [name path=B_2,blue](i')--($(i')!1/4!(m')$);
      \draw [name path=B_2',blue](m')--($(m')!1/3!(i')$);
      \draw [name path=C_1,black](j)--($(j)!1/4!(n)$);
      \draw [name path=C_1',black](o)--(n);
      \draw [name path=C_2,black](j')--($(j')!1/4!(n')$);
      \draw [name path=C_2',black]($(e)!1/2-1/60!(b)$)--(n');
      \node[above,black]at($(f)!1/3!(b)$){$C_1$};
      \node[right,blue]at($(f)!1/2!(d)$){$B_1$};
      \node[red]at($(b)!1/3!(f)$){$A_1$};
\node at ($(a)+(0.5,0)$) {$+$};
      \coordinate (x) at (-1,0);
      \coordinate (a) at ($(0,-1/2)+(x)$);
      \coordinate (b) at ($(1,-1/2)+(x)$);
      \coordinate (c) at ($(0,1/2)+(x)$);
      \coordinate (d) at ($(1,1/2)+(x)$);
      \coordinate (g) at ($(1/2,-1)+(x)$);
      \coordinate (o) at ($(a)!1/2!(c)$);
      \coordinate (f) at ($(b)!1/2!(d)$);
      \draw[red,name path=A_1](a)--(g);
      \draw[green](g)--(b);
      \draw[blue,name path=B_1](c)--(d);            
      \draw[blue,name path=B_2]($(c)+(0,-1/30)$)--($(d)+(0,-1/30)$);
      \draw[green,ultra thick](a)--(c);
      \draw[red,ultra thick](b)--(d);
      \draw[black,name path=C_1](o)--(f);
      \draw[black,name path=C_2]($(o)+(0,-1/30)$)--($(f)+(0,-1/30)$);
      \draw[red,name path=A_2]($(a)+(0,1/30)$)--($(d)!1/2!(c)$);
      \node[red,below]at ($(a)!1/2!(g)$){$A_1$};
      \node[red,right]at ($(b)!1/2!(d)$){$\triangle$};
   \node at ($(f)+(0.6,0)$){$+$};
      \coordinate (x) at (-0.4,0);
      \coordinate (a) at ($(2.5,-1/2)+(x)$);
      \coordinate (c) at ($(2.5,1/2)+(x)$);
      \coordinate (d) at ($(1.5,0)+(x)$);
      \coordinate (e) at ($(3.5,0)+(x)$);
      \coordinate (f) at ($(2,-1)+(x)$);
      \draw[red,ultra thick](c)--(a);
      \draw[black]($(a)!1/2!(d)$)--($(e)!1/2!(c)$);
      \draw[black]($(a)!1/2!(e)$)--($(d)!1/2!(c)$);
      \draw[blue](d)--(c);
      \draw[blue](e)--(c);
      \draw[green](a)--(d);
      \draw[green](a)--(e);
      \draw[green](a)--(f);
\end{tikzpicture}\\

\par\end{center}

We have $\left(K_{\mathcal{Y}}+\mathcal{D}\right)|_{\Sigma_{0}}.C\geq0$
for all the curves $C$ in $\Sigma_{0}$, and in particular $\left(K_{\mathcal{Y}}+\mathcal{D}\right)|_{\mathbb{P}^{1}\times\mathbb{P}^{1}}.C=0$
for all the curves $C$ in $\mathbb{P}^{1}\times\mathbb{P}^{1}$.
Running the minimal model program, we obtain the canonical model by
contracting $A_{1},B_{1},C_{1}$ in $Bl_{4}\mathbb{P}^{2}$, $A_{1},\triangle$
in $Bl_{1}\mathbb{F}_{1}$ and the whole component $\mathbb{P}^{1}\times\mathbb{P}^{1}$,
where $\triangle$ is the double locus. The central fiber of the resulting
canonical model is $\mathbb{P}^{1}\times\mathbb{P}^{1}\cup S$, where
$S$ is obtained from $\mathbb{F}_{2}$ by contracting the (-2)-curve.
So $S$ is a surface with an $A_{1}$-singularity. We call it Case
5.

\begin{center}
\begin{tikzpicture}[xscale=1.5,yscale=1.5][font=\tiny]
  \coordinate (x) at (0,-0.5);
  \coordinate (a) at ($(0,0)+(x)$);
  \coordinate (b) at ($(1,0)+(x)$);
  \coordinate (c) at ($(1,1)+(x)$);
  \coordinate (d) at ($(0,1)+(x)$);
  \draw [name path=B_3,blue](a)--(b);
  \draw [name path=B_0,blue](b)--(c);
  \draw [name path=C_3,black](c)--(d);
  \draw [name path=C_0,black](d)--(a);
  \draw [name path=A_0,red]($(a)!1/2!(b)$)--($(c)!1/2!(d)$);
  \draw [name path=A_3,red]($(b)!1/2!(c)$)--($(d)!1/2!(a)$);
  \draw [name path=triangle,green,ultra thick](b)--(d);
\node at ($(1.5,0.5)+(x)$) {$+$};
\node at ($(1.5,0.5)+(0,-1)+(x)$){$Case$ $5$};

      \coordinate (x) at (2,0);
      \coordinate (a) at ($(0,-1/2)+(x)$);
      \coordinate (b) at ($(1,0)+(x)$);
      \coordinate (c) at ($(0,1/2)+(x)$);
      \coordinate (o) at ($(a)!1/2!(c)$);
      \draw[green](a)--(b);
      \draw[blue,name path=B_1](c)--(b);            
      \draw[blue,name path=B_2]($(c)+(0,-1/30)$)--($(b)+(-2/30,0)$);
      \draw[green,ultra thick](a)--(c);
      \draw[black,name path=C_1](o)--(b);
      \draw[black,name path=C_2]($(o)+(0,-1/30)$)--($(b)+(-2/30,-1/30)$);
      \draw[red,name path=A_2](a)--($(b)!1/2!(c)$);
   \path[fill=red](b)circle[radius=0.03];
    \node [right,red] at (b){$A_1-singularity$};
\end{tikzpicture}\\

\par\end{center}

The following correspond to the non-nodal case with $K^{2}=4$.

\textbf{Case 6.} Similar to case 3, but the point $P_{1}$ is on $B_{0}$
instead. 

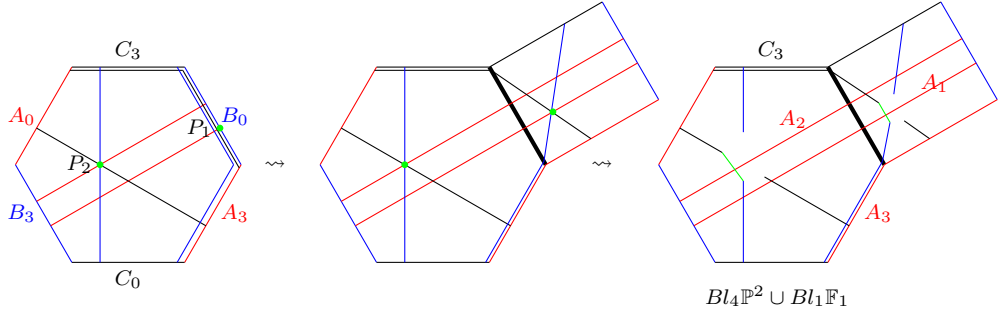
\begin{center}
\begin{tikzpicture}[xscale=1.5,yscale=1.5][font=\tiny]
      \coordinate (a) at (1,0);
      \coordinate (b) at (1/2,3^0.5/2);
      \coordinate (c) at (-1/2,3^0.5/2);
      \coordinate (d) at (-1,0);
      \coordinate (e) at (-1/2,-3^0.5/2);
      \coordinate (f) at (1/2,-3^0.5/2);
      \draw[blue](a)--(b);
      \draw[black](b)--(c);
      \draw[red](c)--(d);
      \draw[blue](d)--(e);
      \draw[black](e)--(f);
      \draw[red](f)--(a);
      \coordinate (g) at ($(a)!5/8!(b)$);
      \coordinate (h) at ($(a)!3/8!(b)$);
      \coordinate (i) at ($(b)!3/4!(c)$);
      \coordinate (j) at ($(c)!1/4!(d)$);
      \coordinate (k) at ($(c)!5/8!(d)$);
      \coordinate (l) at ($(d)!3/8!(e)$);
      \coordinate (m) at ($(e)!1/4!(f)$);
      \coordinate (n) at ($(e)!3/8!(d)$);
      \coordinate (o) at ($(f)!3/8!(a)$);
      \coordinate (p) at ($(f)!3/4!(a)$);
      \draw [name path=A_1,red](g)--(l);
      \draw [name path=A_2,red](h)--(n);
      \draw [name path=C_1,black](k)--(o);
      \draw [name path=C_2,black]($(c)!1/30!(d)$)--($(b)!1/60!(e)$);
      \draw [name path=C_2',black]($(b)!1/60!(e)$)--($(a)!1/30!(f)$);
      \draw [name path=B_1,blue](i)--(m);
      \draw [name path=B_2,blue]($(b)!2/30!(c)$)--($(a)!2/60!(d)$);
      \draw [name path=B_2',blue]($(a)!2/60!(d)$)--($(f)!2/30!(e)$);
      \node[above,black]at($(b)!1/2!(c)$){$C_3$};
      \node[below,black]at($(e)!1/2!(f)$){$C_0$};
      \node[right,blue]at($(a)!1/2!(b)$){$B_0$};
      \node[left,blue]at($(d)!1/2!(e)$){$B_3$};
      \node[left,red]at($(c)!1/2!(d)$){$A_0$};
      \node[right,red]at($(a)!1/2!(f)$){$A_3$};
      \path[fill=green]($(a)!5/8!(d)$)circle[radius=0.03];
      \node [left] at ($(a)!5/8!(d)$){$P_2$};
      \path[fill=green](h)circle[radius=0.03];
      \node[left]at(h){$P_1$};

\node at ($(a)+(0.3,0)$){$\rightsquigarrow$};
      \coordinate (a) at (1+2.7,0);
      \coordinate (b) at (1/2+2.7,3^0.5/2);
      \coordinate (c) at (-1/2+2.7,3^0.5/2);
      \coordinate (d) at (-1+2.7,0);
      \coordinate (e) at (-1/2+2.7,-3^0.5/2);
      \coordinate (f) at (1/2+2.7,-3^0.5/2);
      \draw[black,ultra thick](a)--(b);
      \draw[black](b)--(c);
      \draw[red](c)--(d);
      \draw[blue](d)--(e);
      \draw[black](e)--(f);
      \draw[red](f)--(a);
      \coordinate (g) at ($(a)!5/8!(b)$);
      \coordinate (h) at ($(a)!3/8!(b)$);
      \coordinate (i) at ($(b)!3/4!(c)$);
      \coordinate (j) at ($(c)!1/4!(d)$);
      \coordinate (k) at ($(c)!5/8!(d)$);
      \coordinate (l) at ($(d)!3/8!(e)$);
      \coordinate (m) at ($(e)!1/4!(f)$);
      \coordinate (n) at ($(e)!3/8!(d)$);
      \coordinate (o) at ($(f)!3/8!(a)$);
      \coordinate (p) at ($(f)!3/4!(a)$);
      \draw [name path=A_1,red](g)--(l);
      \draw [name path=A_2,red](h)--(n);
      \draw [name path=C_1,black](k)--(o);
      \draw [name path=C_2,black]($(c)!1/30!(d)$)--($(b)!1/30!(a)$);
      \draw [name path=B_1,blue](i)--(m);
      \draw [name path=B_2,blue]($(a)!1/30!(b)$)--($(f)!1/30!(e)$);
      \path[fill=green]($(a)!5/8!(d)$)circle[radius=0.03];
      \coordinate (r) at ($(d)!1+2/3!(b)$);
      \coordinate (s) at ($(e)!1+2/3!(a)$);
      \coordinate (t) at ($(s)!3/8!(r)$);
      \draw [name path=C_3,black](b)--(r);
      \draw [name path=B_0,blue](r)--(s);
      \draw [name path=A_3,red](a)--(s);
      \draw [name path=A_1,red](g)--($(l)!1+2/3!(g)$);
      \draw [name path=A_2,red](h)--($(s)!3/8!(r)$);
      \draw [name path=C_2,black]($(b)!1/30!(a)$)--($(a)!2/5!(s)$);
      \draw [name path=B_2,blue]($(a)!1/30!(b)$)--($(b)!2/3!(r)$);
      \path[fill=green]($(h)!1/4!(t)$)circle[radius=0.03];
\node at ($(a)+(0.5,0)$){$\rightsquigarrow$};
      \coordinate (a) at (1+5.7,0);
      \coordinate (b) at (1/2+5.7,3^0.5/2);
      \coordinate (c) at (-1/2+5.7,3^0.5/2);
      \coordinate (d) at (-1+5.7,0);
      \coordinate (e) at (-1/2+5.7,-3^0.5/2);
      \coordinate (f) at (1/2+5.7,-3^0.5/2);
      \draw[black,ultra thick](a)--(b);
      \draw[black](b)--(c);
      \draw[red](c)--(d);
      \draw[blue](d)--(e);
      \draw[black](e)--(f);
      \draw[red](f)--(a);
      \coordinate (g) at ($(a)!5/8!(b)$);
      \coordinate (h) at ($(a)!3/8!(b)$);
      \coordinate (i) at ($(b)!3/4!(c)$);
      \coordinate (j) at ($(c)!1/4!(d)$);
      \coordinate (k) at ($(c)!5/8!(d)$);
      \coordinate (l) at ($(d)!3/8!(e)$);
      \coordinate (m) at ($(e)!1/4!(f)$);
      \coordinate (n) at ($(e)!3/8!(d)$);
      \coordinate (o) at ($(f)!3/8!(a)$);
      \coordinate (p) at ($(f)!3/4!(a)$);
      \draw [name path=A_1,red](g)--(l);
      \draw [name path=A_2,red]($(a)!3/8!(b)$)--($(e)!3/8!(d)$);
      \draw [name path=C_1,black](k)--($(k)!1/4!(o)$);
      \draw [name path=C_1',black](o)--($(o)!1/2!(k)$);
      \draw [name path=C_2,black]($(c)!1/30!(d)$)--($(b)!1/30!(a)$);
      \draw [name path=B_1,blue](i)--($(i)!1/3!(m)$);
      \draw [name path=B_1',blue](m)--($(m)!5/12!(i)$);
      \draw [name path=B_2,blue]($(a)!1/30!(b)$)--($(f)!1/30!(e)$);
      \draw [name path=E_1,green]($(k)!1/4!(o)$)--($(m)!5/12!(i)$);
      \node[above,black]at($(b)!1/2!(c)$){$C_3$};
      \node[right,red]at($(a)!1/2!(f)$){$A_3$};
      \node[above,red]at($(g)!1/3!(l)$){$A_2$};
      \coordinate (r) at ($(d)!1+2/3!(b)$);
      \coordinate (s) at ($(e)!1+2/3!(a)$);
      \coordinate (t) at ($(r)!3/8!(s)$);
      \coordinate (t') at ($(s)!3/8!(r)$);
      \coordinate (u) at ($(b)!2/3!(r)$);
      \coordinate (v) at ($(a)!2/5!(s)$);
      \draw [name path=C_3,black](b)--(r);
      \draw [name path=B_0,blue](r)--(s);
      \draw [name path=A_3,red](a)--(s);
      \draw [name path=A_1,red](g)--(t);
      \draw [name path=A_2,red](h)--(t');
      \draw [name path=C_2,black]($(b)!1/30!(a)$)--($(b)!1/2!(v)$);
      \draw [name path=C_2',black](v)--($(v)!1/4!(b)$);
      \draw [name path=B_2,blue]($(a)!1/30!(b)$)--($(a)!3/10!(u)$);
      \draw [name path=B_2',blue](u)--($(u)!1/2!(a)$);
      \draw [name path=E_2,green]($(a)!3/10!(u)$)--($(b)!1/2!(v)$);
      \node [right,red]at($(h)!1/2!(t)$){$A_1$};
\node at($(m)+(0.3,-0.3)$){$Bl_{4}\mathbb{P}^{2}\cup Bl_{1}\mathbb{F}_{1}$};
\end{tikzpicture}
\par\end{center}

The central fiber of the resulting lc model is $\mathbb{P}^{1}\times\mathbb{P}^{1}\cup\mathbb{P}^{1}\times\mathbb{P}^{1}$.
It is not isomorphic to Case 3, which is also $\mathbb{P}^{1}\times\mathbb{P}^{1}\cup\mathbb{P}^{1}\times\mathbb{P}^{1}$,
since the line arrangements are not isomorphic.

\begin{center}
\begin{tikzpicture}[xscale=2,yscale=2][font=\tiny]
      \coordinate (e) at (0,0);\node [below] at (e){1};
      \coordinate (f) at (0,1);\node [above] at (f){4};
      \coordinate (g) at (-3^0.5/2,1/2);
      \coordinate (h) at (3^0.5/2,1/2);
      \draw [name path=triangle,ultra thick,black] (e)--(f);
      \draw [name path=B_1,blue] (g)--(f); 
      \draw [name path=A_0,red] (f)--(h);    
      \draw [name path=C_0,black] (g)--(e);
      \draw [name path=C_1,black] (e)--(h);
      \draw [name path=A_2,red] ($(g)!1/3!(f)$)--($(e)!1/3!(h)$);
      \draw [name path=B_3,blue] ($(e)!2/3!(g)$)--($(h)!2/3!(f)$);
      \draw [name path=E,green]($(f)!1/3!(g)$)--($(h)!1/3!(e)$);
      \node[black] at ($(e)+(0,-1/4)$){$\mathbb{P}^{1}\times\mathbb{P}^{1}$};
\node[black] at ($(h)+(0.4,0)$){$\cup$};
\node[black] at ($(h)+(0.4,-1)$){$Case$ $6$};

      \coordinate (e) at (0+2.5,0);\node [below] at (e){1};
      \coordinate (f) at (0+2.5,1);\node [above] at (f){4};
      \coordinate (g) at (-3^0.5/2+2.5,1/2);
      \coordinate (h) at (3^0.5/2+2.5,1/2);
      \draw [name path=triangle,ultra thick,black] (e)--(f);
      \draw [name path=C_2,black] (g)--(f); 
      \draw [name path=C_3,black] (f)--(h);    
      \draw [name path=A_3,red] (g)--(e);
      \draw [name path=B_2,blue] (e)--(h);
      \draw [name path=B_0,blue] ($(g)!1/3!(f)$)--($(e)!1/3!(h)$);
      \draw [name path=A_1,red] ($(e)!2/3!(g)$)--($(h)!2/3!(f)$);
      \draw [name path=E_1,green] ($(h)!1/3!(f)$)--($(e)!1/3!(g)$);
      \node[black] at ($(e)+(0,-1/4)$){$\mathbb{P}^{1}\times\mathbb{P}^{1}$};
\end{tikzpicture}\\

\par\end{center}

There is a further degeneration as follows, where the central fiber
of the resulting lc model is a union of four copies of $\mathbb{P}^{2}$. 

\begin{tikzpicture}[xscale=2,yscale=2][font=\tiny]
      \coordinate (e) at (0,0);\node [below] at (e){1};
      \coordinate (f) at (0,1);\node [above] at (f){4};
      \coordinate (g) at (-3^0.5/2,1/2);
      \coordinate (h) at (3^0.5/2,1/2);
      \draw [name path=triangle,ultra thick,black] (e)--(f);
      \draw [name path=B_1,blue] (g)--(f); 
      \draw [name path=A_0,red] (f)--(h);    
      \draw [name path=C_0,black] (g)--(e);
      \draw [name path=C_1,black] (e)--(h);
      \draw [name path=A_2,red] ($(g)!1/3!(f)$)--($(e)!1/3!(h)$);
      \draw [name path=B_3,blue] ($(g)!1/30!(e)$)--($(f)!1/30!(e)$);
      \draw [name path=E_2,green]($(f)!1/30!(e)$)--($(h)!1/30!(e)$);
      \node[black] at ($(e)+(0,-1/4)$){$\mathbb{P}^{1}\times\mathbb{P}^{1}$};
\node[black] at ($(h)+(0.4,0)$){$\cup$};

      \coordinate (x) at (2.5,0);
      \coordinate (e) at ($(0,0)+(x)$);\node [below] at (e){1};
      \coordinate (f) at ($(0,1)+(x)$);\node [above] at (f){4};
      \coordinate (g) at ($(-3^0.5/2,1/2)+(x)$);
      \coordinate (h) at ($(3^0.5/2,1/2)+(x)$);
      \draw [name path=triangle,ultra thick,black] (e)--(f);
      \draw [name path=C_2,black] (g)--(f); 
      \draw [name path=C_3,black] (f)--(h);    
      \draw [name path=A_3,red] (g)--(e);
      \draw [name path=B_2,blue] (e)--(h);
      \draw [name path=B_0,blue] ($(g)!1/3!(f)$)--($(e)!1/3!(h)$);
      \draw [name path=A_1,red] ($(g)!1/30!(e)$)--($(f)!1/30!(e)$);
      \draw [name path=E_1,green] ($(h)!1/3!(f)$)--($(e)!1/3!(g)$);
      \node[black] at ($(e)+(0,-1/4)$){$\mathbb{P}^{1}\times\mathbb{P}^{1}$};
     \node at ($(h)+(0.4,0)$){$\rightsquigarrow$};

      \coordinate (x) at (4.7,1/2);
      \coordinate (a) at ($(1/2,0)+(x)$);
      \coordinate (b) at ($(0,1/2)+(x)$);
      \coordinate (d) at ($(-1/2,0)+(x)$);
      \coordinate (e) at ($(0,-1/2)+(x)$);
      \coordinate (O) at ($(0,0)+(x)$);
      \draw [black,name path=C_3] (a)--(b);
      \draw [blue,name path=B_1] (b)--(d);
      \draw [black,name path=C_3] (d)--(e);
      \draw [red,name path=A_3] (e)--(a);
      \draw [black, ultra thick] (e)--(b);
      \draw [blue, ultra thick] (a)--(d);
      \draw [red,name path=A_2] ($(O)!1/3!(b)$)--(a);
      \draw [black,name path=C_2] ($(O)!1/3!(a)$)--(b);
      \draw [blue,name path=B_2] ($(O)!1/3!(a)$)--(e);
      \draw [blue,name path=B_0] ($(O)!1/3!(e)$)--(a);
      \draw [black,name path=C_2] ($(O)!1/3!(d)$)--(e);
      \draw [red,name path=A_2] ($(O)!1/3!(e)$)--(d);
      \draw [red,name path=A_0] ($(O)!1/3!(d)$)--(b);
      \draw [black,name path=B_3] ($(O)!1/3!(b)$)--(d);
      \node at ($(e)+(0,-0.5)$){$Case$ $7$};
\end{tikzpicture}\\

\textbf{Case 8.} For non-nodal case, when $A_{1},A_{2},B_{1},B_{2},C_{1},C_{2}$
intersect at one point. The central fiber of the lc model is $\mathbb{P}^{1}\times\mathbb{P}^{1}\cup Bl_{2}\mathbb{P}^{2}$.

\begin{center}
\begin{tikzpicture}[xscale=1.5,yscale=1.5][font=\tiny]
      \coordinate (a) at (1,0);
      \coordinate (b) at (1/2,3^0.5/2);
      \coordinate (c) at (-1/2,3^0.5/2);
      \coordinate (d) at (-1,0);
      \coordinate (e) at (-1/2,-3^0.5/2);
      \coordinate (f) at (1/2,-3^0.5/2);
      \draw[blue](a)--(b);
      \draw[black](b)--(c);
      \draw[red](c)--(d);
      \draw[blue](d)--(e);
      \draw[black](e)--(f);
      \draw[red](f)--(a);
      \coordinate (g) at ($(a)!1/2-1/60!(b)$);
      \coordinate (h) at ($(a)!1/2+1/60!(b)$);
      \coordinate (i) at ($(b)!1/2!(c)$);
      \coordinate (i') at ($(c)!1/2-1/30!(b)$);
      \coordinate (j) at ($(c)!1/2!(d)$);
      \coordinate (j') at ($(c)!1/2+1/30!(d)$);
      \coordinate (k) at ($(d)!1/2-1/60!(e)$);
      \coordinate (l) at ($(d)!1/2+1/60!(e)$);
      \coordinate (m) at ($(e)!1/2!(f)$);
      \coordinate (m') at ($(e)!1/2-1/30!(f)$);
      \coordinate (n) at ($(f)!1/2!(a)$);
      \coordinate (n') at ($(f)!1/2-1/30!(a)$);
      \draw [name path=A_1,red](h)--(k);
      \draw [name path=A_2,red](g)--(l);
      \draw [name path=B_1,blue](i)--(m);
      \draw [name path=B_2,blue](i')--(m');
      \draw [name path=C_1,black](j)--(n);
      \draw [name path=C_2,black](j')--(n');
      \node[above,black]at($(b)!1/2!(c)$){$C_3$};
      \node[below,black]at($(e)!1/2!(f)$){$C_0$};
      \node[right,blue]at($(a)!1/2!(b)$){$B_0$};
      \node[left,blue]at($(d)!1/2!(e)$){$B_3$};
      \node[left,red]at($(c)!1/2!(d)$){$A_0$};
      \node[right,red]at($(a)!1/2!(f)$){$A_3$};
      \path[fill=green]($(a)!1/2!(d)$)circle[radius=0.03];
\node at ($(a)+(0.5,0)$) {$\rightsquigarrow$};

      \coordinate (a) at (1+3,0);
      \coordinate (b) at (1/2+3,3^0.5/2);
      \coordinate (c) at (-1/2+3,3^0.5/2);
      \coordinate (d) at (-1+3,0);
      \coordinate (e) at (-1/2+3,-3^0.5/2);
      \coordinate (f) at (1/2+3,-3^0.5/2);
      \draw[blue](a)--(b);
      \draw[black](b)--(c);
      \draw[red](c)--(d);
      \draw[blue](d)--(e);
      \draw[black](e)--(f);
      \draw[red](f)--(a);
      \coordinate (g) at ($(a)!3/4+1/30!(b)$);
      \coordinate (h) at ($(a)!3/4!(b)$);
      \coordinate (i) at ($(b)!3/4!(c)$);
      \coordinate (i') at ($(c)!1/4-1/30!(b)$);
      \coordinate (j) at ($(c)!1/2!(d)$);
      \coordinate (j') at ($(c)!1/2+1/30!(d)$);
      \coordinate (k) at ($(d)!1/4!(e)$);
      \coordinate (l) at ($(d)!1/4-1/30!(e)$);
      \coordinate (m) at ($(e)!1/4!(f)$);
      \coordinate (m') at ($(e)!1/4-1/30!(f)$);
      \coordinate (n) at ($(f)!1/2!(a)$);
      \coordinate (n') at ($(f)!1/2-1/30!(a)$);
      \coordinate (o) at ($(c)!1/2!(f)$);
      \draw [name path=A_1,red](h)--($(h)!1/3!(k)$);
      \draw [name path=A_1',red](k)--($(o)!1/2!(d)$);
      \draw [name path=A_2,red](g)--($(g)!1/3!(l)$);
      \draw [name path=A_2',red](l)--($(l)!1/4!(g)$);
      \draw [name path=E,green,ultra thick]($(o)!1/2!(c)$)--($(o)!1/2!(d)$);
      \draw [name path=B_1,blue](i)--($(i)!1/4!(m)$);
      \draw [name path=B_1',blue](m)--($(m)!1/3!(i)$);
      \draw [name path=B_2,blue](i')--($(i')!1/4!(m')$);
      \draw [name path=B_2',blue](m')--($(m')!1/3!(i')$);
      \draw [name path=C_1,black](j)--($(j)!1/4!(n)$);
      \draw [name path=C_1',black](o)--(n);
      \draw [name path=C_2,black](j')--($(j')!1/4!(n')$);
      \draw [name path=C_2',black]($(e)!1/2-1/60!(b)$)--(n');
\node at ($(a)+(0.5,0)$) {$+$};
      \coordinate (x) at (5,0);
      \coordinate (a) at ($(0,-1/2)+(x)$);
      \coordinate (b) at ($(1,-1/2)+(x)$);
      \coordinate (c) at ($(1,1/2)+(x)$);
      \coordinate (d) at ($(0,1/2)+(x)$);
      \draw[red,name path=A_1](a)--(b);
      \draw[red,name path=A_2](a)--(c);
      \draw[green,ultra thick](a)--(d);
      \draw[blue,name path=B_1](d)--(b);
      \draw[blue,name path=B_2](d)--(c);
      \draw[black,name path=C_1](b)--($(a)!1/2!(d)$);
      \draw[black,name path=C_2](c)--($(a)!1/2!(d)$);
      \node[right]at(b){$P_1$};
      \node[right]at(c){$P_2$};
      \path[fill=green](b)circle[radius=0.03];
      \path[fill=green](c)circle[radius=0.03];
    \coordinate (x) at (0,-2.5);
    \coordinate (a) at ($(0,0)+(x)$);
    \coordinate (b) at ($(1,0)+(x)$);
    \coordinate (c) at ($(1,1)+(x)$);
    \coordinate (d) at ($(0,1)+(x)$);
    \draw [name path=B_3,blue](a)--(b);
    \draw [name path=B_0,blue](b)--(c);
    \draw [name path=C_3,black](c)--(d);
    \draw [name path=C_0,black](d)--(a);
    \draw [name path=A_0,red]($(a)!1/2!(b)$)--($(c)!1/2!(d)$);
    \draw [name path=A_3,red]($(b)!1/2!(c)$)--($(d)!1/2!(a)$);
    \draw [name path=triangle,green,ultra thick](b)--(d);
  \node at ($(1.5,0.5)+(x)$) {$+$};
  \node at ($(1.5,0.5)+(0,-1)+(x)$){$Case$ $8$};
      \coordinate (x) at (3.3,-1.7);
      \coordinate (a) at ($(0,-1/2)+(x)$);
      \coordinate (c) at ($(0,1/2)+(x)$);
      \coordinate (d) at ($(-3^0.5/2,0)+(x)$);
      \coordinate (e) at ($(3^0.5/2,0)+(x)$);
      \coordinate (f) at ($(c)!1/2!(e)$);
      \coordinate (g) at ($(c)!1/2!(d)$);
      \coordinate (o) at ($(d)!1/2!(e)$);

      \draw[green,ultra thick](c)--(a);
      \draw[blue,name path=B_1](c)--($(c)!1+1/2!(d)$);
      \draw[blue,name path=B_2](c)--($(c)!1+1/2!(e)$);
      \draw[black,name path=C_1](f)--($(f)!3!(o)$);
      \draw[black,name path=C_2](g)--($(g)!3!(o)$);
      \draw[red,name path=A_1](e)--($(e)!1+1/2!(a)$);
      \draw[red,name path=A_2](d)--($(d)!1+1/2!(a)$);      
      \draw[green,name path=E_1]($(c)!1+1/2!(d)$)--($(e)!1+1/2!(a)$);
      \draw[green,name path=E_2]($(c)!1+1/2!(e)$)--($(d)!1+1/2!(a)$); 
\end{tikzpicture}\\

\par\end{center}

In total, there are 5 types of degenerations with reducible lc models
for $K^{2}=4$ nodal case and 3 types of degenerations for $K^{2}=4$
non-nodal case up to the symmetry group $\mathbb{Z}_{2}$. All of
them could be obtained from the cases with $K^{2}=5$.\\

\begin{center}
\begin{tabular}{|c|c|c|}
\hline 
$K^{2}=5$ & $K^{2}=4$ & further degenerations\tabularnewline
\hline 
\hline 
Case 1 & Case 1 & Case 2\tabularnewline
\hline 
Case 2 & Case 2 & \tabularnewline
\hline 
Case 3 & Case 6 & \tabularnewline
\hline 
Case 4 & Cases 3,6 & Cases 2,7\tabularnewline
\hline 
Case 5  & Case 2,7 & \tabularnewline
\hline 
Case 6 & Case 4,8 & Case 5\tabularnewline
\hline 
\end{tabular}
\par\end{center}

\section{Burniat surfaces with $K^{2}=3$\label{degree3}}

Consider the surface $\mathbb{P}^{2}$ with 9 lines and 3 points $P_{1},P_{2}$
and $P_{3}$ which are the intersections of 3 lines. A Burniat surface
$X$ in $M_{Bur}^{3}$ is the canonical model of a $\mathbb{Z}_{2}^{2}$-cover
of $\Sigma=Bl_{6}\mathbb{P}^{2}$ for the building data $D_{a},D_{b}$
and $D_{c}.$ Here $\Sigma$ is the blown up of $\mathbb{P}^{2}$
at the six points $P_{A},P_{B},P_{c}$ and $P_{1},P_{2},P_{3}$. 

There are three $(-2)$-curves $A_{1},B_{1},C_{1}$ in $\Sigma$ and
$K_{\Sigma}+D=-\dfrac{1}{2}K_{\Sigma}$ is nef but not ample. The
canonical model $\Sigma^{c}$ of $\Sigma$ is obtained from $\Sigma$
by contracting the three $(-2)$-curves. Stable Burniat surfaces $X$
with $K_{X}^{2}=3$ are $\mathbb{Z}_{2}^{2}$-covers of the canonical
models $\Sigma^{c}$ of $\Sigma$ with the building data $\frac{1}{2}D$.
The general fiber of a one-parameter family is $\Sigma^{c}$ and it
contains three $A_{1}$-singularities. We denote the singularities,
by contracting from $A_{1},B_{1},C_{1}$, by $Q_{1},Q_{2},Q_{3}$. 

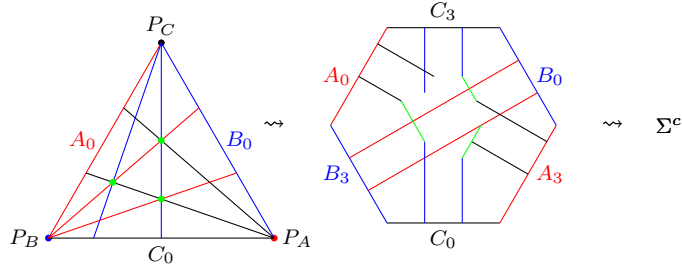
\begin{center}
\begin{tikzpicture}[xscale=1.5,yscale=1.5][font=\tiny]
      \coordinate (a) at (1,0);
      \node[right] at (a){$P_A$};
      \path[fill=red](a)circle[radius=0.03];
      \coordinate (b) at (-1,0);
      \node[left] at (b){$P_B$};
      \path[fill=blue](b)circle[radius=0.03];
      \coordinate (c) at (0,3^0.5);
      \node[above] at (c){$P_C$};
      \path[fill=black](c)circle[radius=0.03];
      \draw[name path=C_0,black](a)--(b);
      \node[below,black]at($(a)!1/2!(b)$){$C_0$};
      \draw[name path=A_0,red](b)--(c);
      \node[left,red]at($(b)!1/2!(c)$){$A_0$};
      \draw[name path=B_0,blue](c)--(a);
      \node[right,blue]at($(c)!1/2!(a)$){$B_0$};
      \draw [name path=A_1,red]($(c)!1/3!(a)$)--(b);
      \draw [name path=A_2,red]($(a)!1/3!(c)$)--(b);
      \draw [name path=B_1,blue]($(a)!1/2!(b)$)--(c);
      \draw [name path=B_2,blue]($(a)!8/10!(b)$)--(c);
      \draw [name path=C_1,black]($(b)!1/3!(c)$)--(a);
      \draw [name path=C_2,black]($(c)!1/3!(b)$)--(a);
      \path [name intersections={of=A_1 and B_1}]; 
      \coordinate (P_1) at (intersection-1);
      \path[fill=green](P_1)circle[radius=0.03];
      \path [name intersections={of=A_1 and B_2}]; 
      \coordinate (P_2) at (intersection-1);
      \path[fill=green](P_2)circle[radius=0.03];
      \path [name intersections={of=A_2 and B_1}]; 
      \coordinate (P_3) at (intersection-1);
      \path[fill=green](P_3)circle[radius=0.03];   

      \coordinate (a) at (1+2.5,0+1);
      \coordinate (b) at (1/2+2.5,3^0.5/2+1);
      \coordinate (c) at (-1/2+2.5,3^0.5/2+1);
      \coordinate (d) at (-1+2.5,0+1);
      \coordinate (e) at (-1/2+2.5,-3^0.5/2+1);
      \coordinate (f) at (1/2+2.5,-3^0.5/2+1);
      \draw[blue](a)--(b);
      \draw[black](b)--(c);
      \draw[red](c)--(d);
      \draw[blue](d)--(e);
      \draw[black](e)--(f);
      \draw[red](f)--(a);
      \coordinate (g) at ($(a)!1/3!(b)$);
      \coordinate (h) at ($(a)!2/3!(b)$);
      \coordinate (i) at ($(b)!1/3!(c)$);
      \coordinate (j) at ($(b)!2/3!(c)$);
      \coordinate (k) at ($(c)!1/6!(d)$);
      \coordinate (l) at ($(c)!1/2!(d)$);
      \coordinate (m) at ($(d)!1/3!(e)$);
      \coordinate (n) at ($(d)!2/3!(e)$);
      \coordinate (q) at ($(e)!1/3!(f)$);
      \coordinate (r) at ($(e)!2/3!(f)$);
      \coordinate (s) at ($(f)!1/2!(a)$);
      \coordinate (t) at ($(a)!1/6!(f)$);
      \draw [name path=A_1,red](h)--(m);
      \draw [name path=A_2,red](g)--(n);
      \draw [name path=B_1,blue](i)--($(i)!1/4!(r)$);
      \draw [name path=B_1',blue](r)--($(r)!1/3!(i)$);
      \draw [name path=B_2,blue](j)--($(j)!1/3!(q)$);
      \draw [name path=B_2',blue](q)--($(q)!5/12!(j)$);
      \draw [name path=C_1,black](k)--($(k)!1/3!(t)$);   
      \draw [name path=C_1',black](t)--($(t)!5/12!(k)$);
      \draw [name path=C_2,black](l)--($(l)!1/4!(s)$);
      \draw [name path=C_2',black](s)--($(s)!1/3!(l)$);
      \draw [name path=E_1,green]($(i)!1/4!(r)$)--($(t)!5/12!(k)$);
      \draw [name path=E_2,green]($(l)!1/4!(s)$)--($(q)!5/12!(j)$);
      \draw [name path=E_3,green]($(g)!1/3!(n)$)--($(r)!1/3!(i)$);
      \node[above,black]at($(b)!1/2!(c)$){$C_3$};
      \node[below,black]at($(e)!1/2!(f)$){$C_0$};
      \node[right,blue]at($(a)!1/2!(b)$){$B_0$};
      \node[left,blue]at($(d)!1/2!(e)$){$B_3$};
      \node[left,red]at($(c)!1/2!(d)$){$A_0$};
      \node[right,red]at($(a)!1/2!(f)$){$A_3$};
\node at ($(d)+(-0.5,0)$){$\rightsquigarrow$};
\node at ($(a)+(0.5,0)$){$\rightsquigarrow$};
\node at ($(a)+(1,0)$){$\Sigma^{c}$};\node at ($(d)+(-0.5,0)$){$\rightsquigarrow$};
\node at ($(a)+(0.5,0)$){$\rightsquigarrow$};
\node at ($(a)+(1,0)$){$\Sigma^{c}$};
\end{tikzpicture}     \\

\par\end{center}

\textbf{Case 1}\emph{.} The three points $P_{1},P_{2},P_{3}$ coincide.
Take a one-parameter family of $\Sigma$ with the general fiber $Bl_{5}\mathbb{P}^{2}$
with the central fiber the degenerating arrangement $\Sigma_{0}$.
Denote the one parameter family space by $\mathcal{Y}$. 

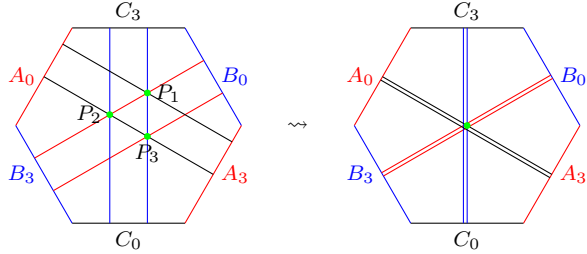
\begin{center}
\begin{tikzpicture}[xscale=1.5,yscale=1.5][font=\tiny]
      \coordinate (a) at (1,0+1);
      \coordinate (b) at (1/2,3^0.5/2+1);
      \coordinate (c) at (-1/2,3^0.5/2+1);
      \coordinate (d) at (-1,0+1);
      \coordinate (e) at (-1/2,-3^0.5/2+1);
      \coordinate (f) at (1/2,-3^0.5/2+1);
      \draw[blue](a)--(b);
      \draw[black](b)--(c);
      \draw[red](c)--(d);
      \draw[blue](d)--(e);
      \draw[black](e)--(f);
      \draw[red](f)--(a);
      \coordinate (g) at ($(a)!1/3!(b)$);
      \coordinate (h) at ($(a)!2/3!(b)$);
      \coordinate (i) at ($(b)!1/3!(c)$);
      \coordinate (j) at ($(b)!2/3!(c)$);
      \coordinate (k) at ($(c)!1/6!(d)$);
      \coordinate (l) at ($(c)!1/2!(d)$);
      \coordinate (m) at ($(d)!1/3!(e)$);
      \coordinate (n) at ($(d)!2/3!(e)$);
      \coordinate (q) at ($(e)!1/3!(f)$);
      \coordinate (r) at ($(e)!2/3!(f)$);
      \coordinate (s) at ($(f)!1/2!(a)$);
      \coordinate (t) at ($(a)!1/6!(f)$);
      \draw [name path=A_1,red](h)--(m);
      \draw [name path=A_2,red](g)--(n);
      \draw [name path=B_1,blue](i)--(r);
      \draw [name path=B_2,blue](j)--(q);
      \draw [name path=C_1,black](k)--(t);
      \draw [name path=C_2,black](l)--(s);
      \node[above,black]at($(b)!1/2!(c)$){$C_3$};
      \node[below,black]at($(e)!1/2!(f)$){$C_0$};
      \node[right,blue]at($(a)!1/2!(b)$){$B_0$};
      \node[left,blue]at($(d)!1/2!(e)$){$B_3$};
      \node[left,red]at($(c)!1/2!(d)$){$A_0$};
      \node[right,red]at($(a)!1/2!(f)$){$A_3$};
      \path [name intersections={of=A_1 and B_1}]; 
      \coordinate (P_1) at (intersection-1);
      \path[fill=green](P_1)circle[radius=0.03];
      \node [right]at (P_1){$P_1$};
      \path [name intersections={of=A_1 and B_2}]; 
      \path[fill=green](intersection-1)circle[radius=0.03];
      \node [left]at (intersection-1){$P_2$};
      \path [name intersections={of=A_2 and B_1}]; 
      \path[fill=green](intersection-1)circle[radius=0.03];
      \node  [below]at (intersection-1){$P_3$};
\node at ($(a)+(0.5,0)$){$\rightsquigarrow$};
      \coordinate (a) at (1+3,0+1);
      \coordinate (b) at (1/2+3,3^0.5/2+1);
      \coordinate (c) at (-1/2+3,3^0.5/2+1);
      \coordinate (d) at (-1+3,0+1);
      \coordinate (e) at (-1/2+3,-3^0.5/2+1);
      \coordinate (f) at (1/2+3,-3^0.5/2+1);
      \draw[blue](a)--(b);
      \draw[black](b)--(c);
      \draw[red](c)--(d);
      \draw[blue](d)--(e);
      \draw[black](e)--(f);
      \draw[red](f)--(a);
      \coordinate (g) at ($(a)!1/2-1/60!(b)$);
      \coordinate (h) at ($(a)!1/2+1/60!(b)$);
      \coordinate (i) at ($(b)!1/2!(c)$);
      \coordinate (i') at ($(c)!1/2-1/30!(b)$);
      \coordinate (j) at ($(c)!1/2!(d)$);
      \coordinate (j') at ($(c)!1/2+1/30!(d)$);
      \coordinate (k) at ($(d)!1/2-1/60!(e)$);
      \coordinate (l) at ($(d)!1/2+1/60!(e)$);
      \coordinate (m) at ($(e)!1/2!(f)$);
      \coordinate (m') at ($(e)!1/2-1/30!(f)$);
      \coordinate (n) at ($(f)!1/2!(a)$);
      \coordinate (n') at ($(f)!1/2-1/30!(a)$);
      \draw [name path=A_1,red](h)--(k);
      \draw [name path=A_2,red](g)--(l);
      \draw [name path=B_1,blue](i)--(m);
      \draw [name path=B_2,blue](i')--(m');
      \draw [name path=C_1,black](j)--(n);
      \draw [name path=C_2,black](j')--(n');
      \node[above,black]at($(b)!1/2!(c)$){$C_3$};
      \node[below,black]at($(e)!1/2!(f)$){$C_0$};
      \node[right,blue]at($(a)!1/2!(b)$){$B_0$};
      \node[left,blue]at($(d)!1/2!(e)$){$B_3$};
      \node[left,red]at($(c)!1/2!(d)$){$A_0$};
      \node[right,red]at($(a)!1/2!(f)$){$A_3$};
      \path[fill=green]($(a)!1/2!(d)$)circle[radius=0.03];
\end{tikzpicture}     
\par\end{center}

We first blow up $\mathcal{Y}$ at the point $P$ on the central fiber.
The central fiber becomes $Bl_{4}\mathbb{P}^{2}\cup\mathbb{P}^{2}$. 

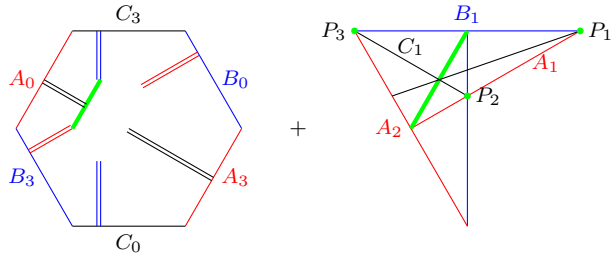
\begin{center}
\begin{tikzpicture}[xscale=1.5,yscale=1.5][font=\tiny]
      \coordinate (a) at (1,0);
      \coordinate (b) at (1/2,3^0.5/2);
      \coordinate (c) at (-1/2,3^0.5/2);
      \coordinate (d) at (-1,0);
      \coordinate (e) at (-1/2,-3^0.5/2);
      \coordinate (f) at (1/2,-3^0.5/2);
      \draw[name path=B_0,blue](a)--(b);
      \draw[name path=C_3,black](b)--(c);
      \draw[name path=A_0,red](c)--(d);
      \draw[name path=B_3,blue](d)--(e);
      \draw[name path=C_0,black](e)--(f);
      \draw[name path=A_3,red](f)--(a);
      \coordinate (g) at ($(a)!3/4+1/30!(b)$);
      \coordinate (h) at ($(a)!3/4!(b)$);
      \coordinate (i) at ($(b)!3/4!(c)$);
      \coordinate (i') at ($(c)!1/4-1/30!(b)$);
      \coordinate (j) at ($(c)!1/2!(d)$);
      \coordinate (j') at ($(c)!1/2+1/30!(d)$);
      \coordinate (k) at ($(d)!1/4!(e)$);
      \coordinate (l) at ($(d)!1/4-1/30!(e)$);
      \coordinate (m) at ($(e)!1/4!(f)$);
      \coordinate (m') at ($(e)!1/4-1/30!(f)$);
      \coordinate (n) at ($(f)!1/2!(a)$);
      \coordinate (n') at ($(f)!1/2-1/30!(a)$);
      \coordinate (o) at ($(c)!1/2!(f)$);
      \draw [name path=A_1,red](h)--($(h)!1/3!(k)$);
      \draw [name path=A_1',red](k)--($(o)!1/2!(d)$);
      \draw [name path=A_2,red](g)--($(g)!1/3!(l)$);
      \draw [name path=A_2',red](l)--($(l)!1/4!(g)$);
      \draw [name path=E,green,ultra thick]($(o)!1/2!(c)$)--($(o)!1/2!(d)$);
      \draw [name path=B_1,blue](i)--($(i)!1/4!(m)$);
      \draw [name path=B_1',blue](m)--($(m)!1/3!(i)$);
      \draw [name path=B_2,blue](i')--($(i')!1/4!(m')$);
      \draw [name path=B_2',blue](m')--($(m')!1/3!(i')$);
      \draw [name path=C_1,black](j)--($(j)!1/4!(n)$);
      \draw [name path=C_1',black](o)--(n);
      \draw [name path=C_2,black](j')--($(j')!1/4!(n')$);
      \draw [name path=C_2',black]($(e)!1/2-1/60!(b)$)--(n');
      \node[above,black]at($(b)!1/2!(c)$){$C_3$};
      \node[below,black]at($(e)!1/2!(f)$){$C_0$};
      \node[right,blue]at($(a)!1/2!(b)$){$B_0$};
      \node[left,blue]at($(d)!1/2!(e)$){$B_3$};
      \node[left,red]at($(c)!1/2!(d)$){$A_0$};
      \node[right,red]at($(a)!1/2!(f)$){$A_3$};
\node at (3/2,0) {$+$};
      \coordinate (a) at (1+3,0+3^0.5/2);
      \coordinate (b) at (-1+3,0+3^0.5/2);
      \coordinate (c) at (0+3,-3^0.5+3^0.5/2);
      \draw[name path=B_1,blue](a)--(b);
      \draw[name path=B_2,blue](c)--($(a)!1/2!(b)$);
      \node[above,blue]at($(a)!1/2!(b)$){$B_1$};
      \draw[name path=A_0,red](b)--(c);
      \node[left,red]at($(b)!1/2!(c)$){$A_2$};
      \draw[name path=A_1,red](a)--($(b)!1/2!(c)$);
      \draw[name path=E,ultra thick,green]($(a)!1/2!(b)$)--($(b)!1/2!(c)$);
      \draw[name path=A_C_1,black](a)--($(b)!1/3!(c)$);   
      \path[name intersections={of=A_1 and B_2}]; 
      \draw[name path=C_1,black](b)--(intersection-1);
      \node[right] at (a){$P_1$};
      \path[fill=green](a)circle[radius=0.03];
      \node[left] at (b){$P_3$};
      \path[fill=green](b)circle[radius=0.03];
      \node[right] at (intersection-1){$P_2$};
      \path[fill=green] (intersection-1)circle[radius=0.03];
      \node[right,red] at ($(a)!1/2!(intersection-1)$){$A_1$};
      \node[above,black] at ($(b)!1/2!(intersection-1)$){$C_1$};
\end{tikzpicture}
\par\end{center}

Now we blow up the total space $\mathcal{Y}_{1}=Bl_{P}\mathcal{Y}$
along the curves $\widetilde{L_{P_{1}}}$, $\widetilde{L_{P_{2}}}$$ $
and $\widetilde{L_{P_{3}}}$, which are the proper transformation
of $L_{P_{1}}$, $L_{P_{2}}$ and $L_{P_{3}}$. The central fiber
turns to be $Bl_{4}\mathbb{P}^{2}\cup Bl_{3}\mathbb{P}^{2}$. The
component $Bl_{3}\mathbb{P}^{2}$ is the blowup of $\mathbb{P}^{2}$
at three points $P_{1},P_{2},P_{3}$. When we run the minimal model
program, in the central fiber, the curves $A_{1},B_{1},C_{2}$ in
the component $Bl_{3}\mathbb{P}^{2}$ are contracted. In the general
fiber, the curves $A_{1},B_{1},C_{1}$ are contracted as well. Clearly
we also have that $Bl_{3}\mathbb{P}^{2}$ goes back to $\mathbb{P}^{2}$
in the central fiber. The general fiber of the lc model is $\Sigma^{c}$,
which we described at the beginning of this section, and the central
fiber is $\mathbb{P}^{1}\times\mathbb{P}^{1}\cup\mathbb{P}^{2}$. 

\begin{center}
\begin{tikzpicture}[xscale=1.5,yscale=1.5][font=\tiny]
  \coordinate (a) at (0,0);
  \coordinate (b) at (1,0);
  \coordinate (c) at (1,1);
  \coordinate (d) at (0,1);
  \draw [name path=B_3,blue](a)--(b);
  \draw [name path=B_0,blue](b)--(c);
  \draw [name path=C_3,black](c)--(d);
  \draw [name path=C_0,black](d)--(a);
  \draw [name path=A_0,red]($(a)!1/2!(b)$)--($(c)!1/2!(d)$);
  \draw [name path=A_3,red]($(b)!1/2!(c)$)--($(d)!1/2!(a)$);
  \draw [name path=E,green,ultra thick](b)--(d);
      \coordinate (a) at (0+2,0+1);
      \coordinate (b) at (1+2,0+1);
      \coordinate (c) at (0+2,-1+1);
      \draw[name path=E,green,ultra thick] (a)--(c);
      \draw[name path=B_2,blue](a)--(b);
      \draw[name path=C_1,black]($(a)!1/2!(c)$)--(b);
      \draw[name path=A_2,red](c)--($(a)!1/2!(b)$);
     \path[name path=E_1,green](a)--(c);
     \path[name path=E_2,green](b)--(c);
     \path[name path=E_3,green](a)--($(a)!1/2!(c)$);
     \node [left,red] at (c){$Q_1$};
     \node [left,black] at ($(a)!1/2!(c)$){$Q_2$};
     \node [left,blue] at (a){$Q_3$};
 \node at (1.4,0.4) {$+$};
 \node at ($(1.4,0.4)+(0,-1)$){$Case$ $1$};

\end{tikzpicture}
\par\end{center}

\textbf{Case 2.}\emph{ }When $P_{1}$ is on $B_{0}$ and $P_{2}$
is on $A_{3}$. We first blow up the total space along the curve $A_{3}$,
then blow up along the strict preimage of $B_{0}$ in $Bl_{3}\mathbb{P}^{2}$.
Finally we blow up the total space along the three curves $\widetilde{L_{P_{1}}}$,
$\widetilde{L_{P_{2}}}$$ $ and $\widetilde{L_{P_{3}}}$. The central
fiber of the minimal model is $Bl_{4}\mathbb{P}^{2}\cup Bl_{2}\mathbb{F}_{1}\cup Bl_{1}\mathbb{F}_{0}$.
We obtain the central fiber $\mathbb{P}^{2}\cup\mathbb{P}^{2}\cup\mathbb{P}^{2}$
of the canonical model by contracting the lines labeled in the fourth
figure below. The general fiber of the lc model is again $\Sigma^{c}$.

\begin{center}
\begin{tikzpicture}[xscale=1.5,yscale=1.5][font=\tiny]
      \coordinate (d) at (1,3^0.5/2);
      \coordinate (e) at (1/2,3^0.5);
      \coordinate (f) at (-1/2,3^0.5);
      \coordinate (a) at (-1,3^0.5/2);
      \coordinate (b) at (-1/2,0);
      \coordinate (c) at (1/2,0);
      \draw[blue](a)--(b);
      \draw[black](b)--(c);
      \draw[red](c)--(d);
      \draw[blue](d)--(e);
      \draw[black](e)--(f);
      \draw[red](f)--(a);
      \draw[name path=C_1,black]($(a)!1/2!(f)$)--($(c)!1/2!(d)$);
      \draw[name path=B_2,blue]($(f)!1/2!(e)$)--($(b)!1/2!(c)$);
      \draw[name path=A_1,red]($(a)!1/2!(b)$)--($(e)!1/2!(d)$);
      \draw[red]($(d)!1/30!(e)$)--($(c)!1/60!(f)$);
      \draw[red]($(c)!1/60!(f)$)--($(b)!1/30!(a)$);
      \draw[blue]($(e)!2/30!(f)$)--($(d)!2/60!(a)$);
      \draw[blue]($(d)!2/60!(a)$)--($(c)!2/30!(b)$);
      \draw[black]($(d)!1/30!(c)$)--($(e)!1/60!(b)$);
      \draw[black]($(e)!1/60!(b)$)--($(f)!1/30!(a)$);
      \path [name intersections={of=A_1 and B_2}]; 
      \coordinate (P_2) at (intersection-1);
      \path[fill=green](P_2)circle[radius=0.03];
      \coordinate (P_1) at ($(e)!1/2!(d)$);
      \path[fill=green](P_1)circle[radius=0.03]; 
      \coordinate (P_3) at ($(c)!1/2!(d)$);
      \path[fill=green](P_3)circle[radius=0.03];  
  \node at ($(d)+(0.3,0)$){$\rightsquigarrow$};
      \coordinate (d) at (1+2.7,3^0.5/2);
      \coordinate (e) at (1/2+2.7,3^0.5);
      \coordinate (f) at (-1/2+2.7,3^0.5);
      \coordinate (a) at (-1+2.7,3^0.5/2);
      \coordinate (b) at (-1/2+2.7,0);
      \coordinate (c) at (1/2+2.7,0);
      \draw[blue](a)--(b);
      \draw[black](b)--(c);
      \draw[ultra thick,blue](c)--(d);
      \draw[blue](d)--(e);
      \draw[black](e)--(f);
      \draw[red](f)--(a);
      \draw[name path=C_1,black]($(a)!1/2!(f)$)--($(c)!1/2!(d)$);
      \draw[name path=B_2,blue]($(f)!1/2!(e)$)--($(b)!1/2!(c)$);
      \draw[name path=A_1,red]($(a)!1/2!(b)$)--($(e)!1/2!(d)$);
      \draw[red]($(c)!1/30!(d)$)--($(b)!1/30!(a)$);
      \draw[blue]($(e)!2/30!(f)$)--($(d)!2/30!(c)$);
      \draw[black]($(d)!1/30!(c)$)--($(e)!1/60!(b)$);
      \draw[black]($(e)!1/60!(b)$)--($(f)!1/30!(a)$);
     \path [name intersections={of=A_1 and B_2}]; 
     \coordinate (P_2) at (intersection-1);
     \path[fill=green](P_2)circle[radius=0.03];
      \coordinate (w) at ($(a)!1+2/3!(c)$);
      \coordinate (x) at ($(f)!1+2/3!(d)$);
      \coordinate (u) at ($(c)!1/2!(d)$);
      \coordinate (v) at ($(w)!1/2!(x)$);
      \draw[black](c)--(w);
      \draw[black]($(u)$)--($(v)$);
      \draw[black]($(d)!1/30!(c)$)--($(x)!1/30!(w)$);
      \draw[red](w)--(x);
      \draw[red]($(c)!1/30!(d)$)--($(d)!2/3!(x)$);
      \draw[blue]($(d)!2/30!(c)$)--($(c)!2/3!(w)$);
      \draw[blue](d)--(x);
     \coordinate (P_3) at ($(u)!1/3!(v)$);
     \path[fill=green](P_3)circle[radius=0.03];  
     \coordinate (P_1) at ($(e)!1/2!(d)$);
     \path[fill=green](P_1)circle[radius=0.03];   
 \node at ($(d)+(0.7,0)$){$\rightsquigarrow$};
      \coordinate (d) at (1+6,3^0.5/2);
      \coordinate (e) at (1/2+6,3^0.5);
      \coordinate (f) at (-1/2+6,3^0.5);
      \coordinate (a) at (-1+6,3^0.5/2);
      \coordinate (b) at (-1/2+6,0);
      \coordinate (c) at (1/2+6,0);
      \draw[blue](a)--(b);
      \draw[black](b)--(c);
      \draw[ultra thick,blue](c)--(d);
      \draw[black](e)--(f);
      \draw[red](f)--(a);
      \draw[name path=C_1,black]($(a)!15/28!(f)$)--($(c)!15/28!(d)$);
      \draw[name path=B_2,blue]($(f)!15/28!(e)$)--($(b)!15/28!(c)$);
      \draw[name path=A_1,red]($(a)!1/2!(b)$)--($(e)!1/2!(d)$);
      \draw[red]($(c)!1/30!(d)$)--($(b)!1/30!(a)$);
      \draw[black]($(e)!1/60!(b)$)--($(f)!1/30!(a)$);
      \draw[ultra thick,black](d)--(e);
     \path [name intersections={of=A_1 and B_2}]; 
     \coordinate (P_2) at (intersection-1);
     \path[fill=green](P_2)circle[radius=0.03];
      \coordinate (w) at ($(a)!1+2/3!(c)$);
      \coordinate (x) at ($(f)!1+2/3!(d)$);
      \coordinate (u) at ($(c)!1/2!(d)$);
      \coordinate (v) at ($(w)!1/2!(x)$);
      \coordinate (y) at ($(a)!1+1/2!(d)$);
      \draw[name path=C_0,black](c)--(w);
      \draw[name path=B_1,blue]($(c)!3/8!(w)$)--($(d)!1/2!(y)$);
      \draw[name path=A_3,red](w)--(x);
      \draw[name path=A_2,red]($(c)!1/30!(d)$)--($(x)!1/2!(y)$);
\draw[name path=C_2,black]($(y)!1/30!(d)$)--($(x)!1/15!(w)$);
 
      \draw[red,ultra thick](d)--(y);
      \draw[name path=B_0,blue](x)--(y);
      \draw[name path=C_1,black]($(c)!15/28!(d)$)--($(w)!15/28!(x)$);
     \path [name intersections={of=C_1 and B_1}]; 
     \coordinate (P_3) at (intersection-1);
     \path[fill=green](P_3)circle[radius=0.03];
      \coordinate (z) at ($(f)!2!(e)$);
      \draw[name path=B_0,blue](y)--(z);
      \draw[name path=C_3,black](e)--(z);
      \draw[name path=C_2,black]($(y)!1/30!(d)$)--($(e)!1/30!(d)$);
      \draw[name path=B_1,blue]($(d)!1/2!(y)$)--($(z)!1/2!(e)$);
      \draw[name path=A_1,red]($(d)!1/2!(e)$)--($(y)!1/2!(z)$);
     \path [name intersections={of=A_1 and B_1}]; 
     \coordinate (P_1) at (intersection-1);
     \path[fill=green](P_1)circle[radius=0.03];

\end{tikzpicture}     
\par\end{center}

\begin{center}
\begin{tikzpicture}[xscale=1.5,yscale=1.5][font=\tiny]
      \coordinate (d) at (1,3^0.5/2);
      \coordinate (e) at (1/2,3^0.5);
      \coordinate (f) at (-1/2,3^0.5);
      \coordinate (a) at (-1,3^0.5/2);
      \coordinate (b) at (-1/2,0);
      \coordinate (c) at (1/2,0);
      \coordinate (g) at ($(a)!1/2!(b)$);
      \coordinate (h) at ($(b)!15/28!(c)$);
      \coordinate (i) at ($(c)!15/28!(d)$);
      \coordinate (j) at ($(e)!1/2!(d)$);
      \coordinate (k) at ($(f)!15/28!(e)$);
      \coordinate (l) at ($(a)!15/28!(f)$);
      \draw[blue](a)--(b);
      \draw[black](b)--(c);
      \draw[ultra thick,blue](c)--(d);
      \draw[black](e)--(f);
      \draw[red](f)--(a);
      \draw[name path=C_1,black](l)--(i);
      \draw[name path=B_2,blue](k)--($(k)!1/3!(h)$);
      \draw[name path=B_2',blue](h)--($(h)!1/3!(k)$);
      \draw[name path=A_1,red](g)--($(g)!1/3!(j)$);
      \draw[name path=A_1',red](j)--($(j)!1/3!(g)$);
      \draw[name path=E_2,green]($(k)!1/3!(h)$)--($(g)!1/3!(j)$);
      \draw[name path=A_2,red]($(c)!1/30!(d)$)--($(b)!1/30!(a)$);
      \draw[name path=C_2,black]($(e)!1/60!(b)$)--($(f)!1/30!(a)$);
      \draw[ultra thick,black](d)--(e);
      \node[above,black] at ($(e)!1/2!(f)$){$C_3$};
      \node[below,black] at ($(b)!1/2!(c)$){$C_0$};
      \node[below,black] at ($(i)!1/3!(l)$){$C_1$};
      \node[above,red] at ($(j)!1/3!(g)$){$A_1$};
      \coordinate (w) at ($(a)!1+2/3!(c)$);
      \coordinate (x) at ($(f)!1+2/3!(d)$);
      \coordinate (u) at ($(w)!15/28!(x)$);
      \coordinate (v) at ($(c)!3/8!(w)$);
      \coordinate (y) at ($(a)!1+1/2!(d)$);
      \coordinate (m) at ($(d)!1/2!(y)$);      
      \draw[name path=C_0,black](c)--(w);
      \draw[name path=B_1,blue](v)--($(v)!1/4!(m)$);
      \draw[name path=B_1',blue](m)--($(m)!1/3!(v)$);
      \draw[name path=A_3,red](w)--(x);
      \draw[name path=A_2,red]($(c)!1/30!(d)$)--($(x)!1/2!(y)$);
      \draw[red,ultra thick](d)--(y);
      \draw[name path=B_0,blue](x)--(y);
      \draw[name path=C_2,black]($(y)!1/30!(d)$)--($(x)!1/15!(w)$);
      \draw[name path=C_1,black](i)--($(i)!1/4!(u)$);
      \draw[name path=C_1',black](u)--($(u)!1/3!(i)$);
      \draw[name path=E_3,green]($(m)!1/3!(v)$)--($(u)!1/3!(i)$);
      \node[left,blue]  at ($(m)!1/3!(v)$){$B_1$};
      \node[below,black] at ($(u)!1/3!(i)$){$C_1$};
      \coordinate (z) at ($(f)!2!(e)$);
      \coordinate (n) at ($(z)!1/2!(e)$);
      \coordinate (o) at ($(y)!1/2!(z)$);
      \draw[name path=B_0,blue](y)--(z);
      \draw[name path=C_3,black](e)--(z);
      \draw[name path=C_2,black]($(y)!1/30!(d)$)--($(e)!1/30!(d)$);
      \draw[name path=B_1,blue](m)--($(m)!1/3!(n)$);
      \draw[name path=B_1',blue](n)--($(n)!1/4!(m)$);
      \draw[name path=A_1,red](j)--($(j)!1/4!(o)$);
      \draw[name path=A_1',red](o)--($(o)!1/3!(j)$);
      \draw[name path=E_1,green]($(j)!1/4!(o)$)--($(n)!1/4!(m)$);
      \node[right,blue] at ($(n)!1/4!(m)$){$B_1$};
      \node[below,red] at ($(j)!1/4!(o)$){$A_1$};  
      \node at ($(y)+(0.5,0)$){$\rightsquigarrow$}; 
      \coordinate (a) at (-1+3,3^0.5/2);
      \coordinate (d) at (1+3,3^0.5/2);
      \coordinate (e) at (1/2+3,3^0.5);
      \coordinate (c) at (1/2+3,0);
      \coordinate (g) at ($(a)!1/2!(b)$);
      \coordinate (i) at ($(c)!15/28!(d)$);
      \coordinate (j) at ($(d)!1/2!(e)$);
      \draw[ultra thick,blue](c)--(d);
      \draw[name path=A_0,red](e)--(i);
      \draw[name path=B_2,blue](c)--(e);
      \draw[name path=B_3,blue](c)--(j);
      \draw[name path=E_2,green](i)--(j);
      \draw[ultra thick,black](d)--(e);
      \coordinate (y) at ($(a)!1+1/2!(d)$);
      \coordinate (m) at ($(d)!1/2!(y)$);      
      \draw[name path=C_0,black](c)--(m);
      \draw[name path=A_3,red](i)--(y);
      \draw[name path=A_2,red](c)--(y);
      \draw[red,ultra thick](d)--(y);
      \draw[name path=E_3,green](i)--(m);
      \draw[name path=B_0,blue](y)--(j);
      \draw[name path=C_3,black](e)--(m);
      \draw[name path=C_2,black](y)--(e);
      \draw[name path=E_1,green](m)--(j);
\node at ($(c)+(1/2,-1/2)$){$Case$ $2$};

\end{tikzpicture}     
\par\end{center}

There are only 2 types of degenerations with reducible lc models for
$K^{2}=3$. Both of them could be obtained from cases with $K^{2}=4$
.

We summarize the above computations in the following statement:\textbf{}\\
\textbf{}\\
\textbf{Theorem \ref{maintheory}.} \emph{The main component of the
compactified coarse moduli space $\overline{M}_{Bur}^{d}$ of stable
Burniat surfaces, or equivalently, of stable pairs $\left(Y,\frac{1}{2}D\right)$,
is of dimension $d-2$, irreducible for $d\neq4$, and with two components
for $d=4$. The types of degenerations, up to symmetry, are listed
as below.}

\emph{(i) There are 6 types of degenerate configurations of stable
pairs with reduced log canonical models in the moduli space of stable
pairs $\left(Y,\frac{1}{2}D\right)$ for $K^{2}=5$ case up to the
symmetry group $\mathbb{Z}_{6}$ described in Section \ref{degree5}.}

\emph{(ii) There are 5 types of degenerations with reducible lc models
in the moduli space of stable pairs $\left(Y,\frac{1}{2}D\right)$
for $K^{2}=4$ nodal case and 3 types of degeneration for $K^{2}=4$
non-nodal case up to the symmetry group $\mathbb{Z}_{2}$ described
in Section \ref{degree4}.}

\emph{(iii) There are only 2 types of degenerations with reducible
lc models in the moduli space of stable pairs $\left(Y,\frac{1}{2}D\right)$
for $K^{2}=3$ described in Section \ref{degree3}. }

There is only one surface with $K^{2}=2$, thus the moduli space of
Burniat surfaces with $K^{2}=2$ is just a single point.\emph{ }

\section{Matroid tilings of polytopes $\triangle_{Bur}^{d},d\leq5$ \label{matroid}}

According to the general theory of \cite{Ale08}, the unweighted stable
hyperplane arrangements are described by matroid tilings of the hypersimplex
$\triangle(r,n)$. Their weighted counterparts are described by partial
tilings of the hypersimplex $\triangle(r,n)$ in $\mathbb{R}^{n}$.
In this section, we will discuss the matroid tiling of the certain
polytopes $\triangle_{Bur}^{d},d\leq5$ corresponding to Burniat surfaces
with $K^{2}=d$.

In \cite{AP09}, Alexeev-Pardini defined the polytope $\triangle_{Bur}^{6}$
corresponding to Burniat surfaces with $K^{2}=6$, which is a subpolytope
of a hypersimplex $\triangle(3,9)$. They computed all stable Burniat
surfaces with $K^{2}=6$ by computing matroid tilings of a certain
polytope $\triangle_{Bur}^{6}$. We define the corresponding polytopes
$\triangle_{Bur}^{d},d\leq5$ similarly to $\triangle_{Bur}^{6}$.
We restrict the matroid tilings of the polytope $\triangle_{Bur}^{d}$
to $\triangle_{Bur}^{d-1}$ for $d\leq6$ to find all possible stable
surfaces in the main component of the compactified moduli space of
Burniat surfaces with $K^{2}=5$. 

Let's recall some definitions and results in \cite{AP09,Ale08}. A
\emph{hypersimplex} $\triangle(r,n)$ is defined to be a convex hull

\begin{eqnarray*}
\triangle(r,n) & = & \mbox{Conv}(e_{I}|I\in\overline{n},|I|=r)\\
 & = & \{(x_{1},...,x_{n})\in\mathbb{R}^{n}|0\leq x_{i}\leq1,\Sigma x_{i}=r\}
\end{eqnarray*}

A \emph{matroid polytope $ $$BP_{V}\subset\triangle(r,n)$ }is the
polytope corresponding to the toric variety $\overline{T.V}$ for
some geometric point $[V\subset\mathbb{A}^{n}]\in G(r,n)(k)$. One
can also describe the matroid polytopes in terms of hyperplane arrangements.
Let $\mathbb{P}V\simeq\mathbb{P}^{r}$ and assume that it is not contained
in the $n$ coordinate hyperplane $H_{i}$ (i.e. all $z_{i}\neq0$
on $\mathbb{P}V$); let $B_{1},...,B_{n}\subset\mathbb{P}V$ be $H_{i}\cap\mathbb{P}V$.
Then for the hyperplane arrangement $\left(\mathbb{P}V,\sum B_{i}\right)$,
the matroid polytope $BP_{V}$ is the convex hull of the points $v_{I}\in\mathbb{Z}^{n}$
for all $I\subset\overline{n}$ such that $\cap_{i\in I}B_{i}=\varnothing$,
or in terms of inequalities as 

\[
BP_{V}=\left\{ (x_{1},...,x_{n})\in\triangle(r,n)|\underset{i\in I}{\sum}x_{i}\leq\mbox{codim}\cap_{i\in I}B_{i},\mbox{ }\forall I\subset\overline{n}\right\} .
\]

For a hyperplane arrangement in general position, one has $BP_{V}=\triangle(r,n)$. 

Let $b=(b_{1},...,b_{n})$ be a weight, a $b$-\emph{cut hypersimplex}
is 

\begin{eqnarray*}
\triangle_{b}(r,n) & = & \left\{ (x_{1},...,x_{n})|0\leq x_{i}\leq b_{i},\Sigma x_{i}=r\right\} \\
 & = & \left\{ \alpha\in\triangle(r,n)|\alpha\leq b\right\} 
\end{eqnarray*}

We have the theorem in \cite{Ale08}
\begin{thm*}
\cite[2.12]{Ale08} The matroid polytope $BP_{V}$ is the set of points
$(x_{i})\in\mathbb{R}^{n}$ such that the pair $(\mathbb{P}V,\sum x_{i}B_{i})$
is lc and $K_{\mathbb{P}V}+\sum x_{i}B_{i}=0$; the interior $\mbox{Int}BP_{V}$
is the set of points such that $(\mathbb{P}V,\sum x_{i}B_{i})$ is
klt and $K_{\mathbb{P}V}+\sum x_{i}B_{i}=0$.
\end{thm*}
A tiling of the $b$-cut hypersimplex $\triangle_{b}$ is a partial
matroid tiling of $\triangle(r,n)$ such that $\cup BP_{M_{j}}\supset\triangle_{b}$
and such that all base polytopes $BP_{M_{j}}$ intersect the interior
of $\triangle_{b}$. 

Let $(\mathbb{P}V,\sum b_{i}B_{i})$ be a hyperplane arrangement.
For a point $p\in\mathbb{P}V$, we denote by $I(p)$ the set of $i\in\overline{n}$
such that $p\in B_{i}$. We define $\triangle_{b}^{p}$ to be the
face (possibly empty) of $\triangle_{b}$, where $x_{i}=b_{i}$ for
all $i\in I(p)$.
\begin{thm*}
\cite[6.6]{Ale08} Let $(\mathbb{P}V,\sum b_{i}B_{i})$ be a hyperplane
arrangement of general type. Suppose $BP_{M}\cap\triangle_{b}\neq\emptyset$.
Then $(\mathbb{P}V,\Sigma b_{i}B_{i})$ is lc at $p$ if and only
if $BP_{M}\cap\triangle_{b}^{p}\neq\emptyset$.
\end{thm*}
Now let us look at Burniat surfaces with $K^{2}=5$. 

A Burniat surface with $K^{2}=5$ is a $\mathbb{Z}_{2}^{2}$-cover
of $Bl_{4}\mathbb{P}^{2}$ for the data 
\begin{eqnarray*}
D & = & \sum\left(a_{i}A_{i}+b_{i}B_{i}+c_{i}C_{i}+eE\right)
\end{eqnarray*}

where $D_{a},D_{b},D_{c}$ are branched divisors of the Galois cover
and $E$ is not.

Denote by $\triangle_{Bur}^{d}$ the polytope corresponding to Burniat
surfaces with $K^{2}=d$. This is a subpolytope of the hypersimplex
$\triangle(3,9)$ with weight $b=(\frac{1}{2},...,\frac{1}{2})$.
In \cite{AP09}, the polytope $\triangle_{Bur}^{6}$ is defined to
be 

\begin{eqnarray*}
\triangle(3,9)\supset\triangle_{Bur}^{6} & = & \{(a_{0},a_{1},a_{2},b_{0},b_{2},b_{3},c_{0},c_{2},c_{3})\in\mathbb{R}^{9}\mbox{ satisfying }\\
 &  & 0\leq a_{i},b_{i},c_{i}\leq\dfrac{1}{2},e\leq0;\\
 &  & \sum_{i=0}^{2}(a_{i}+b_{i}+c_{i})=3;\\
 &  & 0\leq a_{3}=c_{0}+c_{1}+c_{2}+b_{0}-1\leq1/2;\\
 &  & 0\leq b_{3}=a_{0}+a_{1}+a_{2}+c_{0}-1\leq1/2;\\
 &  & 0\leq c_{3}=b_{0}+b_{1}+b_{2}+a_{0}-1\leq1/2;\}
\end{eqnarray*}

For the case $K^{2}=5$ , the divisor $D$ on $\Sigma$ satisfies
$K_{\Sigma}+D=0$ and we got an extra equation $e=a_{1}+b_{1}+c_{1}-1$
comparing to $\triangle_{Bur}^{6}$. Since the cover $\pi:X\rightarrow Y$
is unramified over $E$, the coefficient $e\leq\frac{r-1}{r}=0$,
where $r=1$ is the ramification index. Then we define 

\begin{eqnarray*}
\triangle(3,9)\supset\triangle_{Bur}^{5} & = & \{(a_{0},a_{1},a_{2},b_{0},b_{2},b_{3},c_{0},c_{2},c_{3})\in\mathbb{R}^{9}\mbox{ satisfying }\\
 &  & 0\leq a_{i},b_{i},c_{i}\leq\dfrac{1}{2},e\leq0;\\
 &  & \sum_{i=0}^{2}(a_{i}+b_{i}+c_{i})=3;\\
 &  & 0\leq a_{3}=c_{0}+c_{1}+c_{2}+b_{0}-1\leq1/2;\\
 &  & 0\leq b_{3}=a_{0}+a_{1}+a_{2}+c_{0}-1\leq1/2;\\
 &  & 0\leq c_{3}=b_{0}+b_{1}+b_{2}+a_{0}-1\leq1/2;\\
 &  & e=a_{1}+b_{1}+c_{1}-1\leq0\}
\end{eqnarray*}
We need to classify all matroid tilings of the polytope $\triangle_{Bur}^{5}$.
In \cite{AP09}, the authors listed all the nonempty intersection
of maximal-dimensional matroid polytopes $BP_{M}$ with the interior
of $\triangle_{Bur}^{6}$ in Table 1. It is easy to see that all the
base polytopes $BP_{M}$ listed in Table 1 \cite{AP09} still intersect
$\triangle_{Bur}^{5}$. 

Now we restrict the matroid tilings of the polytope $\triangle_{Bur}^{6}$
to $\triangle_{Bur}^{5}$ and list the tilings corresponding to degenerations
of stable Burniat surfaces of degree $5$, with reducible lc models
$(Y,\frac{1}{2}D)$. We will give the explanation below the table
where we use $a_{i}b_{j}c_{j}\leq1$ as the abbreviation of $a_{i}+b_{j}+c_{k}\leq1$.
\\
\\
\\

\begin{center}
Table 1
\par\end{center}

\begin{center}
\begin{tabular}{|c|c|c|}
\hline 
\# & Tilings for $K^{2}=5$ & From $K^{2}=6$ \tabularnewline
\hline 
\hline 
1 & $a_{0}a_{2}b_{2}\leq1,b_{2}b_{3}c_{2}\leq1$; $a_{1}c_{0}c_{2}\leq1,a_{1}a_{3}b_{1}\leq1$;  & Case 2\tabularnewline
 & $a_{2}c_{1}c_{3}\leq1,b_{0}b_{3}c_{1}\leq1$ & \tabularnewline
\hline 
2  & $a_{0}a_{2}b_{2}\leq1,b_{2}b_{3}c_{2}\leq1$;$a_{2}c_{1}c_{3}\leq1,b_{0}b_{3}c_{1}\leq1$;  & Case 3\tabularnewline
 & $a_{1}a_{3}b_{1}\leq1,a_{1}c_{0}c_{2}\leq1,b_{2}b_{3}c_{2}\leq1$; & \tabularnewline
 & $a_{0}a_{1}b_{1}\leq1,a_{1}c_{1}c_{3}\leq1,b_{1}b_{3}c_{2}\leq1$;
$ $ & \tabularnewline
\hline 
3 & $a_{0}a_{2}b_{2}\leq1,b_{2}b_{3}c_{2}\leq1$; & Case 4\tabularnewline
 & $a_{1}a_{3}b_{1}\leq1,a_{1}c_{0}c_{2}\leq1,b_{2}b_{3}c_{2}\leq1$; & \tabularnewline
 & $a_{0}a_{1}b_{1}\leq1,a_{1}c_{1}c_{3}\leq1,b_{1}b_{3}c_{2}\leq1$;  & \tabularnewline
 & $a_{1}c_{1}c_{3}\leq1,b_{0}b_{1}c_{1}\leq1,a_{1}a_{3}b_{1}\leq1$; & \tabularnewline
 & $a_{0}a_{2}b_{2}\leq1,a_{2}c_{1}c_{3}\leq1,b_{0}b_{1}c_{1}\leq1$; & \tabularnewline
\hline 
4 & $a_{0}a_{2}b_{0}\leq1$; $a_{1}a_{3}b_{1}\leq1$; & Case 6\tabularnewline
\hline 
5 & $a_{1}a_{3}b_{1}\leq1,a_{1}c_{0}c_{1}\leq1$; $a_{0}a_{2}b_{2}\leq1,a_{2}c_{2}c_{3}\leq1$;  & Case 7\tabularnewline
 & $a_{1}c_{2}c_{3}\leq1,b_{0}b_{1}c_{2}\leq1,a_{1}a_{3}b_{1}\leq1$; & \tabularnewline
 & $a_{1}a_{2}b_{2}\leq1,b_{2}b_{3}c_{1}\leq1,a_{0}c_{0}c_{1}\leq1$; & \tabularnewline
\hline 
6 & $a_{1}a_{2}b_{1}b_{2}c_{1}c_{2}\leq2$; $a_{0}b_{0}c_{0}\leq1$;  & Cases 9,10\tabularnewline
\hline 
\end{tabular}
\par\end{center}

In Table 1, tiling \#1 for $K^{2}=5$ is the union of 3 matroid polytopes
$BP_{M_{1}}\cup BP_{M_{2}}\cup BP_{M_{3}}$, where 
\begin{eqnarray*}
BP_{M_{1}} & = & \{a_{0}+a_{2}+b_{2}\leq1,b_{2}+b_{3}+c_{2}\leq1\}\cap\triangle(3,9)\\
BP_{M_{2}} & = & \{a_{1}+c_{0}+c_{2}\leq1,a_{1}+a_{3}+b_{1}\leq1\}\cap\triangle(3,9)\\
BP_{M_{3}} & = & \{a_{2}+c_{1}+c_{3}\leq1,b_{0}+b_{3}+c_{1}\leq1\}\cap\triangle(3,9)
\end{eqnarray*}
Tiling \#1 for $K^{2}=5$ is the same as the tiling \#2 in \cite{AP09}
for $K^{2}=6$.\\

We compare all the tilings of $\triangle_{Bur}^{6}$ in \cite{AP09}
with tilings of $\triangle_{Bur}^{5}$ listed above. The tiling \#2
of $\triangle_{Bur}^{6}$ is

\begin{eqnarray*}
BP_{M_{1}} & = & \left\{ a_{0}a_{1}a_{2}\leq1,c_{3}a_{1}a_{2}\leq1\right\} \\
BP_{M_{2}} & = & \left\{ b_{0}b_{1}b_{2}\leq1,a_{3}b_{1}b_{2}\leq1\right\} \\
BP_{M_{3}} & = & \left\{ c_{0}c_{1}c_{2}\leq1,b_{3}c_{1}c_{2}\leq1\right\} 
\end{eqnarray*}

$ $ The restriction of the tiling \#2 of $\triangle_{Bur}^{6}$ is
the tiling \#1 of $\triangle_{Bur}^{5}$. All the tilings of $\triangle_{Bur}^{5}$
above come from the restriction of the tilings of $\triangle_{Bur}^{6}$.
However, not all the restrictions of the tilings of $\triangle_{Bur}^{6}$
correspond to stable Burniat surfaces of degree 5. There are several
special cases to consider. \\

The tiling \#8 of $\triangle_{Bur}^{6}$ is the union of two base
polytopes 
\begin{eqnarray*}
BP_{M_{1}} & = & \left\{ a_{1}+a_{2}+b_{1}+b_{2}+c_{1}\leq2\right\} \cap\triangle(3,9)\\
BP_{M_{2}} & = & \left\{ a_{0}+b_{0}+c_{0}+c_{2}\leq1\right\} \cap\triangle(3,9)
\end{eqnarray*}

For the polytope $\triangle_{Bur}^{5}$, we have the inequalities
$a_{1}+b_{1}+c_{1}\leq1$ and $0\leq a_{2},b_{2}\leq\frac{1}{2}$.
These two inequalities imply $a_{1}+a_{2}+b_{1}+b_{2}+c_{1}\leq2$.
Hence $\triangle_{Bur}^{5}\subset BP_{M_{1}}$ and the corresponding
line arrangement for $BP_{M_{1}}$ is lc. This coincides with what
we got in Section \ref{degree5}. \\

The tiling \#10 of $\triangle_{Bur}^{6}$ composes of 3 matroid polytopes,
\begin{eqnarray*}
BP_{M_{1}} & = & \{a_{1}+a_{2}+b_{1}+b_{2}+c_{1}+c_{2}\leq2\}\cap\triangle(3,9)\\
BP_{M_{2}} & = & \{a_{0}+b_{0}+c_{0}\leq1,a_{1}+a_{2}+b_{1}+b_{2}+c_{1}\leq2\}\cap\triangle(3,9)\\
BP_{M_{3}} & = & \{a_{0}+b_{0}+c_{0}+c_{2}\leq1\}\cap\triangle(3,9)
\end{eqnarray*}
 The hypersimplex $\triangle(3,9)$ lies in the $\left\{ \underset{i=0}{\overset{2}{\Sigma}}a_{i}+b_{i}+c_{i}=3\right\} $,
and the complement of $BP_{M_{3}}$ in $\triangle(3,9)$ is $\left\{ a_{1}+a_{2}+b_{1}+b_{2}+c_{1}\geq2\right\} $.
For the polytope $\triangle_{Bur}^{5}$, the conditions $a_{1}+b_{1}+c_{1}\leq1$
and $a_{2}\leq\frac{1}{2},\mbox{ }b_{2}\leq\frac{1}{2}$ imply $ $$a_{1}+a_{2}+b_{1}+b_{2}+c_{1}\geq2$.
Therefore, $BP_{M_{3}}\cap\mbox{int}(\triangle_{Bur}^{5})=\emptyset$.
The restriction of the tiling \#10 of $\triangle_{Bur}^{6}$ and the
tiling \#9 of $\triangle_{Bur}^{6}$ are the same and is our tiling
\#6 of $\triangle_{Bur}^{5}$. \\

The tiling \#1 of $\triangle_{Bur}^{6}$ composes of 3 matroid polytopes,
\begin{eqnarray*}
BP_{M_{1}} & = & \left\{ a_{0}+a_{1}+a_{2}\leq1,a_{1}+a_{2}+c_{3}\leq1\right\} \cap\triangle(3,9)\\
BP_{M_{2}} & = & \left\{ b_{0}+b_{1}+b_{2}\leq1,a_{3}+b_{1}+b_{2}\leq1\right\} \cap\triangle(3,9)\\
BP_{M_{3}} & = & \left\{ c_{0}+c_{1}+c_{2}\leq1,b_{3}+c_{1}+c_{2}\leq1\right\} \cap\triangle(3,9)
\end{eqnarray*}

This three matroid polytopes correspond to three degenerations with
$K^{2}=6$. The restriction of the tiling is a tiling of $\triangle_{Bur}^{5}$
as well, but it does not correspond to any stable Burniat surface
of degree 5. To get the further degeneration with $A_{1},B_{1},C_{1}$
intersecting at a point, we can degenerate for instance $B_{1}$ to
$A_{0}+C_{3}$. The corresponding base polytope is 
\begin{eqnarray*}
BP_{M} & = & \{a_{0}+a_{1}+a_{2}+b_{1}\leq1,a_{1}+a_{2}+c_{3}\leq1,b_{1}+b_{0}\leq1\}\cap\triangle(3,9)
\end{eqnarray*}

Since we have $\Sigma_{i=0}^{3}a_{i}+b_{i}+c_{i}=3$, the inequality
$a_{0}+a_{1}+a_{2}+b_{1}\leq1$ is equivalent to $b_{0}+b_{2}+c_{0}+c_{1}+c_{2}\geq2$.
For $\triangle_{Bur}^{5}$, we have $a_{3}=c_{0}+c_{1}+c_{2}+b_{0}-1$,
so $b_{0}+b_{2}+c_{0}+c_{1}+c_{2}\geq2$ is the same as $b_{2}+a_{3}\geq1$.
But $0\leq b_{2},a_{3}\leq\frac{1}{2}$ for $\triangle_{Bur}^{5}$,
hence $BP_{M}\cap\mbox{int}(\triangle_{Bur}^{5})=\emptyset$ and the
base polytope $BP_{M}$ does not correspond to a degeneration with
$K^{2}=5$. This means the further degeneration is a lc degeneration
for $K^{2}=5$ with reducible lc model. We check all the possibility
of the degenerations coming from Case 1 in \cite{AP09} and find out
that the restriction of the tiling \#1 does not correspond to any
degenerations with $K^{2}=5$. \\

The tiling \#5 of $\triangle_{Bur}^{6}$ consists of 6 matroid polytopes.
WOLG, we pick one matroid polytope in the tiling

$BP_{M}=\{a_{0}+a_{1}+b_{2}\leq1,a_{1}+a_{2}+c_{3}\leq1,b_{2}+b_{3}+c_{1}\leq1\}\cap\triangle(3,9)$\\

If we force the condition that $A_{1},B_{1},C_{1}$ intersect at one
point, to the line arrangement corresponding to the polytope with
$K^{2}=6$, then the resulting degeneration has an finite automorphism
group. Therefore the tiling \#5 of $\triangle_{Bur}^{5}$ is still
a tiling of $\triangle_{Bur}^{5}$, but it does not correspond to
any degenerations with $K^{2}=5$.\\

Table 1 tells us that the 6 types of degenerations listed in Section
\ref{degree5} are all the degenerations up to symmetry in the main
component of the compactified moduli space of Burniat surfaces with
$K^{2}=5$.\\

We perform the same process for the cases $K^{2}=4$ and $K^{2}=3$.
 The following tables are for $K^{2}=4$ and $K^{2}=3$ cases. There
is no need to look at tilings for $K^{2}=2$, as the moduli space
for $K^{2}=2$ is just a point.

\begin{center}
Table 2
\par\end{center}

\begin{center}
\begin{tabular}{|c|c|c|}
\hline 
\# & Tilings for $K^{2}=4$ & From $K^{2}=5$ \tabularnewline
\hline 
\hline 
1 & $a_{1}c_{0}c_{2}\leq1,a_{1}a_{3}b_{1}\leq1$; $a_{1}a_{2}b_{2}\leq1,b_{2}b_{3}c_{2}\leq1$; & Case 1\tabularnewline
 & $b_{0}b_{1}c_{1}\leq1,a_{2}c_{1}c_{3}\leq1$ & \tabularnewline
\hline 
2  & $b_{0}b_{1}c_{1}\leq1,a_{1}a_{3}b_{1}\leq1,a_{1}c_{0}c_{2}\leq1$; & Case 2\tabularnewline
 & $a_{1}a_{3}b_{1}\leq1,a_{1}c_{0}c_{2}\leq1,b_{2}b_{3}c_{2}\leq1$; & \tabularnewline
\hline 
3 & $b_{0}b_{1}c_{1}\leq1$; $b_{2}b_{3}c_{2}\leq1$ & Case 4\tabularnewline
\hline 
4 & $a_{1}a_{2}b_{1}b_{2}c_{1}c_{2}\leq2$; $a_{0}b_{0}c_{0}\leq1$ & Case 6\tabularnewline
\hline 
5 & $a_{1}a_{2}b_{1}b_{2}c_{1}c_{2}\leq2$ & \tabularnewline
 & $a_{0}b_{0}c_{0}\leq1,a_{1}b_{1}b_{2}c_{1}c_{2}\leq2$ & Case 6\tabularnewline
\hline 
6 & $b_{0}b_{1}c_{1}\leq1$; $b_{2}b_{3}c_{2}\leq1$ & Case 3,4\tabularnewline
\hline 
7 & $b_{0}b_{1}c_{1}\leq1,a_{1}a_{3}b_{1}\leq1$; $a_{0}a_{2}b_{2}\leq1,b_{2}b_{3}c_{2}\leq1$; & Case 5\tabularnewline
 & $a_{0}a_{2}b_{2}\leq1,a_{2}c_{1}c_{3}\leq1,b_{0}b_{1}c_{1}\leq1$; & \tabularnewline
 & $a_{1}a_{3}b_{1}\leq1,a_{1}c_{0}c_{2}\leq1,b_{2}b_{3}c_{2}\leq1$; & \tabularnewline
\hline 
8 & $a_{1}a_{2}b_{1}b_{2}c_{1}c_{2}\leq2$; $a_{0}b_{0}c_{0}\leq1$; & Case 6\tabularnewline
\hline 
\end{tabular}
\par\end{center}

\pagebreak{}

\begin{center}
Table 3
\par\end{center}

\begin{center}
\begin{tabular}{|c|c|c|}
\hline 
\# & Tilings for $K^{2}=3$ & From $K^{2}=4$ \tabularnewline
\hline 
\hline 
1 & $a_{1}a_{2}b_{1}b_{2}c_{1}c_{2}\leq2$; $a_{0}b_{0}c_{0}\leq1$ & Case 4\tabularnewline
\hline 
2  & $b_{0}b_{2}c_{2}\leq1,a_{2}a_{3}b_{2}\leq1$; $a_{1}c_{1}c_{2}\leq1,b_{1}b_{3}c_{1}\leq1$; & \tabularnewline
 & $a_{0}a_{1}b_{1}\leq1,a_{1}c_{1}c_{3}\leq1$ & Case 1,3\tabularnewline
\hline 
\end{tabular}
\par\end{center}

\address{Department of Mathematics, University of Georgia, Athens, GA 30602,
USA}

\email{xhu@math.uga.edu}
\end{document}